\documentclass[11pt]{article}

\usepackage[small,bf]{caption}
\usepackage{amsmath, amsthm, amssymb}
\usepackage{mathtools}
\usepackage{algorithmic}
\usepackage{algorithm}
\usepackage[parfill]{parskip}
\usepackage{soul}

\usepackage{epsfig,fullpage,psfrag,url,verbatim}
\usepackage{graphics}
\usepackage{graphicx}
\usepackage{epstopdf}

\usepackage{amsfonts}
\usepackage{color}

\usepackage{array}
\usepackage{booktabs}
\usepackage[dvipsnames]{xcolor}
\usepackage{tikz}

\newcolumntype{M}[1]{>{\centering\arraybackslash}m{#1}}

\newtheorem{thm}{Theorem}[section]
\newtheorem{lem}{Lemma}[section]

\newtheorem{mydef}{Definition}[section]
\newtheorem{rem}{Remark}[section]
\newtheorem{cond}{Technical Condition}[section]

\newcommand{\bl}{{\bf{L}}}
\newcommand{\hh}{h(\cdot)}
\newcommand{\ff}{f(\cdot)}
\newcommand{\EE}{\mathbb{E}}

\newcommand{\Dh}{D_h}

\newcommand{\RR}{\mathbb{R}}

\begin{document}
\title{Accelerating Greedy Coordinate Descent Methods}

\author{Haihao Lu\thanks{MIT Department of Mathematics and MIT Operations Research Center, 77 Massachusetts Avenue, Cambridge, MA   02139
%(\href{mailto:rfreund@mit.edu}{rfreund@mit.edu}).  This author's research is supported by AFOSR Grant No. FA9550-11-1-0141 and the MIT-Chile-Pontificia Universidad
({mailto:  haihao@mit.edu}). This author's research is supported by AFOSR Grant No. FA9550-15-1-0276 and  the MIT-Belgium Universit\'{e} Catholique de Louvain Fund.}
\and
Robert M. Freund\thanks{MIT Sloan School of Management, 77 Massachusetts Avenue, Cambridge, MA   02139
%(\href{mailto:pgrigas@mit.edu}{pgrigas@mit.edu}).  This author's research has been partially supported through an NSF Graduate Research Fellowship and the
({mailto:  rfreund@mit.edu}).  This author's research is supported by AFOSR Grant No. FA9550-15-1-0276 and  the MIT-Belgium Universit\'{e} Catholique de Louvain Fund.}
\and
Vahab Mirrokni\thanks{Google Research, 111 8th Avenue New York, NY 10011 ({mailto:  mirrokni@google.com}).}
}
\date{Revised June 6, 2018 (original dated Feb. 9, 2018)}
%\date{} % Activate to display a given date or no date (if empty),
         % otherwise the current date is printed
\maketitle

\begin{abstract}
We study ways to accelerate greedy coordinate descent in theory and in practice, where ``accelerate'' refers either to  $O(1/k^2)$ convergence in theory, in practice, or both. We introduce and study two algorithms: Accelerated Semi-Greedy Coordinate Descent (ASCD) and Accelerated Greedy Coordinate Descent (AGCD).  While ASCD  takes greedy steps in the $x$-updates and randomized steps in the  $z$-updates, AGCD is a straightforward extension of standard greedy coordinate descent that only takes greedy steps.  On the theory side, our main results are for ASCD: we show that ASCD achieves $O(1/k^2)$ convergence, and it also achieves accelerated linear convergence for strongly convex functions.  On the empirical side, we observe that both AGCD and ASCD outperform Accelerated Randomized Coordinate Descent on a variety of instances.  In particular, we note that AGCD significantly outperforms the other accelerated coordinate descent methods in numerical tests, in spite of a lack of theoretical guarantees for this method.  To complement the empirical study of AGCD, we present a Lyapunov
energy function argument that points to an explanation for why a direct extension of the acceleration proof for AGCD
does not work; and we also introduce a technical condition under which AGCD is guaranteed to have accelerated convergence. Last of all, we confirm that this technical condition holds in our empirical study.
\end{abstract}

\section{Introduction: Related Work and Accelerated Coordinate Descent Framework}
Coordinate descent methods have received much-deserved attention recently due to their capability for solving large-scale optimization problems (with sparsity) that arise in machine learning applications and elsewhere. With inexpensive updates at each iteration, coordinate descent algorithms obtain faster running times than similar full gradient descent algorithms in order to reach the same near-optimality tolerance; indeed some of these algorithms now comprise the state-of-the-art in machine learning algorithms for loss minimization.

Most recent research on coordinate descent has focused on versions of randomized coordinate descent, which can essentially recover the same results (in expectation) as full gradient descent, including obtaining ``accelerated'' (i.e., $O(1/k^2)$) convergence guarantees. On the other hand, in some important machine learning applications, greedy coordinate methods demonstrate superior numerical performance while also delivering much sparser solutions.  For example, greedy coordinate descent is one of the fastest algorithms for the graphical LASSO implemented in DP-GLASSO \cite{mazumder2012graphical}.  And sequence minimization optimization (SMO) (a variant of greedy coordinate descent) is widely known as the best solver for kernel SVM \cite{joachims1998making}\cite{zeng2008fast} and is implemented in LIBSVM and SVMLight.

In general, for smooth convex optimization the standard first-order methods converge at a rate of $O(1/k)$ (including greedy coordinate descent). In 1983, Nesterov \cite{nesterov1983method} proposed an algorithm that achieved a rate of $O(1/k^2)$ -- which can be shown to be the optimal rate achievable by any first-order method \cite{NemirovskyYudin83}.  This method (and other similar methods) is now referred to as Accelerated Gradient Descent (AGD).  

However, there has not been much work on accelerating the standard Greedy Coordinate Descent (GCD) due to the inherent difficulty in demonstrating $O(1/k^2)$ computational guarantees (we discuss this difficulty further in Section \ref{sec:AGCDfails}). One work that might be close is \cite{song2017accelerated}, which updates the $z$-sequence using the full gradient and thus should not be considered as a coordinate descent method in the standard sense. There is a very related concurrent work \cite{locatello2018revisiting} and we will discuss the connections to our results in Section \ref{sec:AGCDfails}.

In this paper, we study ways to accelerate greedy coordinate descent in theory and in practice.  We introduce and study two algorithms: Accelerated Semi-Greedy Coordinate Descent (ASCD) and Accelerated Greedy Coordinate Descent (AGCD).  While ASCD  takes greedy steps in the $x$-updates and randomized steps in the $z$-updates, AGCD is a straightforward extension of GCD that only takes greedy steps. On the theory side, our main results are for ASCD: we show that ASCD achieves $O(1/k^2)$ convergence, and it also achieves accelerated linear convergence when the objective function is furthermore strongly convex. However, a direct extension of convergence proofs for ARCD does not work for ASCD due to the different coordinates we use to update $x$-sequence and $z$-sequence. Thus, we present a new proof technique -- which shows that a greedy coordinate step yields better objective function value than a full gradient step with a modified smoothness condition. %Indeed, this proof technique provides an alternative proof of the convergence rate of the ARCD algorithm as well and might be of independent interest.

On the empirical side, we first note that in most of our experiments ASCD outperforms Accelerated Randomized Coordinate Descent (ARCD) in terms of running time. On the other hand, we note that AGCD significantly outperforms the other accelerated coordinate descent methods in all instances, in spite of a lack of theoretical guarantees for this method. To complement the empirical study of AGCD, we present a Lyapunov energy function argument that points to an explanation for why a direct extension of the proof for AGCD does not work. This argument inspires us to introduce a technical condition under which AGCD is guaranteed to converge at an accelerated rate. Interestingly, we confirm that technical condition holds in a variety of instances in our empirical study, which in turn justifies our empirical observation that AGCD works very well in practice.

\subsection{Related Work}\label{sec:relatedwork}
{\bf Coordinate Descent.} Coordinate descent methods have a long history in optimization, and convergence of these methods has been extensively studied in the optimization community in the 1980s-90s, see \cite{bertsekas1989parallel}, \cite{luo1992convergence}, and \cite{luo1993error}. There are roughly three types of coordinate descent methods depending on how the coordinate is chosen: randomized coordinate descent (RCD), cyclic coordinate descent (CCD), and greedy coordinate descent (GCD).  RCD has received much attention since the seminal paper of Nesterov \cite{nesterov2012efficiency}. In RCD, the coordinate is chosen randomly from a certain fixed distribution. \cite{richtarik2014iteration} provides an excellent review of theoretical results for RCD. CCD chooses the coordinate in a cyclic order, see \cite{beck2013convergence} for basic convergence results.  More recent results show that CCD is inferior to RCD in the worst case \cite{sun2016worst}, while it is better than RCD in certain situations \cite{gurbuzbalaban2017cyclic}.  In GCD, we select the coordinate yielding the largest reduction in the objective function value.  GCD usually delivers better function values at each iteration in practice, though this comes at the expense of having to compute the full gradient in order to select the gradient coordinate with largest magnitude. The recent work \cite{nutini2015coordinate} shows that GCD has faster convergence than RCD in theory, and also provides several applications in machine learning where the full gradient can be computed cheaply.  A parallel GCD method is proposed in \cite{you2016asynchronous} and numerical results show its advantage in practice.

{\bf Accelerated Randomized Coordinate Descent.} Since Nesterov's paper on RCD \cite{nesterov2012efficiency} there has been significant focus on accelerated versions of RCD. In particular, \cite{nesterov2012efficiency} developed the first accelerated randomized coordinate gradient method for minimizing unconstrained smooth functions.  \cite{lu2015complexity} present a sharper convergence analysis of Nesterov's method using a randomized estimate sequence framework.  \cite{fercoq2015accelerated} proposed the APPROX (Accelerated, Parallel and PROXimal) coordinate descent method and obtained an accelerated sublinear convergence rate, and \cite{lee2013efficient} developed an efficient implementation of ARCD. 

%{\bf Regarding Connections with a Concurrent Paper \cite{locatello2018revisiting}}
%In a concurrent paper \cite{locatello2018revisiting}, the authors develop theories for matching pursuit algorithms, which can be viewed a generalized greedy coordinate descent where the allowed directions do not need to form an orthogonal basis. They also provided an accelerated version for matching pursuit algorithms, which turns out to be equivalent to ASCD discussed here in the special case where the chosen directions are orthogonal. Although the focuses of \cite{locatello2018revisiting} and our paper are different -- \cite{locatello2018revisiting} focuses on (accelerated) greedy direction update along a certain linear subspace and our method focuses on studying when and how we can accelerate greedy coordinate update -- both of the work share similar spirits when developing accelerating methods, that is, decoupling the coordinate update for $\{x^k\}$ sequence (with greedy rule) and for $\{z^k\}$ sequence (with randomized rule). This is also consistent with our argument on why we cannot accelerate greedy coordinate descent directly in Section \ref{sec:AGCDfails}.

\subsection{Accelerated Coordinate Descent Framework}

Our optimization problem of interest is:
\begin{equation}\label{poi1}
 P:  f^* := \ \  \mbox{minimum}_x \ \  f(x) \ , 
\end{equation}
where $f(\cdot) : \mathbb{R}^n \to \mathbb{R}$ is a differentiable convex function.\medskip

\begin{mydef} $f(\cdot)$ is coordinate-wise $L$-smooth for the vector of parameters $L:=(L_1,L_2,\ldots,L_n)$ if $\nabla f(\cdot)$ is coordinate-wise Lipschitz continuous for the corresponding coefficients of $L$, i.e., for all $x\in \RR^n$ and $h\in \RR$ it holds that:
\begin{equation}\label{eq:smoothness}
|\nabla_i f(x+h e_i)-\nabla_i f(x)| \le L_i |h|  \ , \ i=1, \ldots, n \ , 
\end{equation}

where $\nabla_i f(\cdot)$ denotes the $i^{\mathrm{th}}$ coordinate of $\nabla f(\cdot)$ and $e_i$ is $i^{\mathrm{th}}$ unit coordinate vector, for $i=1, \ldots, n$.  \end{mydef}

We presume throughout that $L_i > 0$ for $i=1, \ldots, n$.  Let $\bl$ denote the diagonal matrix whose diagonal coefficients correspond to the respective coefficients of $L$.   Let $\langle \cdot, \cdot \rangle$ denote the standard coordinate inner product in $\mathbb{R}^n$, namely $\langle x, y \rangle = \sum_{i=1}^n x_i y_i $, and let $\|\cdot\|_p$ denote the $\ell_p$ norm for $1 \le p \le \infty$.  \medskip
Let $\langle x, y \rangle_L := \sum_{i=1}^n L_i x_i y_i = \langle x, \bl y \rangle = \langle \bl x,  y \rangle$ denote the $L$-inner product.  Define the norm $\|x\|_L := \sqrt{\langle x, \bl x \rangle}$.  Letting $\bl^{-1}$ denote the inverse of $\bl$, we will also use the norm $\|\cdot\|_{L^{-1}}$ defined by $\|v\|_{L^{-1}} := \sqrt{\langle v, \bl^{-1} v \rangle} = \sqrt{\sum_{i=1}^n \bl_i^{-1}v_i^2}$.

%%\eqref{integration} can be written in an equivalent format using $\bl$: for any $x,y\in \mathbb{R}^n$, we have
%%\begin{equation}\label{eq:vector_integration}
%%f(x+y) \le f(x) + \langle y, \nabla f(x)\rangle + \tfrac{1}{2}\|y\|_{\bl}^2 .
%%\end{equation}
%\medskip

Algorithm \ref{al:nonframe} presents a generic framework for accelerated coordinate descent methods that is flexible enough to encompass deterministic as well as randomized methods.  One specific case is the standard Accelerated Randomized Coordinate Descent (ARCD).  In this paper we propose and study two other cases.  The first is Accelerated Greedy Coordinate Descent (AGCD), which is a straightforward extension of greedy coordinate descent to the acceleration framework and which, surprisingly, has not been previously studied (that we are aware of).  The second is a new algorithm which we call Accelerated Semi-Greedy Coordinate Descent (ASCD) that takes greedy steps in the $x$-updates and randomized steps in the $z$-update.  
\begin{algorithm}[tb]
\caption{Accelerated Coordinate Descent Framework without Strong Convexity}\label{al:nonframe}

\begin{algorithmic}
\STATE {\bf Initialize.}  Initialize with $x^0$, set $z^0 \leftarrow x^0$.  Assume $f(\cdot)$ is coordinate-wise $L$-smooth for known and given $L$.  Define the sequence $\{\theta_k\}$ as follows: $\theta_{0}=1$, and define $\theta_{k}$ recursively by the relationship $\tfrac{1-\theta_{k}}{\theta_{k}^{2}}=\tfrac{1}{\theta_{k-1}^{2}}$ for $k=1, 2, \ldots$. \\

$ \ $

At iteration $k$ :\\
\STATE  {\bf Perform Updates.}  

Define $y^{k}:=(1-\theta_{k})x^{k}+\theta_{k}z^{k}$ \\\medskip

Choose coordinate $j_{k}^{1}$ (by some rule) \\\medskip

Compute $x^{k+1}:=y^{k}-\tfrac{1}{L_{j_k^1}}\nabla_{j_{k}^{1}}f(y^{k}) e_{j_k^1}$ \\\medskip

Choose coordinate $j_{k}^{2}$ (by some rule) \\\medskip

Compute $z^{k+1}:=z^{k}-\tfrac{1}{nL_{j_{k}^2}\theta_{k}}\nabla_{j_{k}^{2}}f(y^{k}) e_{j_k^2}$ \ . \\\medskip

\end{algorithmic}
\end{algorithm}\medskip

In the framework of Algorithm \ref{al:nonframe} we choose a coordinate $j_{k}^{1}$ of the gradient $\nabla f(y^{k})$ to perform the update of the $x$-sequence, and we choose (a possibly different) coordinate $j_{k}^{2}$ of the gradient $\nabla f(y^{k})$ to perform the update of the $z$-sequence. Herein we will study three different rules for choosing the coordinates $j_k^1$ and $j_k^2$ which then define three different specific algorithms:

\begin{itemize}
\item ARCD (Accelerated Randomized Coordinate Descent): use the rule
\begin{equation}\label{eq:non_rule_ARCD}
j_{k}^{2}=j_{k}^{1}:\sim {\mathcal U}[1,\cdots,n]
\end{equation}
\item AGCD (Accelerated Greedy Coordinate Descent): use the rule
\begin{equation}\label{eq:non_rule_AGCD}
j_{k}^{2}=j_{k}^{1}:=\arg\max_{i}\tfrac{1}{\sqrt{L_i}}|\nabla_{i}f(y^{k})|
\end{equation}
\item ASCD (Accelerated Semi-Greedy Coordinate Descent): use the rule
\begin{equation}\label{eq:non_rule_ASCD}\begin{array}{l}
j_{k}^{1}:=\arg\max_{i}\tfrac{1}{\sqrt{L_i}}|\nabla_{i}f(y^{k})|  \\
j_{k}^{2}:\sim {\mathcal U}[1,\cdots,n] \ . 
\end{array}\end{equation}
\end{itemize}

In ARCD a random coordinate $j_k^1$ is chosen at each iteration $k$, and this coordinate is used to update both the $x$-sequence and the $z$-sequence.  ARCD is well studied, and is known to have the following convergence guarantee in expectation (see \cite{fercoq2015accelerated} for details):
\begin{equation}
E\left[f(x^{k})-f(x^{*})\right] \le \tfrac{2n^{2}}{(k+1)^{2}}\|x^{*}-x^{0}\|_L^{2} \ ,
\end{equation}
where the expectation is on the random variables used to define the first $k$ iterations.  

In AGCD we choose the coordinate $j_k^1$ in a ``greedy'' fashion, i.e., corresponding to the maximal (weighted) absolute value coordinate of the the gradient $\nabla f(y^k)$.  This greedy coordinate is used to update both the $x$-sequence and the $z$-sequence.  As far as we know AGCD has not appeared in the first-order method literature.  One reason for this is that while AGCD is the natural accelerated version of greedy coordinate descent, the standard proof methodologies for establishing acceleration guarantees (i.e., $O(1/k^2)$ convergence) fail for AGCD.  Despite this lack of worst-case guarantee, we show in Section \ref{sec:numerical} that AGCD is extremely effective in numerical experiments on synthetic linear regression problems as well as on practical logistic regression problems, and dominates other coordinate descent methods in terms of numerical performance.  Furthermore, we observe that AGCD attains $O(1/k^2)$ convergence (or better) on these problems in practice.  Thus AGCD is worthy of significant further study, both computationally as well as theoretically.  Indeed, in Section \ref{sec:AGCD} we will present a technical condition that implies $O(1/k^2)$ convergence when satisfied, and we will argue that this condition ought to be satisfied in many settings.  

ASCD, which we consider to be the new theoretical contribution of this paper, combines the salient features of AGCD and ARCD.  In ASCD we choose the greedy coordinate of the gradient to perform the $x$-update, while we choose a random coordinate to perform the $z$-update.  In this way we achieve the practical advantage of greedy $x$-updates, while still guaranteeing $O(1/k^2)$ convergence in expectation by virtue of choosing the random coordinate used in the $z$-update, see Theorem \ref{thm:nonstrong}.  And under strong convexity, ASCD achieves linear convergence as will be shown in Section \ref{sec:strongconvex}.

The paper is organized as follows.  In Section \ref{sec:ASCD} we present the $O(1/k^2)$ convergence guarantee (in expectation) for ASCD.  In Section \ref{sec:strongconvex} we present an extension of the accelerated coordinate descent framework to the case of strongly convex functions, and we present the associated linear convergence guarantee for ASCD under strong convexity. In Section \ref{sec:AGCD} we study AGCD; we present a Lyapunov energy function argument that points to why standard analysis of accelerated gradient descent methods fails in the analysis of AGCD. In Section \ref{sec:AGCDworks} we present a technical condition under which AGCD will achieve $O(1/k^2)$ convergence. In Section \ref{sec:numerical}, we present results of our numerical experiments using AGCD and ASCD on synthetic linear regression problems as well as practical logistic regression problems.

% we propose an Accelerated Semi-Greedy Coordinate Descent (ASCD) and the convergence result associated to it in both strongly convex case and non strongly-convex case. In section 3 {\color{red} (label)} we propose an Accelerated Greedy Coordinate Descent (AGCD), and introduce a technical assumption under which AGCD is guaranteed to converge with accelerated rate. In section 4 {\color{red} (label)}, we present two numerical studies: linear regression and logistic regression.)

\section{Accelerated Semi-Greedy Coordinate Descent (ASCD)}\label{sec:ASCD}

In this section we present our computational guarantee for the Accelerated Semi-Greedy Coordinate Descent (ASCD) method in the non-strongly convex case.  Algorithm \ref{al:nonframe} with rule \eqref{eq:non_rule_ASCD} presents the Accelerated Semi-Greedy Coordinate Descent method (ASCD) for the non-strongly convex case.  At each iteration $k$ the ASCD method choose the greedy coordinate $j_k^1$ to do the $x$-update, and chooses a randomized coordinate $j_k^2\sim  {\mathcal U}[1,\cdots,n]$ to do the $z$-update. Unlike ARCD where the same randomized coordinate is used in both the $x$-update and the $z$-update --  in ASCD $j_k^1$ is chosen in a deterministic greedy way, $j_k^1$ and $j_k^2$ are likely to be different.

At each iteration $k$ of ASCD the random variable $j_k^2$ is introduced, and therefore $x^{k}$ depends on the realization of the random variable
$$
\xi_{k} := \{j_0^2, \ldots, j_{k-1}^2\} \ .
$$ For convenience we also define $\xi_{0} := \emptyset$.

The following theorem presents our computational guarantee for ASCD for the non-strongly convex case:\medskip

\begin{thm}\label{thm:nonstrong}
Consider the Accelerated Semi-Greedy Coordinate Descent method (Algorithm \ref{al:nonframe} with rule \eqref{eq:non_rule_ASCD}). If $\ff$ is coordinate-wise $L$-smooth, it holds for all $k \ge 1$ that:

\begin{equation}\label{eq:conv_ASCD_non}
E_{\xi_{k}}\left[f(x^{k})-f(x^{*})\right]\le\tfrac{n^{2}\theta_{k-1}^{2}}{2}\|x^{*}-x^{0}\|_L^{2} \le \tfrac{2n^{2}}{(k+1)^{2}}\|x^{*}-x^{0}\|_L^{2} \ .
\end{equation}
 \qed
\end{thm}

In the interest of both clarity and a desire to convey some intuition on proofs of accelerated methods in general, we will present the proof of Theorem \ref{thm:nonstrong} after first establishing some intermediary results along with some explanatory comments.  We start with the ``Three-Point Property'' of Tseng \cite{tsengaccelerated}.  Given a differentiable convex function $h(\cdot)$, the Bregman distance for $h(\cdot)$ is $D_h(y,x) := h(y)-h(x) -\langle \nabla h(x), y-x \rangle $.  The Three-Point property can be stated as follows:\medskip\medskip

\begin{lem}\label{paul} {\bf (Three-Point Property (Tseng \cite{tsengaccelerated}))} Let $\phi(\cdot)$ be a convex function, and let $\Dh(\cdot, \cdot)$ be the Bregman distance for $\hh$. For a given vector $z$, let
\begin{equation*}
z^+ := \arg\min_{x\in \mathbb{R}^n} \left\{ \phi(x) + \Dh(x,z) \right\} \ .
\end{equation*}
Then
\begin{equation*}
\phi(x) + \Dh(x,z) \ge \phi(z^+) + \Dh(z^+, z) + \Dh(x,z^+)\ \   \mbox{ for\ all} \ \ x \in \mathbb{R}^n \ ,
\end{equation*}
with equality holding in the case when $\phi(\cdot)$ is a linear function and $h(\cdot)$ is a quadratic function. \qed
\end{lem}\medskip\medskip

Also, it follows from elementary integration and the coordinate-wise Lipschitz condition \eqref{eq:smoothness} that
\begin{equation}\label{integration} f(x+h e_i) \le f(x) + h \cdot \nabla_i f(x) + \tfrac{h^2 L_i }{2  } \ \ \mbox{for~all} \ x 
\in \mathbb{R}^n \ \mbox{and~} h \in \mathbb{R} \ . \end{equation}

At each iteration $k =0, 1, \ldots$ of ASCD, notice that $x^{k+1}$ is one step of greedy coordinate descent from $y^k$ in the norm $\| \cdot \|_L$. Now define $s^{k+1}:=y^{k}-\tfrac{1}{n}\bl^{-1}\nabla f(y^{k})$, which is a full steepest-descent step from $y^k$ in the norm $\|\cdot\|_{nL}$ . We first show that the greedy coordinate descent step yields a good objective function value as compared to the quadratic model that yields $s^{k+1}$.\medskip

\begin{lem}\label{robomart}
$$f(x^{k+1}) \le f(y^{k})+\langle\nabla f(y^{k}),s^{k+1}-y^{k}\rangle+\tfrac{n}{2}\|s^{k+1}-y^{k}\|_L^{2} \ .$$
\end{lem}
{\bf Proof: }
\begin{equation}\label{lem:non-first}
\begin{array}{lcl}
f(x^{k+1}) & \le &  f(y^{k})-\tfrac{1}{2L_{j_k^1}}\left(\nabla_{j_{k}^{1}}f(y^{k})\right)^{2}\\ \\
 & \le & f(y^{k})-\tfrac{1}{2n}\|\nabla f(y^{k})\|_{L^{-1}}^{2}\\ \\
 & = & f(y^{k})+\langle\nabla f(y^{k}),s^{k+1}-y^{k}\rangle+\tfrac{n}{2}\|s^{k+1}-y^{k}\|_L^{2} \ ,\\ \\  
 \end{array}
\end{equation}
where the first inequality of \eqref{lem:non-first} derives from the smoothness of $f(\cdot)$, and is a simple instance of \eqref{integration} using $x = y^k$, $i = j_k^1$, and $h = -\tfrac{1}{L_{j_k^1}}\nabla_{j_{k}^{1}}f(y^{k})$. The second inequality of \eqref{lem:non-first} follows from the definition of $j_k^1$ which yields:
$$n\left[\tfrac{1}{L_{j_k^1}}\left(\nabla_{j_{k}^{1}}f(y^{k})\right)^{2} \right] \ge \sum_{i=1}^n \tfrac{1}{L_{i}}\left(\nabla_{i}f(y^{k})\right)^{2}=  \left\| \nabla f(y^k) \right\|^2_{L^{-1}} \ . $$ The last equality of \eqref{lem:non-first} follows by using the definition of $s^{k+1}$ and rearranging terms.  \qed

\medskip

Utilizing the interpretation of $s^{k+1}$ as a gradient descent step from $y^k$ but with a larger smoothness descriptor ($nL$ as opposed to $L$), we can invoke the standard proof for accelerated gradient descent derived in \cite{tsengaccelerated} for example. We define $t^{k+1}:= z^{k} - \tfrac{1}{n\theta_k}\bl^{-1}\nabla f(y^k)$, or equivalently we can define $t^{k+1}$ by: 
\begin{equation}\label{eq:3pointtk}
t^{k+1}=\arg\min_z \  \langle \nabla f(y^k), z-z^k \rangle + \tfrac{n\theta_k}{2}\|z-z^k\|_{L}^2\ 
\end{equation}
(which corresponds to $z^{k+1}$ in \cite{tsengaccelerated} for standard accelerated gradient descent). Then we have:\medskip

\begin{lem}\label{lem:non-second}
\begin{equation}\label{eq:non-second}
f(x^{k+1}) \le (1-\theta_{k})f(x^{k})+\theta_{k}f(x^{*})+\tfrac{n\theta_{k}^2}{2}\|x^*-z^{k}\|_L^{2}-\tfrac{n\theta_k^2}{2}\|x^*-t^{k+1}\|_L^{2} \ . 
\end{equation}
\end{lem}
{\bf Proof.} Following from Lemma \ref{robomart}, we have
\begin{equation}\label{kreso}
\begin{array}{lcl}
f(x^{k+1}) & \le & f(y^{k})+\langle\nabla f(y^{k}),s^{k+1}-y^{k}\rangle+\tfrac{n}{2}\|s^{k+1}-y^{k}\|_L^{2} \\ \\
& = & f(y^{k})+\theta_k\left(\langle\nabla f(y^{k}),t^{k+1}-z^{k}\rangle+\tfrac{n\theta_k}{2}\|t^{k+1}-z^{k}\|_L^{2}\right)  \\ \\
& = & f(y^{k})+\theta_k \left( \langle\nabla f(y^{k}),x^*-z^{k}\rangle+\tfrac{n\theta_k}{2}\|x^*-z^{k}\|_L^{2} -  \tfrac{n\theta_k}{2}\|x^*-t^{k+1}\|_L^{2}\right) \\ \\
& = &  (1-\theta_{k})\left(f(y^{k})+\langle\nabla f(y^{k}),x^{k}-y^{k}\rangle\right)+\theta_{k}\left(f(y^k)+\langle\nabla f(y^{k}),x^*-y^{k}\rangle\right)\\ \\
 & &  \,\,+\tfrac{n\theta_{k}^2}{2}\|x^*-z^{k}\|_L^{2}-\tfrac{n\theta_k^2}{2}\|x^*-t^{k+1}\|_L^{2} \\ \\
 & \le & (1-\theta_{k})f(x^{k})+\theta_{k}f(x^{*})+\tfrac{n\theta_{k}^2}{2}\|x^*-z^{k}\|_L^{2}-\tfrac{n\theta_k^2}{2}\|x^*-t^{k+1}\|_L^{2}
\end{array}
\end{equation}
where the first equality of \eqref{kreso} utilizes $s^{k+1}-y^k=\theta_k(t^{k+1} - z^k)$.  The second equality of \eqref{kreso} follows as an application of the Three-Point-Property (Lemma \ref{paul}) together with \eqref{eq:3pointtk}, where we set $\phi(x) = \langle\nabla f(y^{k}),x-z^{k}\rangle$ and $h(x) = \tfrac{n\theta_k}{2}\|x\|_L^{2}$ (whereby $D_h(x,v) = \tfrac{n\theta_k}{2}\|x-v\|_L^{2}$).  The third equality of \eqref{kreso} is derived from $y^k = (1-\theta_k)x^k+\theta_k z^k$ and rearranging the terms. And the last inequality of \eqref{kreso} is an application of the gradient inequality at $y^k$ applied to $x^k$ and also to $x^*$.  \qed

Notice that $t^{k+1}$ is an all-coordinate update of $z^k$, and computing $t^{k+1}$ can be very expensive.  Instead we will use $z^{k+1}$ to replace $t^{k+1}$ in \eqref{eq:non-second} by using the equality in the next lemma.\medskip

\begin{lem}\label{lem:non-third}
\begin{equation}\label{eq:non-third}
\tfrac{n}{2}\|x^*-z^{k}\|_L^{2}-\tfrac{n}{2}\|x^*-t^{k+1}\|_L^{2} = \tfrac{n^{2}}{2}\|x^*-z^{k}\|_L^{2}-\tfrac{n^{2}}{2}E_{j_{k}^{2}}\left[\|x^*-z^{k+1}\|_L^{2}\right] \ .
\end{equation}
\end{lem}
{\bf Proof:}
 \begin{equation}\label{pats}\begin{array}{rcl}
\tfrac{n}{2}\|x^*-z^{k}\|_L^{2}-\tfrac{n}{2}\|x^*-t^{k+1}\|_L^{2} & = & \tfrac{n}{2}\left\langle t^{k+1}-z^{k},2x^*-2z^{k}\right\rangle_L -\tfrac{n}{2}\left\Vert t^{k+1}-z^{k}\right\Vert_L ^{2}\\ \\
 & = & \tfrac{n^{2}}{2}E_{j_{k}^{2}}\left[\left\langle z^{k+1}-z^{k},2x^*-2z^{k}\right\rangle_L -\left\Vert z^{k+1}-z^{k}\right\Vert_L ^{2}\right]\\ \\
 & = & \tfrac{n^{2}}{2}\|x^*-z^{k}\|_L^{2}-\tfrac{n^{2}}{2}E_{j_{k}^{2}}\left[\|x^*-z^{k+1}\|_L^{2}\right],
\end{array}\end{equation}
where the first and third equations above are straightforward arithmetic rearrangements,
and the second equation follows from the two easy-to-verify identities $t^{k+1}-z^{k}=nE_{j_{k}^{2}}\left[z^{k+1}-z^{k}\right]$ and $\left\Vert t^{k+1}-z^{k}\right\Vert_L ^{2}=nE_{j_{k}^{2}}\left[\left\Vert z^{k+1}-z^{k}\right\Vert_L ^{2}\right]$ .  \qed

We now have all of the ingredients needed to prove Theorem \ref{thm:nonstrong}.

{\bf Proof of Theorem \ref{thm:nonstrong}} Substituting \eqref{eq:non-third} into \eqref{eq:non-second}, we obtain:
\begin{equation}
f(x^{k+1})  \le  (1-\theta_{k})f(x^{k})+\theta_{k}f(x^{*})+\tfrac{n^{2}\theta_{k}^2}{2}\|x^{*}-z^{k}\|_L^{2}-\tfrac{n^{2}\theta_k^2}{2}E_{j_{k}^{2}}\left[\|x^{*}-z^{k+1}\|_L^{2}\right] \ .
\end{equation}

Rearranging and substituting $\tfrac{1-\theta_{k+1}}{\theta_{k+1}^{2}}=\tfrac{1}{\theta_{k}^{2}}$,
we arrive at:
\begin{equation*}
  \tfrac{1-\theta_{k+1}}{\theta_{k+1}^{2}}\left(f(x^{k+1})-f\left(x^{*}\right)\right)+\tfrac{n^{2}}{2}E_{j_{k}^{2}}\left\Vert x^{*}-z^{k+1}\right\Vert_L ^{2} \le  \left[\tfrac{1-\theta_{k}}{\theta_{k}^{2}}\left(f(x^{k})-f\left(x^{*}\right)\right)+\tfrac{n^{2}}{2}\left\Vert x^{*}-z^{k}\right\Vert_L ^{2}\right] \ .
\end{equation*}
Taking the expectation over the random variables $j_1^2, j_2^2, \ldots, j_k^2$, it follows that: 
\begin{equation*}
 E_{\xi_{k+1}}\left[ \tfrac{1-\theta_{k+1}}{\theta_{k+1}^{2}}\left(f(x^{k+1})-f\left(x^{*}\right)\right)+\tfrac{n^{2}}{2}\left\Vert x^{*}-z^{k+1}\right\Vert_L ^{2}\right] \le  E_{\xi_{k}} \left[\tfrac{1-\theta_{k}}{\theta_{k}^{2}}\left(f(x^{k})-f\left(x^{*}\right)\right)+\tfrac{n^{2}}{2}\left\Vert x^{*}-z^{k}\right\Vert_L ^{2}\right] \ .
\end{equation*}
Applying the above inequality in a telescoping manner for $k =1, 2, \ldots$, yields:
\begin{equation*}\begin{array}{rcl}
 E_{\xi_{k}}\left[ \tfrac{1-\theta_{k}}{\theta_{k}^{2}}\left(f(x^{k})-f\left(x^{*}\right)\right) \right] & \le &  E_{\xi_{k}}\left[ \tfrac{1-\theta_{k}}{\theta_{k}^{2}}\left(f(x^{k})-f\left(x^{*}\right)\right)+\tfrac{n^{2}}{2}\left\Vert x^{*}-z^{k}\right\Vert_L ^{2}\right] \\ & \vdots & \\ & \le & E_{\xi_{0}} \left[\tfrac{1-\theta_{0}}{\theta_{0}^{2}}\left(f(x^{0})-f\left(x^{*}\right)\right)+\tfrac{n^{2}}{2}\left\Vert x^{*}-z^{0}\right\Vert_L ^{2}\right] \\ \\
 & = & \tfrac{n^{2}}{2}\left\Vert x^{*}-x^{0}\right\Vert_L ^{2} \  .
\end{array}\end{equation*}
Note from an induction argument that $\theta_i \le \tfrac{2}{i+2}$ for all $=0, 1, \ldots$, whereby the above inequality rearranges to:
$$E_{\xi_{k}}\left[ \left(f(x^{k})-f\left(x^{*}\right)\right) \right] \le \tfrac{\theta_{k}^{2}}{1-\theta_{k}}\tfrac{n^{2}}{2}\left\Vert x^{*}-x^{0}\right\Vert_L ^{2}  = \tfrac{n^{2}\theta_{k-1}^{2}}{2}\left\Vert x^{*}-x^{0}\right\Vert_L ^{2} \le \tfrac{2n^{2}}{(k+1)^2}\left\Vert x^{*}-x^{0}\right\Vert_L ^{2}  \ . \qed $$

%%%%%%%%%%%%%%%%%%%%%%%%%%%%%%%%%%%%%%%%%%%%%%%%%%%%%

\section{Accelerated Coordinate Descent Framework under Strong Convexity}\label{sec:strongconvex}

We begin with the definition of strong convexity as developed in \cite{lu2015complexity}:\medskip

\begin{mydef} $f(\cdot)$ is $\mu$-strongly convex with respect to $\|\cdot\|_L$ if for all $x, y\in \RR^n$ it holds that:
$$
f(y) \ge f(x) + \langle \nabla f(x), y-x \rangle + \tfrac{\mu}{2}\|y-x\|_L^2 \ .
$$
\end{mydef}Note that $\mu$ can be viewed as an extension of the condition number of $\ff$ in the traditional sense since $\mu$ is defined relative to the coordinate smoothness coefficients through $\|\cdot\|_L$, see \cite{lu2015complexity}.  Algorithm \ref{al:strframe} presents the generic framework for accelerated coordinate descent methods in the case when $\ff$ is $\mu$-strongly convex for known $\mu$.

\begin{algorithm}
\caption{Accelerated Coordinate Descent Framework ($\mu$-strongly convex case)}\label{al:strframe}

$ \ $

\begin{algorithmic}
\STATE {\bf Initialize.}  Initialize with $z^0=x^0$.  Assume $f(\cdot)$ is coordinate-wise $L$-smooth and $\mu$-strongly convex for known and given $L$ and $\mu$, and define the parameters $a=\frac{\sqrt{\mu}}{n+\sqrt{\mu}}$ and $b=\frac{\mu a}{n^2}$ .   \\

$ \ $

At iteration $k$ :\\
\STATE  {\bf Perform Updates.}  

Define $y^{k}=(1-a)x^{k}+az^{k}$ \\\medskip

Choose $j_{k}^{1}$ (by some rule) \\\medskip

Compute $x^{k+1}=y^{k}-\frac{1}{L_{j_k^1}}\nabla_{j_{k}^{1}}f(y^{k}) e_{j_{k}^{1}}$ \\\medskip

Compute  $u^k=\frac{a^2}{a^2+b}z^k + \frac{b}{a^2+b}y^k$ \\\medskip

Choose $j_{k}^{2}$ (by some rule) \\\medskip

Compute $z^{k+1}=u^k -\frac{a}{a^{2}+b}\frac{1}{nL_{j_{k}^2}}\nabla f_{j_{k}^2}(y^{k})e_{j_{k}^{2}}$\\\medskip

\end{algorithmic}
\end{algorithm}\medskip

Just as in the non-strongly convex case, we extend the three algorithms ARCD, AGCD, and ASCD to the strongly convex case by using the rules \eqref{eq:non_rule_ARCD}, \eqref{eq:non_rule_AGCD}, and \eqref{eq:non_rule_ASCD} in Algorithm \ref{al:strframe}.  The following theorem presents our computational guarantee for ASCD for strongly convex case:\medskip

\begin{thm}\label{thm:strong}
Consider the Accelerated Semi-Greedy Coordinate Descent method for strongly convex case (Algorithm \ref{al:strframe} with rule \eqref{eq:non_rule_ASCD}). If $\ff$ is coordinate-wise $L$-smooth and $\mu$-strongly convex with respect to $\|\cdot\|_L$, it holds for all $k \ge 1$ that:
\begin{equation}\label{eq:conv_ASCD}
E_{\xi_{k}}\left[f(x^{k})-f^*+\tfrac{n^2}{2}(a^{2}+b)\|z^{k}-x^*\|_{L}^{2}\right]\le\left(1-\tfrac{\sqrt{\mu}}{n+\sqrt{\mu}}\right)^k\left(f(x^{0})-f^*+\tfrac{n^2}{2}(a^{2}+b)\|x^{0}-x^*\|_{L}^{2}\right).
\end{equation}
 \qed
\end{thm}

We provide a concise proof of Theorem \ref{thm:strong} in the Appendix.

%\section{Extensions}\label{sec:extension}
\section{Accelerated Greedy Coordinate Descent}\label{sec:AGCD}
In this section we discuss accelerated greedy coordinate descent (AGCD), which is Algorithm \ref{al:nonframe} with rule \eqref{eq:non_rule_AGCD}.  In the interest of clarity we limit our discussion to the non-strongly convex case.  We present a Lyapunov function argument which shows why the standard type of proof of accelerated gradient methods fails for AGCD, and we propose  a technical condition under which AGCD is guaranteed to have an $O(1/k^2)$ accelerated convergence rate.  Although there are no guarantees that the technical condition will hold for a given function $\ff$, we provide intuition as to why the technical condition ought to hold in most cases.

\subsection{Why AGCD fails (in theory)}\label{sec:AGCDfails}
The mainstream research community's interest in Nesterov's accelerated method \cite{nesterov1983method} started around 15 years ago; and yet even today most researchers struggle to find basic intuition as to what is really going on in accelerated methods.  Indeed, Nesterov's estimation sequence proof technique seems to work out arithmetically but with little fundamental intuition. There are many recent work trying to explain this acceleration phenomenon \cite{su2016differential}\cite{wilson2016lyapunov}\cite{hu2017dissipativity}\cite{lin2015universal}\cite{frostig2015regularizing}\cite{allen2014linear}
\cite{bubeck2015geometric}. A line of recent work has attempted to give a physical explanation of acceleration techniques by studying the continuous-time interpretation of accelerated gradient descent via dynamical systems, see \cite{su2016differential}, \cite{wilson2016lyapunov}, and \cite{hu2017dissipativity}.  In particular, \cite{su2016differential} introduced the continuous-time dynamical system model for accelerated gradient descent, and presented a convergence analysis using a Lyapunov energy function in the continuous-time setting. \cite{wilson2016lyapunov} studied discretizations of the continuous-time dynamical system, and also showed that Nesterov's estimation sequence analysis is equivalent to the Lyapunov energy function analysis in the dynamical system in the discrete-time setting. And \cite{hu2017dissipativity} presented an energy dissipation argument from control theory for understanding accelerated gradient descent.

In the discrete-time setting, one can construct a Lyapunov energy function of the form \cite{wilson2016lyapunov}:
\begin{equation}\label{eq:energy}
E_k = A_k(f(x^k)-f^*) + \tfrac{1}{2}\|x^*-z^k\|_L^2 \ ,
\end{equation}
where $A_k$ is a parameter sequence with $A_k \sim O(k^2)$, and one shows that $E_k$ is non-increasing in $k$, thereby yielding:$$
f(x^k)-f^*\le \frac{E_k}{A_k}\le \frac{E_0}{A_k} \sim O\left(\frac{1}{k^2}\right) \ .
$$
The proof techniques of acceleration methods such as  \cite{nesterov1983method}, \cite{tsengaccelerated} and \cite{allen2014linear}, as well as the recent proof techniques for accelerated randomized coordinate descent (such as \cite{nesterov2012efficiency}, \cite{lu2015complexity}, and \cite{fercoq2015accelerated}) can all be rewritten in the above form (up to expectation) each with slightly different parameter sequences $\{A_k\}$.

Now let us return to accelerated greedy coordinate descent.  Let us assume for simplicity that $L_1 = \cdots = L_n$ (as we can always do rescaling to achieve this condition).  Then the greedy coordinate $j_k^1$ is chosen as the coordinate of the gradient with the largest magnitude, which corresponds to the coordinate yielding the greatest guaranteed decrease in the objective function value.  However, in the proof of acceleration using the Lyapunov energy function, one needs to prove a decrease in $E_k$ \eqref{eq:energy} instead of a decrease in the objective function value $f(x^k)$. The coordinate $j_k^1$ is not necessarily the greedy coordinate for decreasing the energy function $E_k$ due to the presence of the second term $\|x^*-z^k\|_L^2$ in  \eqref{eq:energy}.  This explains why the greedy coordinate can fail to decrease $E_k$, at least in theory.  And because $x^*$ is not known when running AGCD, there does not seem to be any way to find the greedy descent coordinate for the energy function $E_k$.

%Based on this argument, we propose a new algorithm ASCD ( Algorithm Framework \ref{al:nonframe} with rule \eqref{eq:non_rule_ASCD}) in the next section in order to overcome the above difficulty. 

That is why in ASCD we use the greedy coordinate to perform the $x$-update (which corresponds to the fastest coordinate-wise decrease for the first term in energy function), while we choose a random coordinate to perform the $z$-update (which corresponds to the second term in the energy function); thereby mitigating the above problem in the case of ASCD.

In a concurrent paper \cite{locatello2018revisiting}, the authors develop computational theory for matching pursuit algorithms, which can be viewed as a generalized version of greedy coordinate descent where the directions do not need to form an orthogonal basis. The paper also develops an accelerated version of the matching pursuit algorithms, which turns out to be equivalent to the algorithm ASCD discussed here in the special case where the chosen directions are orthogonal. Although the focus in \cite{locatello2018revisiting} and in our paper are different -- \cite{locatello2018revisiting} is more focused on (accelerated) greedy direction updates along a certain linear subspace whereas our focus is on when and how one can accelerate greedy coordinate updates -- both of the works share a similar spirit and similar approaches in developing accelerated greedy methods.  Moreover, both works use a decoupling of the coordinate update for the $\{x^k\}$ sequence (with a greedy rule) and the $\{z^k\}$ sequence (with a randomized rule). In fact, \cite{locatello2018revisiting} is consistent with the argument in our paper as to why one cannot accelerate greedy coordinate descent in general.

\subsection{How to make AGCD work (in theory)}\label{sec:AGCDworks}

Here we propose the following technical condition under which the proof of acceleration of AGCD can be made to work.\medskip

\begin{cond}\label{cond:nonstr} There exists a positive constant $\gamma$ and an iteration number $K$ such that for all $k\ge K$ it holds that: 
\begin{equation}\label{dmd}\sum_{i=0}^k\frac{1}{\theta_i}\langle \nabla f(y^i), z^i-x^* \rangle \le \sum_{i=0}^k\frac{n\gamma}{\theta_i} \nabla_{j_i} f(y^i)(z^i_{j_i}-x^*_{j_i}) \ , \end{equation} where $j_i=\arg\max_{i}\tfrac{1}{\sqrt{L_i}}|\nabla_{i}f(y^{k})|$ is the greedy coordinate at iteration $i$.
\end{cond}

One can show that this condition is sufficient to prove an accelerated convergence rate $O(1/k^2)$ for AGCD.  Therefore let us take a close look at Technical Condition \ref{cond:nonstr}.  The condition considers the weighted sum (with weights $\frac{1}{\theta_i}\sim O(i^2)$) of the inner product of $\nabla f(y^k)$ and $z^k-x^*$, and the condition states that the inner product corresponding to the greedy coordinate (the right side above) is larger than the average of all coordinates in the inner product, by a factor of $\gamma$.  In the case of ARCD and ASCD, it is easy to show that Technical Condition  \ref{cond:nonstr} holds automatically up to expectation, with $\gamma = 1$.

Here is an informal explanation of why Technical Condition  \ref{cond:nonstr} ought to hold for most convex functions and most iterations of AGCD.  When $k$ is sufficiently large, the three sequence $\{x^k\}$, $\{y^k\}$ and $\{z^k\}$ ought to all converge to $x^*$ (which always happens in practice though lack of theoretical justification), whereby $z^k$ is close to $y^k$.  Thus we can instead consider the inner product $\langle \nabla f(y^k), y^k-x^* \rangle$ in \eqref{dmd}. Notice that for any coordinate $j$ it holds that $|y^k_j-x^*_j| \ge \frac{1}{L_j}|\nabla _j f(y^k)|$, and therefore $|\nabla_j f(y^k)\cdot (y^k_j-x^*_j)| \ge \frac{1}{L_j}|\nabla _j f(j^k)|^2$.  Now the greedy coordinate is chosen by $j_i := \arg\max_j \frac{1}{L_j}|\nabla _j f(j^k)|^2$, and therefore it is reasonably likely that in most cases the greedy coordinate will yield a better product than the average of the components of the inner product.

% By gradient inequality and smoothness, we can bound the inner product by 
%$$ f(y^k)-f(x^*)+\tfrac{1}{2}\|y^k-x^*\|_L^2 \ge \langle \nabla f(y^k), y^k-x^* \rangle \ge f(y^k)-f(x^*) ,$$
%Thus $\langle \nabla f(y^k), y^k-x^* \rangle$ is close to $f(y^k)-f(x^*)$ as the difference is bounded by $\|y^k-x^*\|_L^2$ which is in second order. Meanwhile, the greedy coordinate corresponds to the coordinate for $f(y^k)-f(x^*)$ to change, thus greedy coordinate of the inner product should be larger than the average intuitively. 

The above is not a rigorous argument, and we can likely design some worst-case functions for which Technical Condition  \ref{cond:nonstr} fails.  But the above argument provides some intuition as to why the condition ought to hold in most cases, thereby yielding the observed improvement of AGCD as compared with ARCD that we will shortly present in Section \ref{sec:numerical}, where we also observe that Technical Condition  \ref{cond:nonstr} holds empirically on all of our problem instances.

With a slight change in the proof of Theorem \ref{thm:nonstrong}, we can show the following result:\medskip
\begin{thm}\label{thm:nonstrongtech}
Consider the Accelerated Greedy Coordinate Descent (Algorithm \ref{al:nonframe} with rule \eqref{eq:non_rule_AGCD}). If $\ff$ is coordinate-wise $L$-smooth and satisfies Technical Condition \ref{cond:nonstr} with constant $\gamma \le 1$ and iteration number $K$, then it holds for all $k \ge K$ that:
\begin{equation}\label{eq:conv_ASCD_non_tech}
f(x^{k})-f(x^{*})\le \tfrac{2n^{2}\gamma}{(k+1)^{2}}\|x^{*}-x^{0}\|_L^{2} \ . 
\end{equation}
\end{thm}
We note that if $\gamma<1$ (which we always observe in practice), then AGCD will have a better convergence guarantee than ARCD.  
\begin{rem}
The arguments in Section \ref{sec:AGCDfails} and Section \ref{sec:AGCDworks} also work for strongly convex case, albeit with suitable minor modifications.
\end{rem}

%{\color{red} We can get similar condition and result for strongly convex case as well. The linear rate becomes $1-\tfrac{\sqrt{\mu}}{n\sqrt{\gamma}+\sqrt{\mu}}$, but the algorithm also need to be changed with constant $\gamma$. Do we want to include that result too, or just mention similar results can be got for strongly convex case?}
%
%We can also show the following result with some change of the proof of Theorem \ref{thm:strong}:
%
%\begin{thm}\label{thm:strong}
%Consider the Accelerated Greedy Coordinate Descent for Strongly Convex Case (Algorithm \ref{al:strframe} with rule \eqref{eq:str_rule_AGCD}) and with parameter $a=\frac{\sqrt{\mu}}{n\sqrt{\gamma}+\sqrt{\mu}}$ and $b=\frac{\mu a}{n^2 \gamma}$. If $\ff$ is coordinate-wise $L$-smooth, $\mu$-strongly convex with respect to $\|\cdot\|_L$, and it satisfies Technical Condition \ref{cond:nonstr} with constant $\frac{1}{n}\le\gamma\le 1$ and iteration number $K$ then for all $k \ge 1$ and $x\in Q$, and the following inequality holds:
%
%\begin{equation}\label{eq:conv_ASCD_non}
%f(x^{k})-f^*+\tfrac{n^2}{2}(a^{2}+b)\|z^{k}-x^*\|_{L}^{2}\le\left(1-\tfrac{\sqrt{\mu}}{n\sqrt{\gamma}+\sqrt{\mu}}\right)^k\left(f(x^{0})-f^*+\tfrac{n^2\gamma}{2}(a^{2}+b)\|x^{0}-x^*\|_{L}^{2}\right).
%\end{equation}
% \qed
%\end{thm}
%

\section{Numerical Experiments}\label{sec:numerical}

\subsection{Linear Regression}
We consider solving synthetic instances of the linear regression model with least-squares objective function:
\begin{equation*}
f^* \ := \ \min_{\beta \in \mathbb{R}^p} f(\beta) := \|y-X\beta\|_2^2 \ 
\end{equation*}

using ASCD, ARCD and AGCD, where the mechanism for generating the data $(y,X)$ and the algorithm implementation details are described in the supplementary materials.  Figure \ref{Fig:LinearR} shows the optimality gap versus time (in seconds) for solving different instances of linear regression with different condition numbers of the matrix $X^TX$ using ASCD, ARCD and AGCD. In each plot, the vertical axis is the objective value optimality gap $f(\beta^k)-f^*$ in log scale, and the horizontal axis is the running time in seconds. Each column corresponds to an instance with the prescribed condition number $\kappa$ of $X^T X$, where $\kappa=\infty$ means that the minimum eigenvalue of $X^T X$ is $0$. The first row of plots is for Algorithm Framework \ref{al:nonframe} which is ignorant of any strong convexity information. The second row of plots is for Algorithm Framework \ref{al:strframe}, which uses given strong convexity information.  And because the linear regression optimization problem is quadratic, it is straightforward to compute $\kappa$ as well as the true parameter $\mu$ for the instances where $\kappa>0$.  The last column of the figure corresponds to $\kappa = \infty$, and in this instance we set $\mu$ using the smallest positive eigenvalue of $X^TX$, which can be shown to work in theory since all relevant problem computations are invariant in the nullspace of $X$. 

Here we see in Figure \ref{Fig:LinearR} that AGCD and ASCD consistently have superior performance over ARCD for both the non-strongly convex case and the strongly convex case, with ASCD performing almost as well as AGCD in most instances.

We remark that the behavior of any convex function near the optimal solution is similar to the quadratic function defined by the Hessian at the optimum, and therefore the above numerical experiments show promise that AGCD and ASCD are likely to outperform ARCD asymptotically for any twice-differentiable convex function.

%The following are some observation.
%\begin{itemize}
%\item AGCD and ASCD always win ARCD. ASCD is comparable to AGCD in some cases while not all.
%\item With the increase of $\kappa$, the running time to get certain optimality gap increases correspondingly.
%\item In general, the strongly convex methods (Algorithm Framework \ref{al:strframe}) work better than non-strongly convex methods (Algorithm Framework \ref{al:nonframe}).
%\item For the case $\kappa = \infty$, both approaches work well. Actually the decreasing rate in this case is related to the positive condition number, i.e. the ratio between minimal non-zero eigenvalue to the maximal eigenvalue.
%\item The behavior of any convex function near the optimal solution is similar to a quadratic function, thus this shows the asymptotic convergence results for these three algorithms.
%
%\end{itemize}

\begin{figure}
\centering
\begin{tabular}{c|M{30mm}M{30mm}M{30mm}M{30mm}}
\toprule
        & $\kappa=10^{2}$ & $\kappa=10^{3}$ & $\kappa=10^4$ & $\kappa=\infty$ \\
\midrule
$\begin{array}{c} \mbox{Algorithm} \cr \mbox{Framework~1} \cr \mbox{(non-strongly} \cr \mbox{convex)} \end{array}$& \includegraphics[width=35mm,height=35mm]{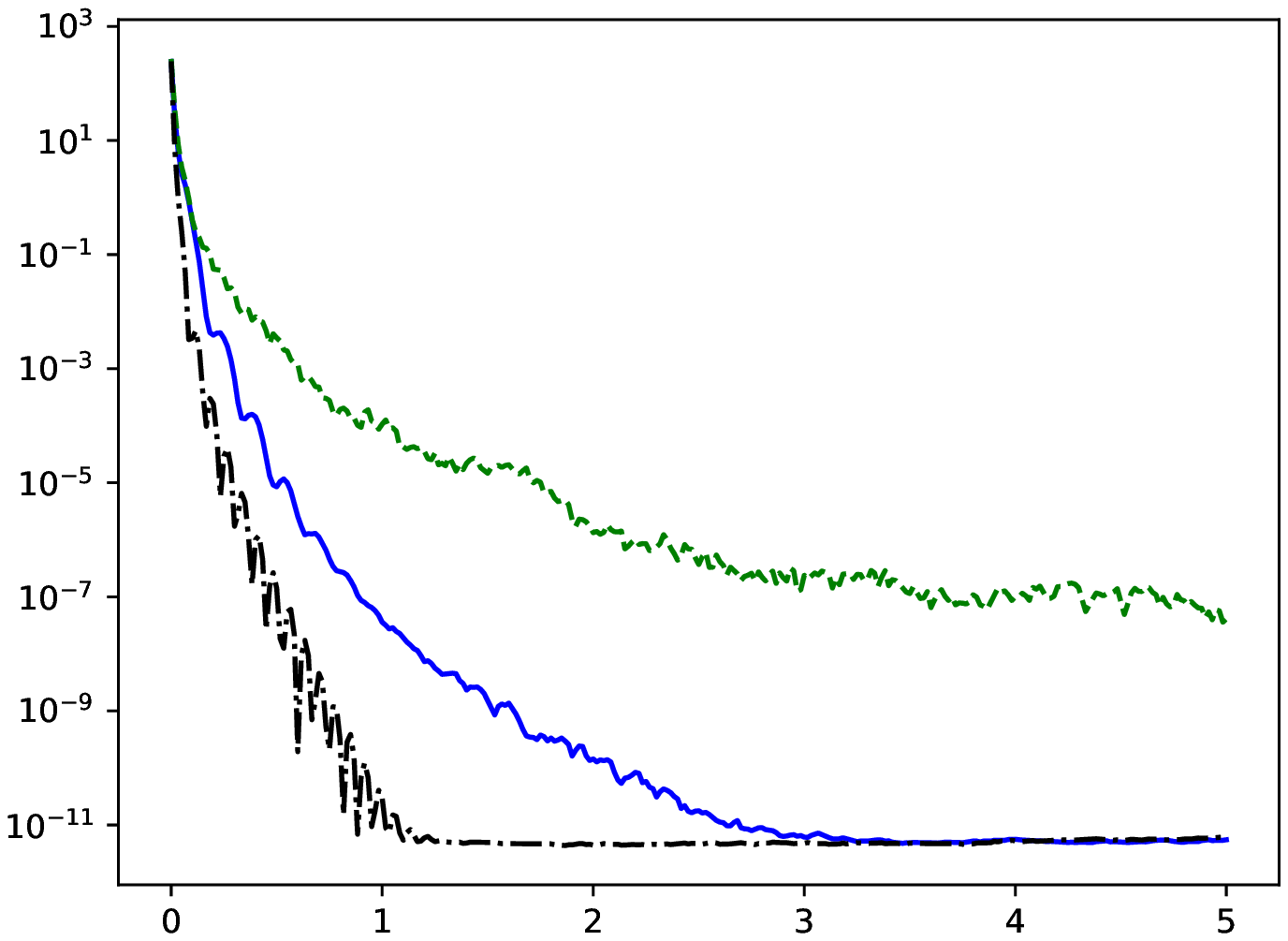} & \includegraphics[width=35mm,height=35mm]{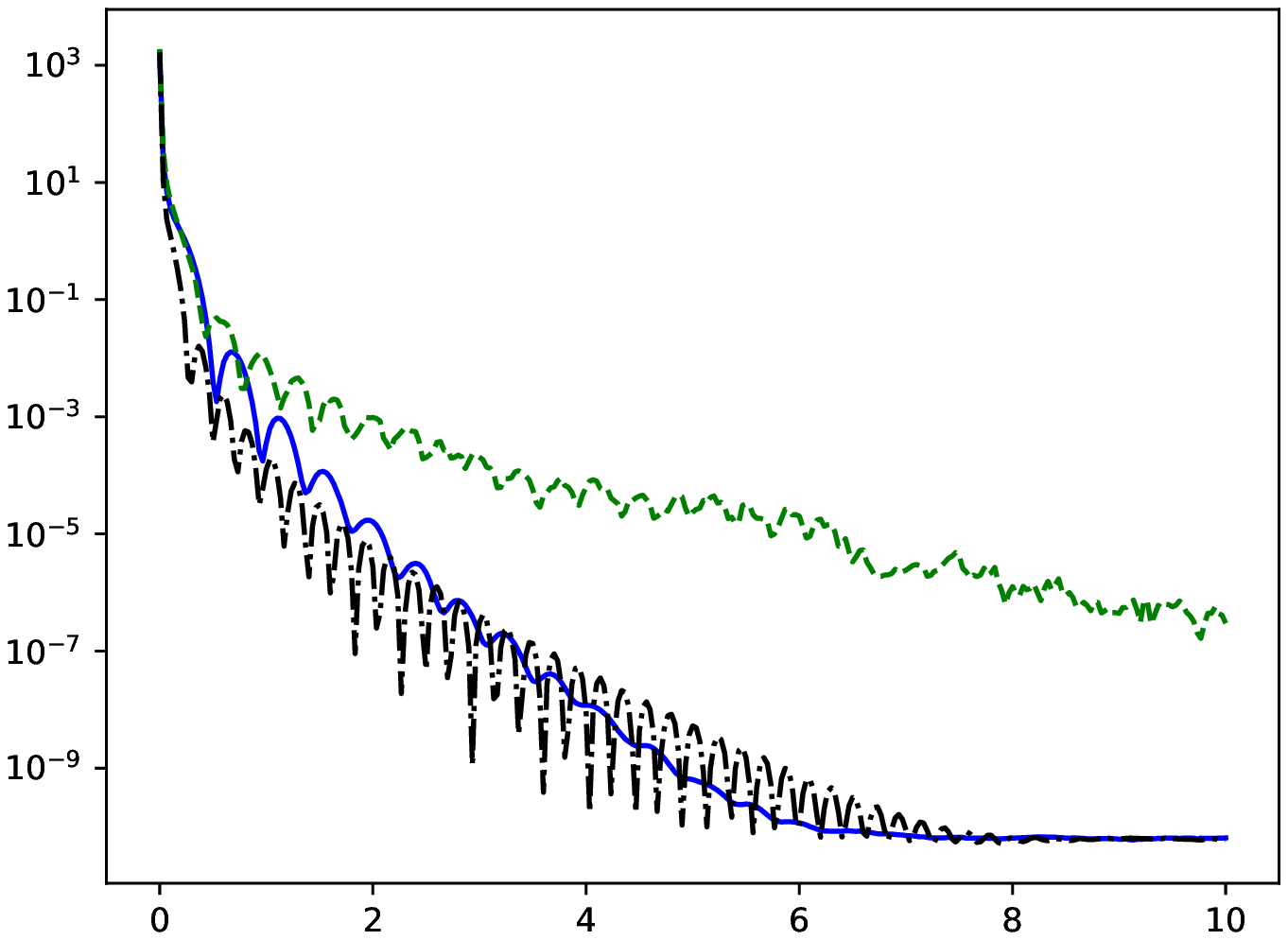} & \includegraphics[width=35mm,height=35mm]{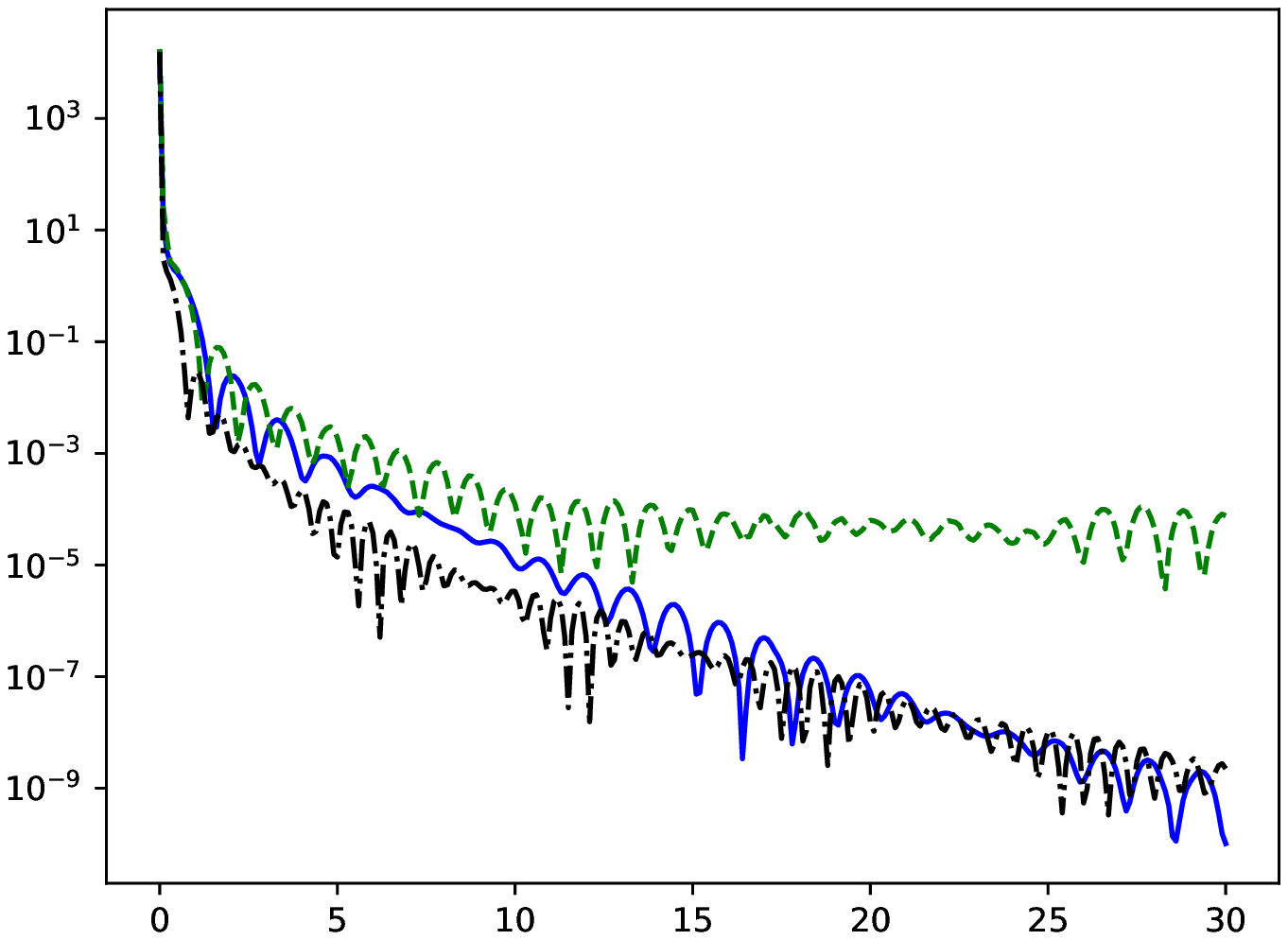} & \includegraphics[width=35mm,height=35mm]{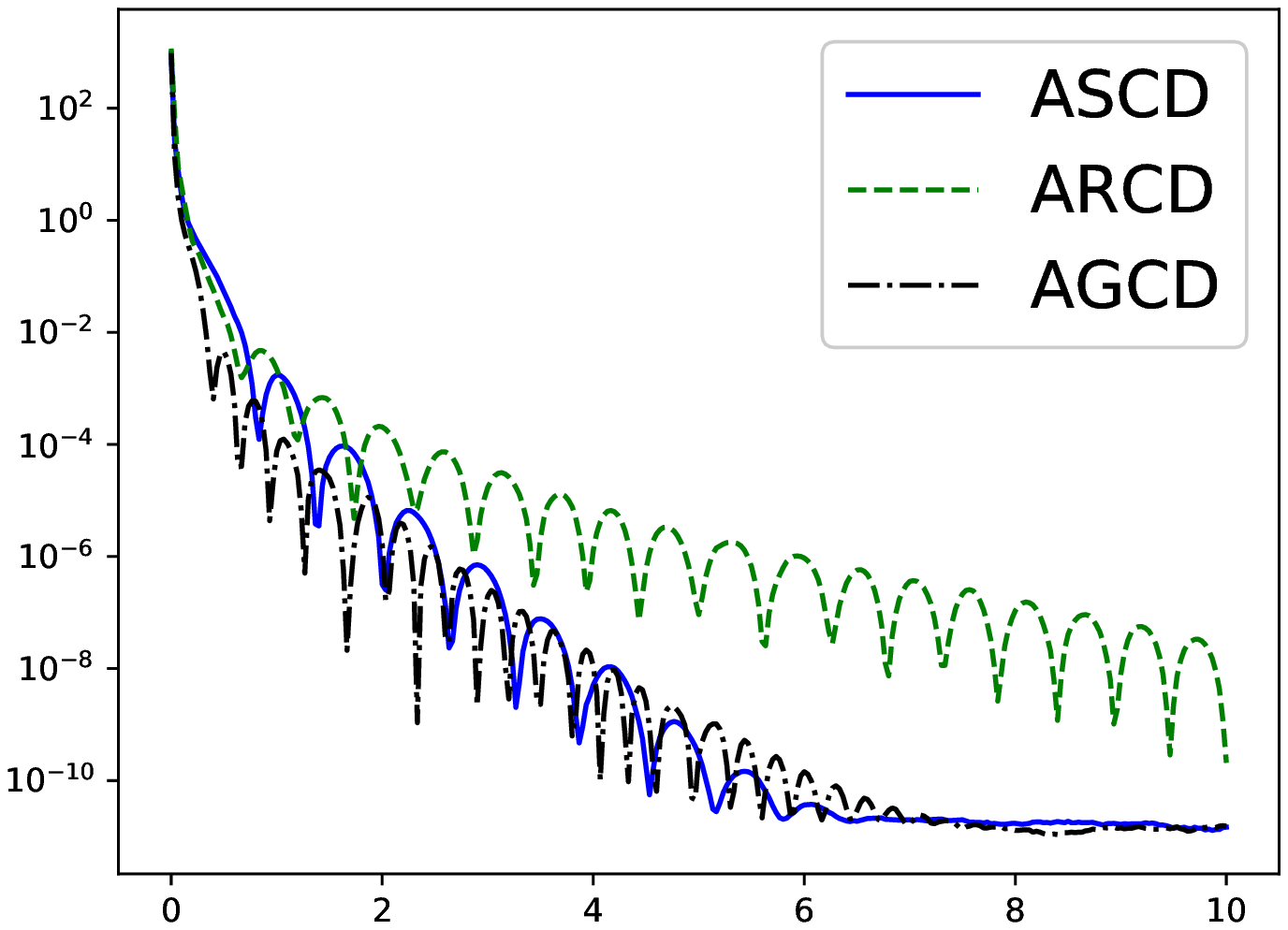} \\
$\begin{array}{c} \mbox{Algorithm} \cr \mbox{Framework~2} \cr \mbox{(strongly} \cr \mbox{convex)} \end{array}$ & \includegraphics[width=35mm,height=35mm]{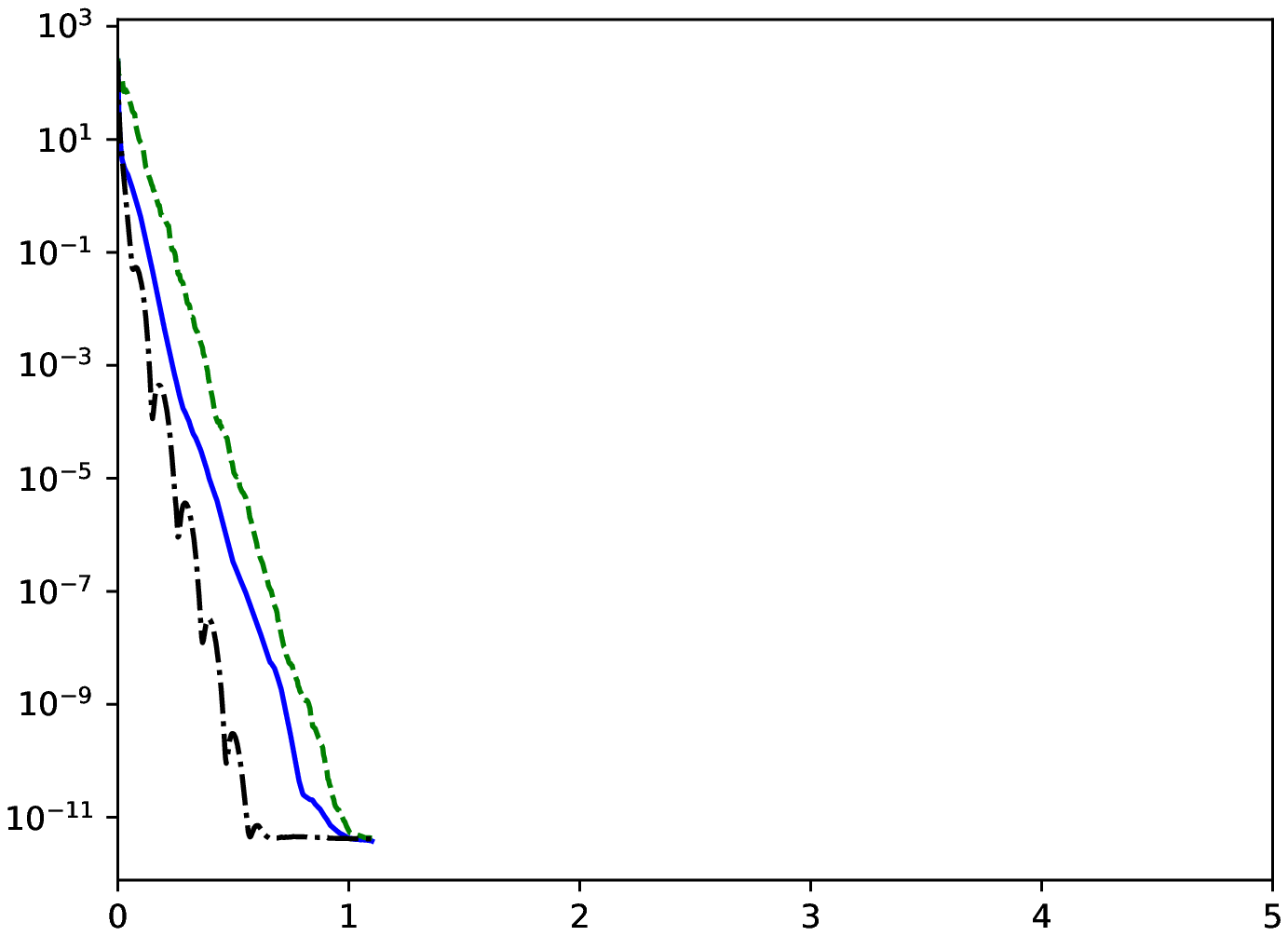} & \includegraphics[width=35mm,height=35mm]{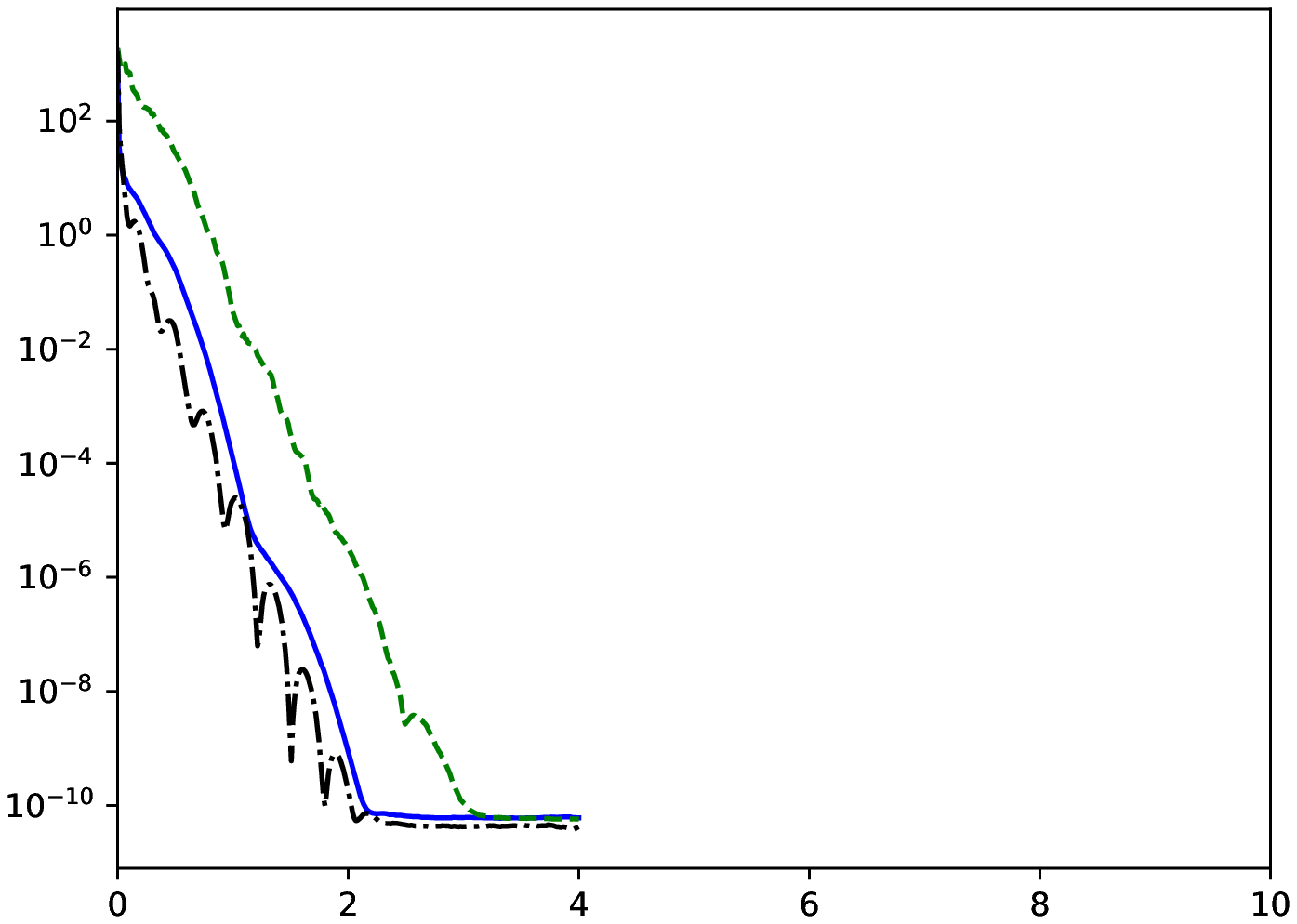} & \includegraphics[width=35mm,height=35mm]{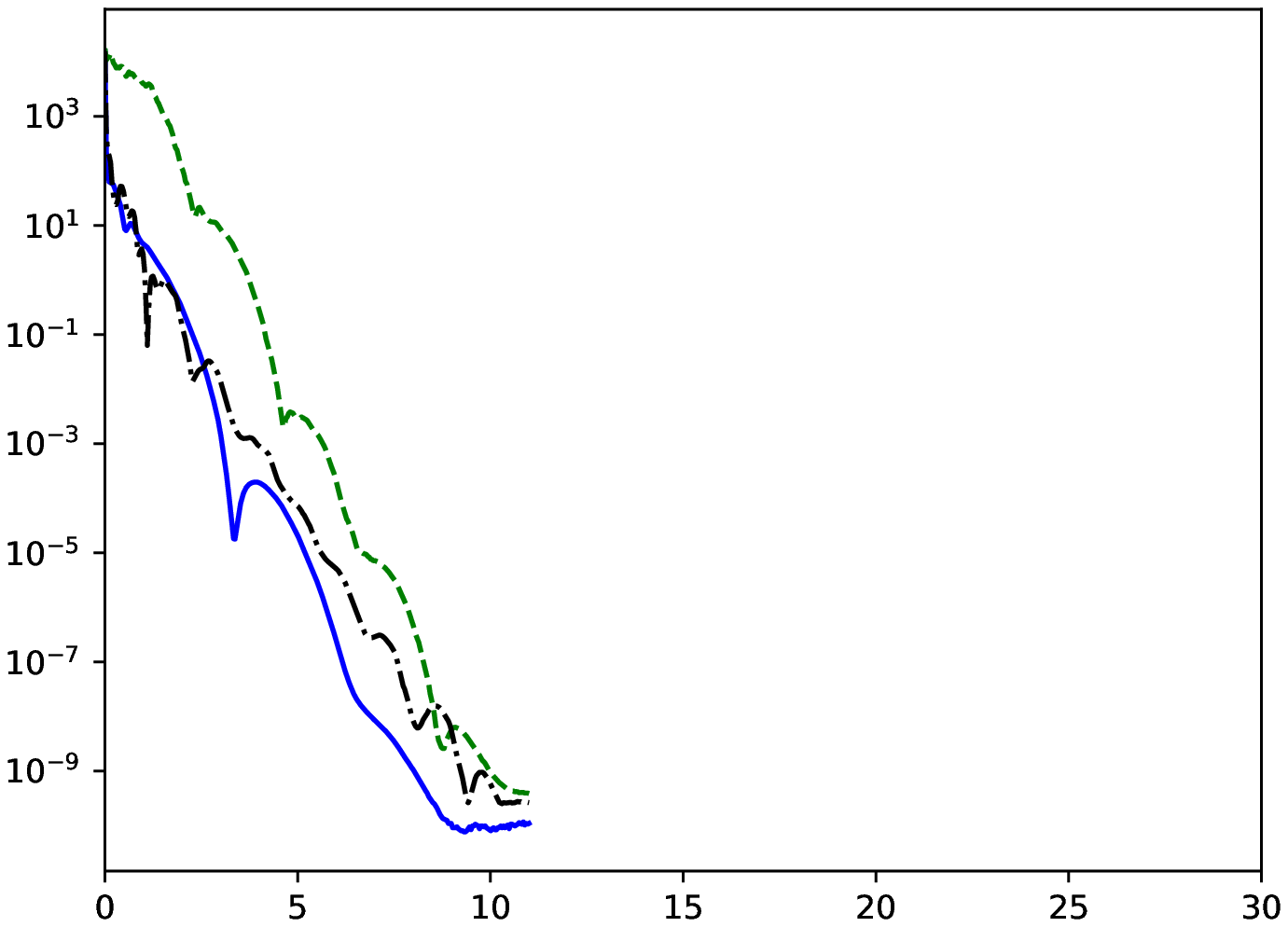} & \includegraphics[width=35mm,height=35mm]{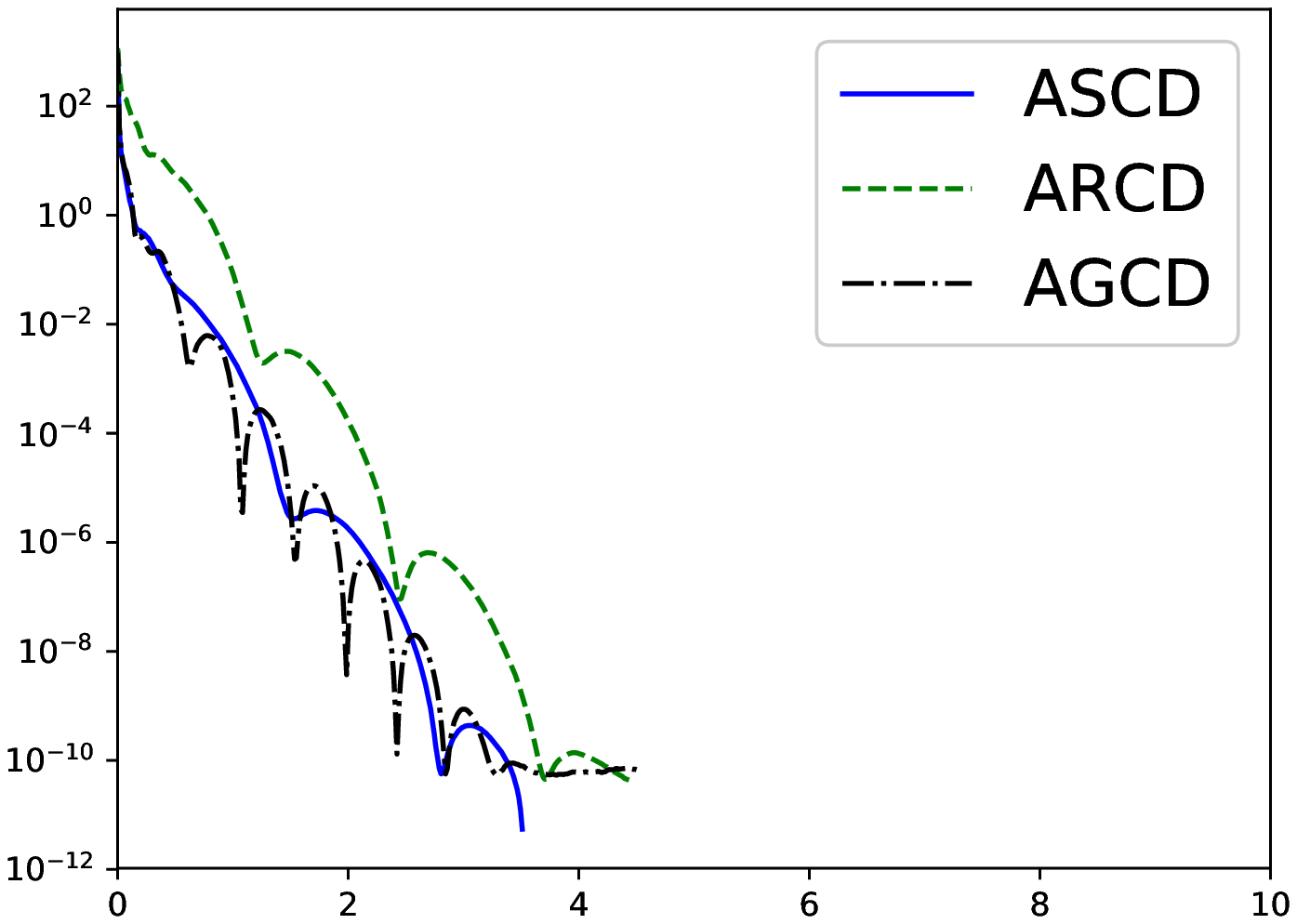} \\

\bottomrule
\end{tabular}
\caption{Plots showing the optimality gap versus run-time (in seconds) for synthetic linear regression problems solved by ASCD, ARCD and AGCD. }
\label{Fig:LinearR}
\end{figure}

\subsection{Logistic Regression}
\begin{figure}
\centering
\begin{tabular}{c|M{30mm}M{30mm}M{30mm}M{30mm}}
\toprule
Dataset & $\bar\mu=10^{-3}$ & $\bar\mu=10^{-5}$ & $\bar\mu=10^{-7}$ & $\bar\mu=0$ \\
\midrule
w1a & \includegraphics[width=35mm,height=35mm]{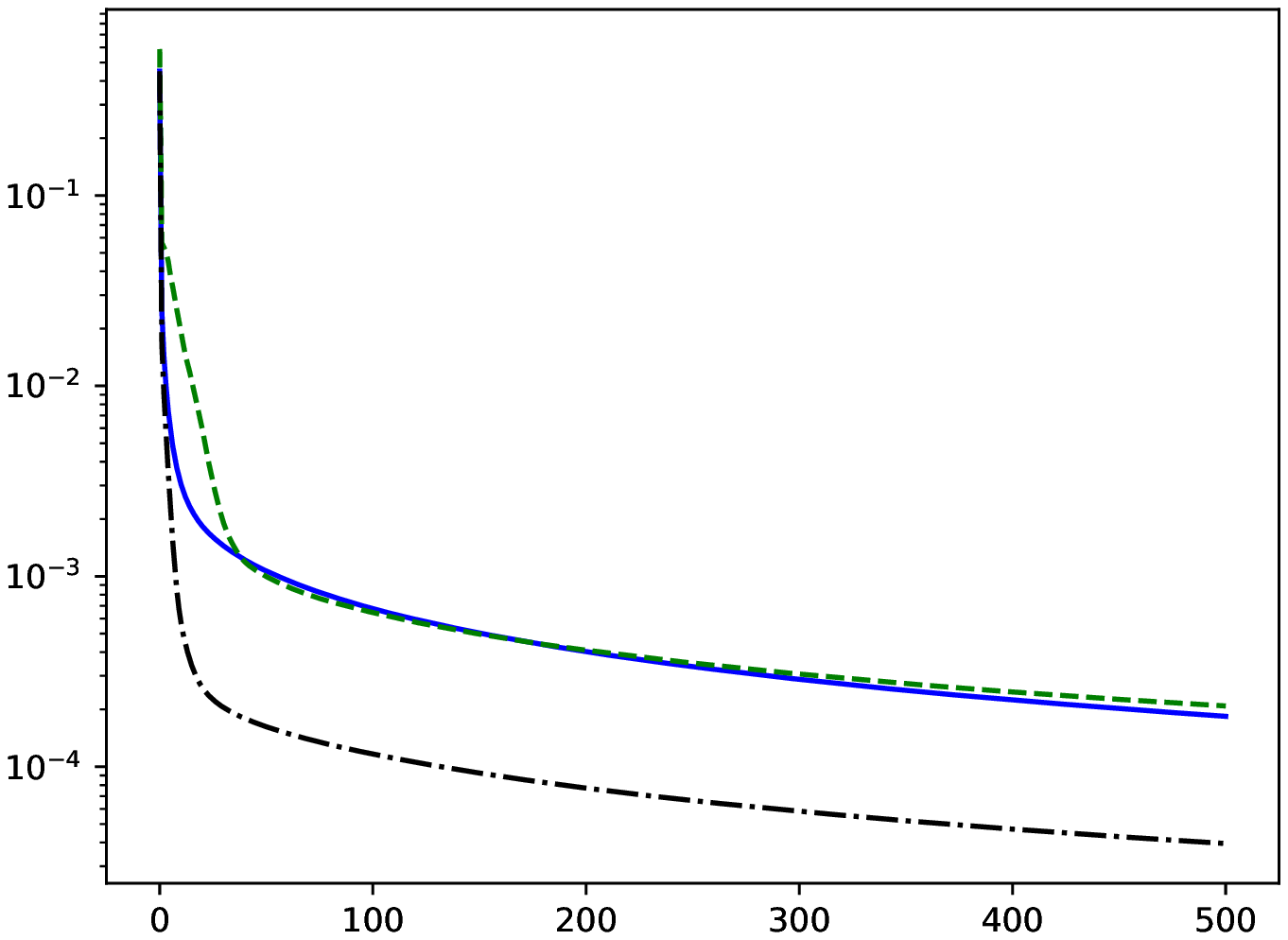} & \includegraphics[width=35mm,height=35mm]{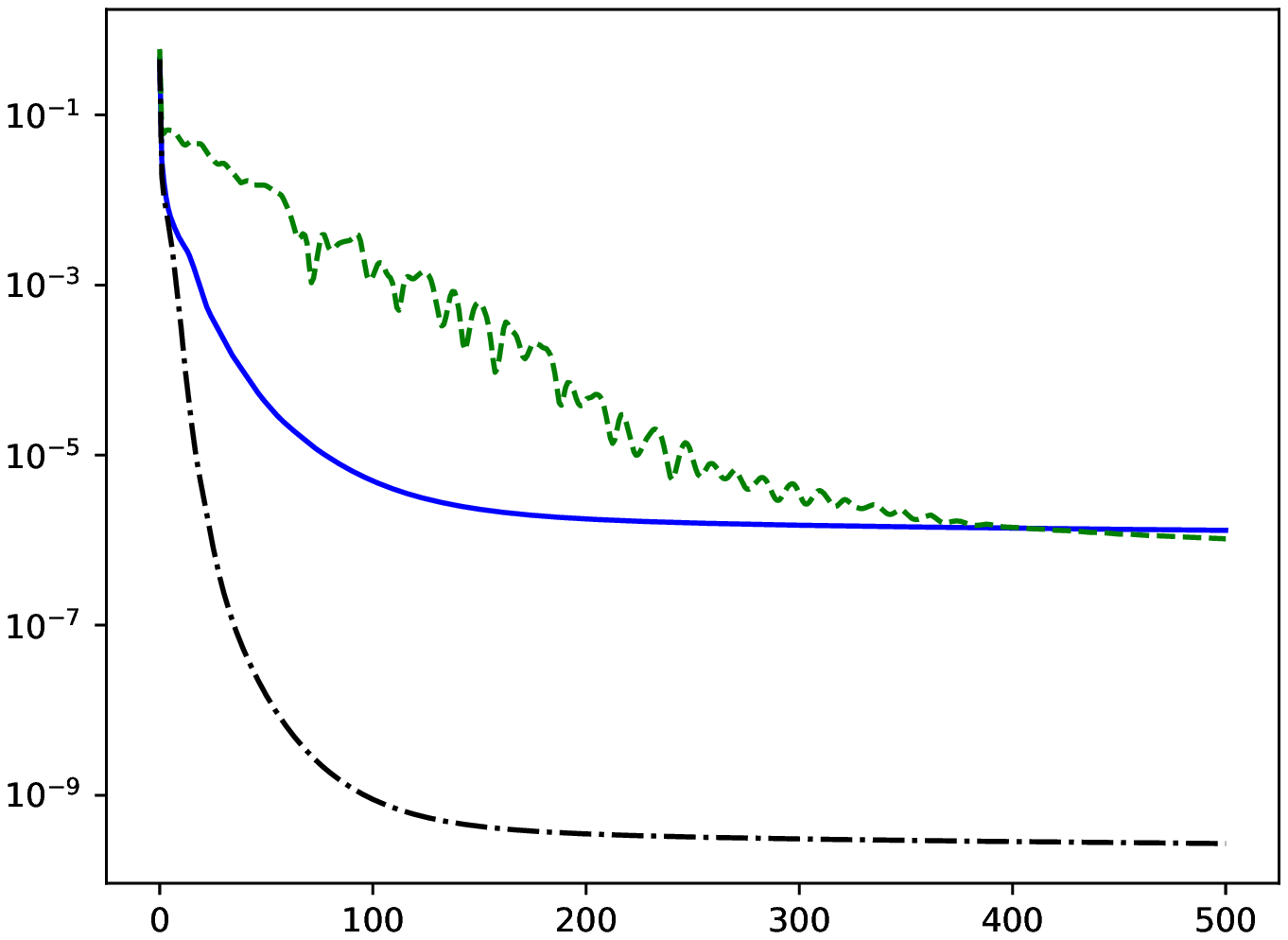} & \includegraphics[width=35mm,height=35mm]{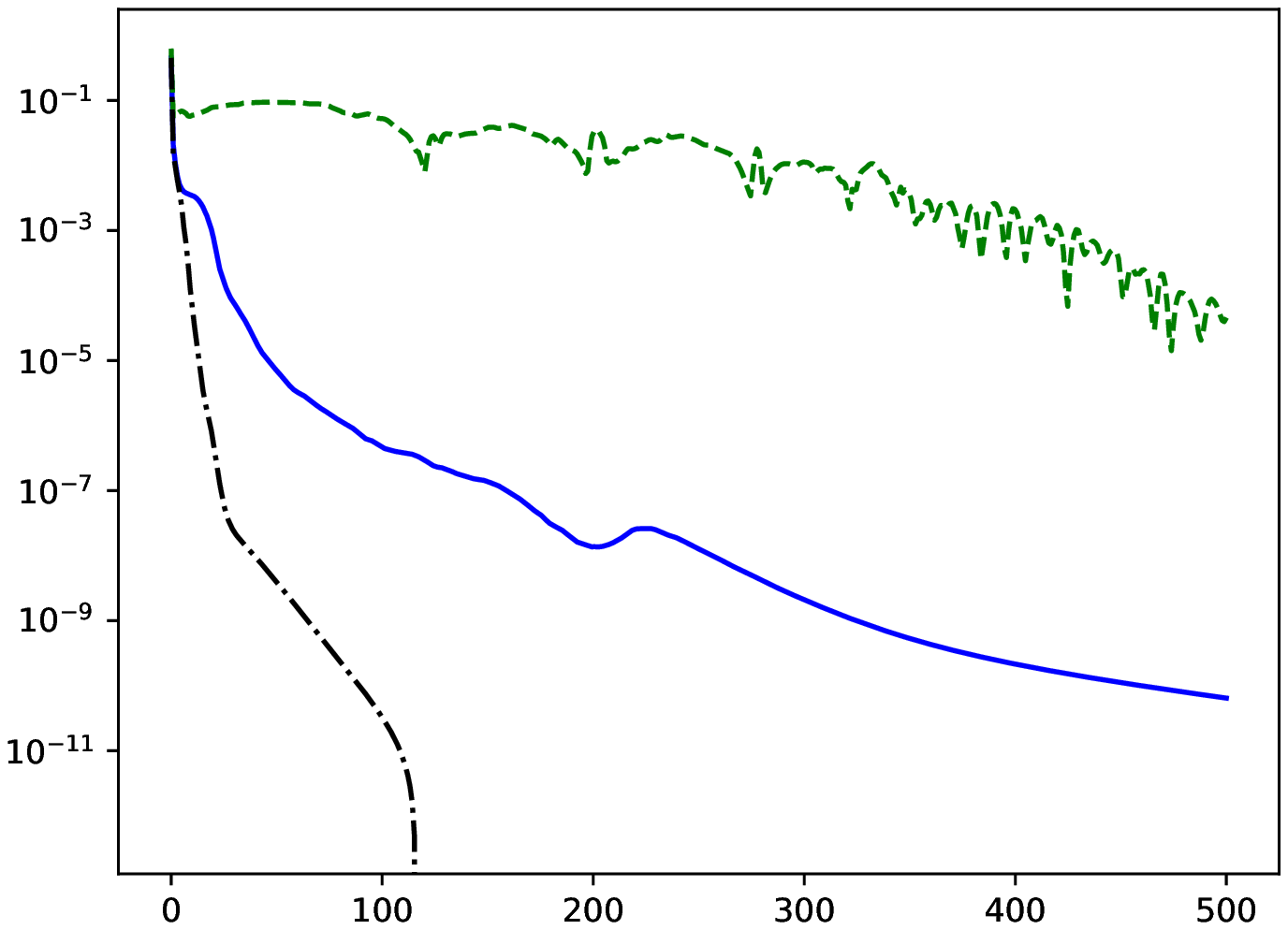} & \includegraphics[width=35mm,height=35mm]{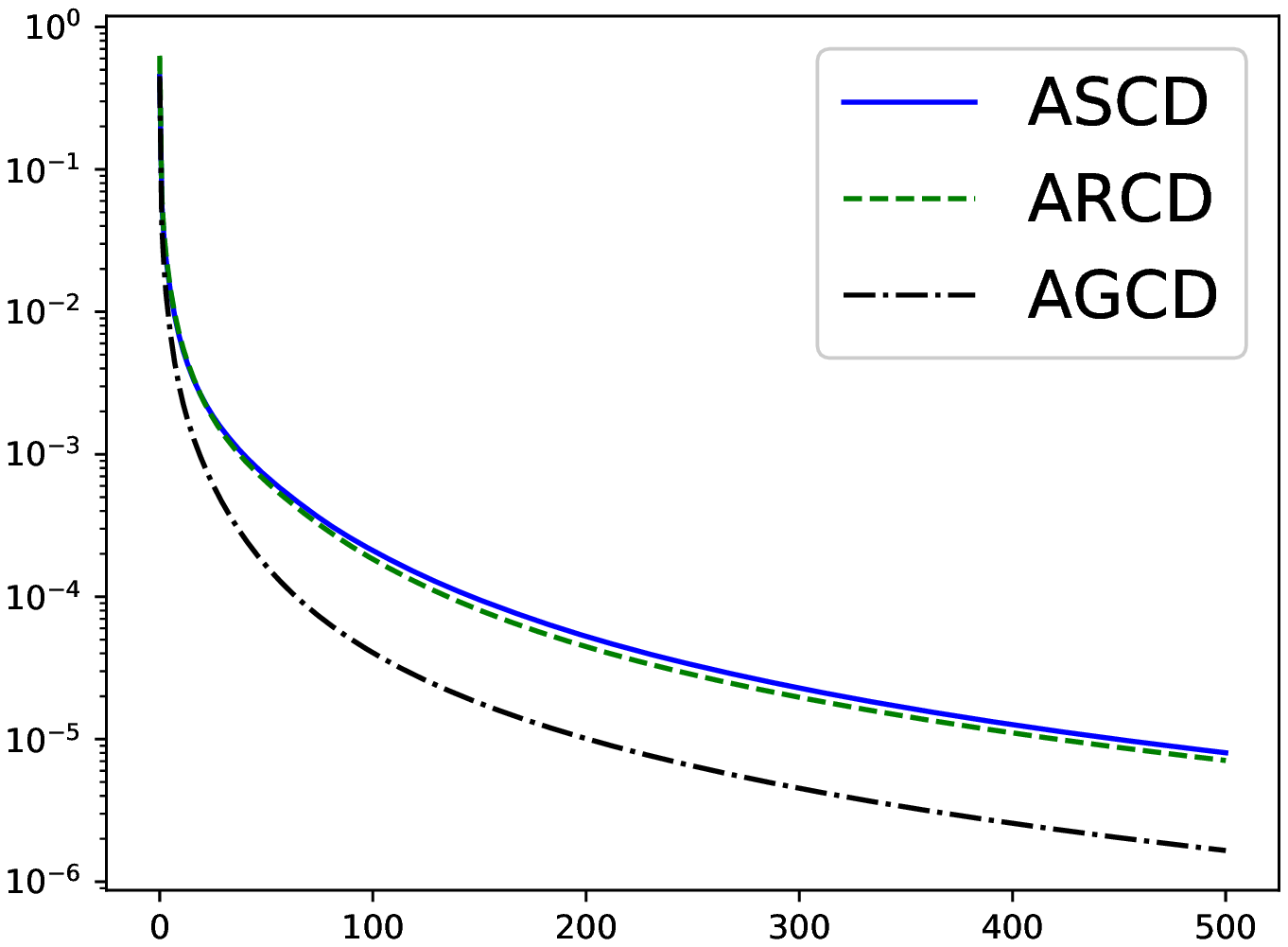} \\
a1a & \includegraphics[width=35mm,height=35mm]{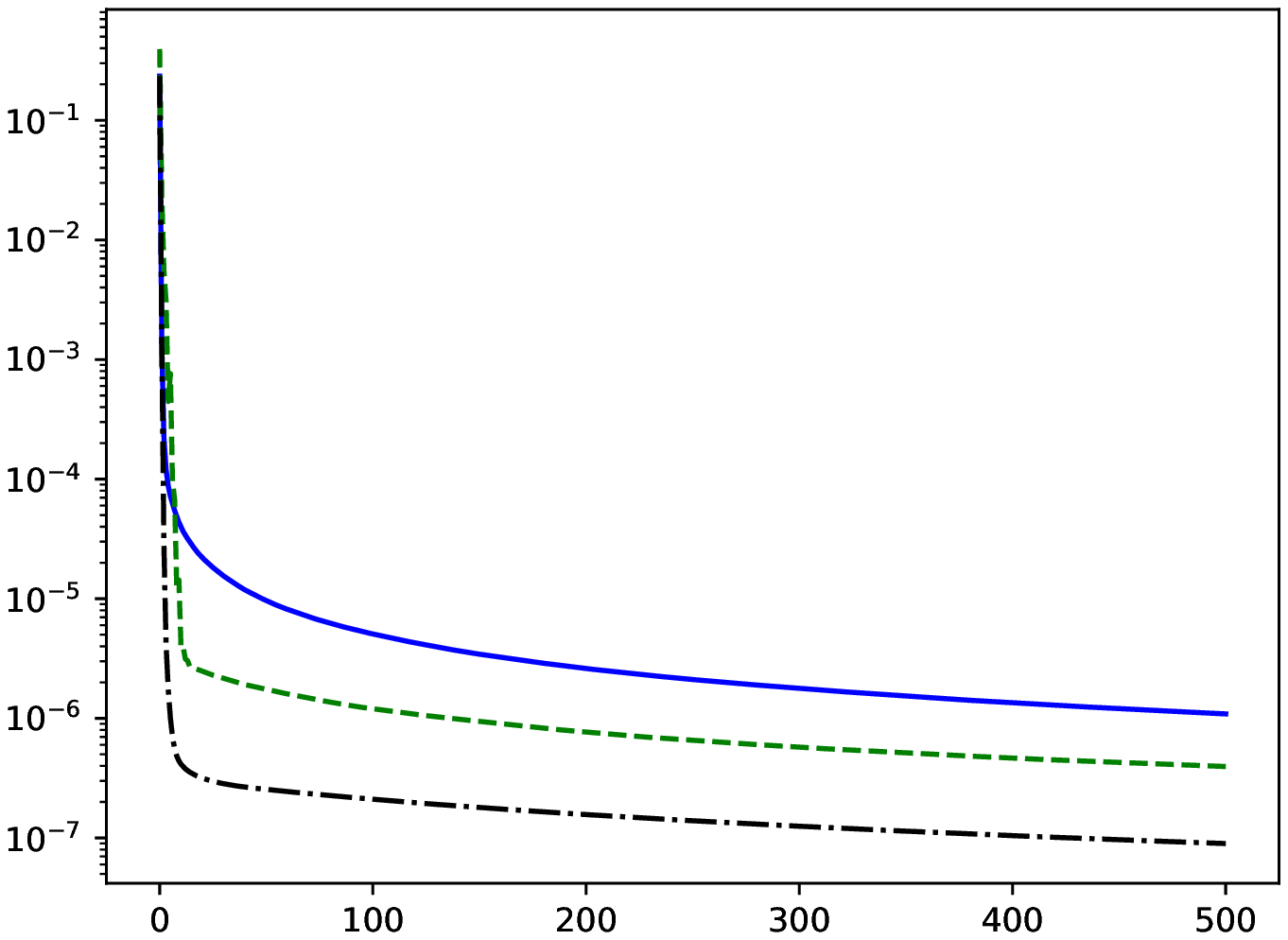} & \includegraphics[width=35mm,height=35mm]{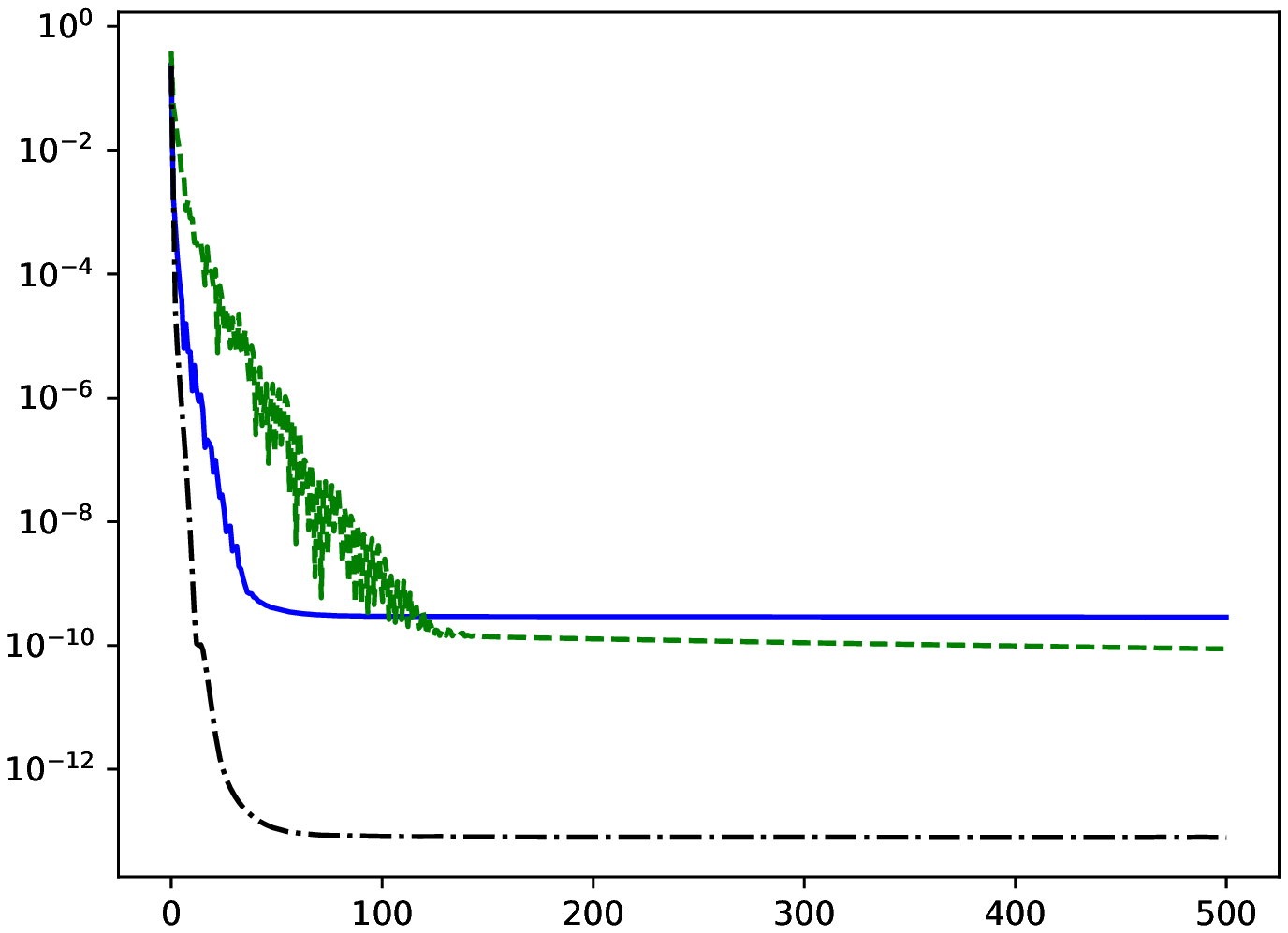} & \includegraphics[width=35mm,height=35mm]{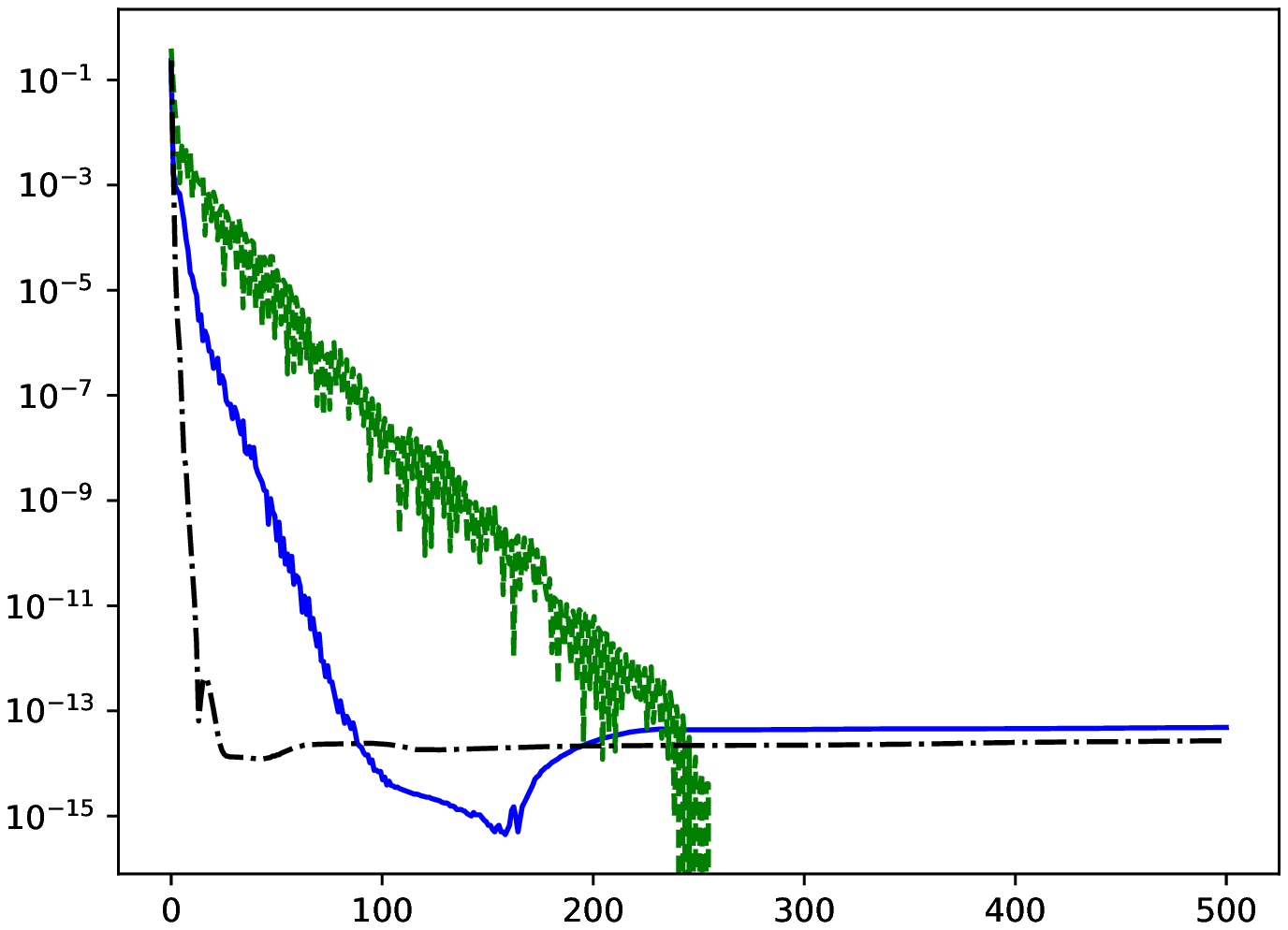} & \includegraphics[width=35mm,height=35mm]{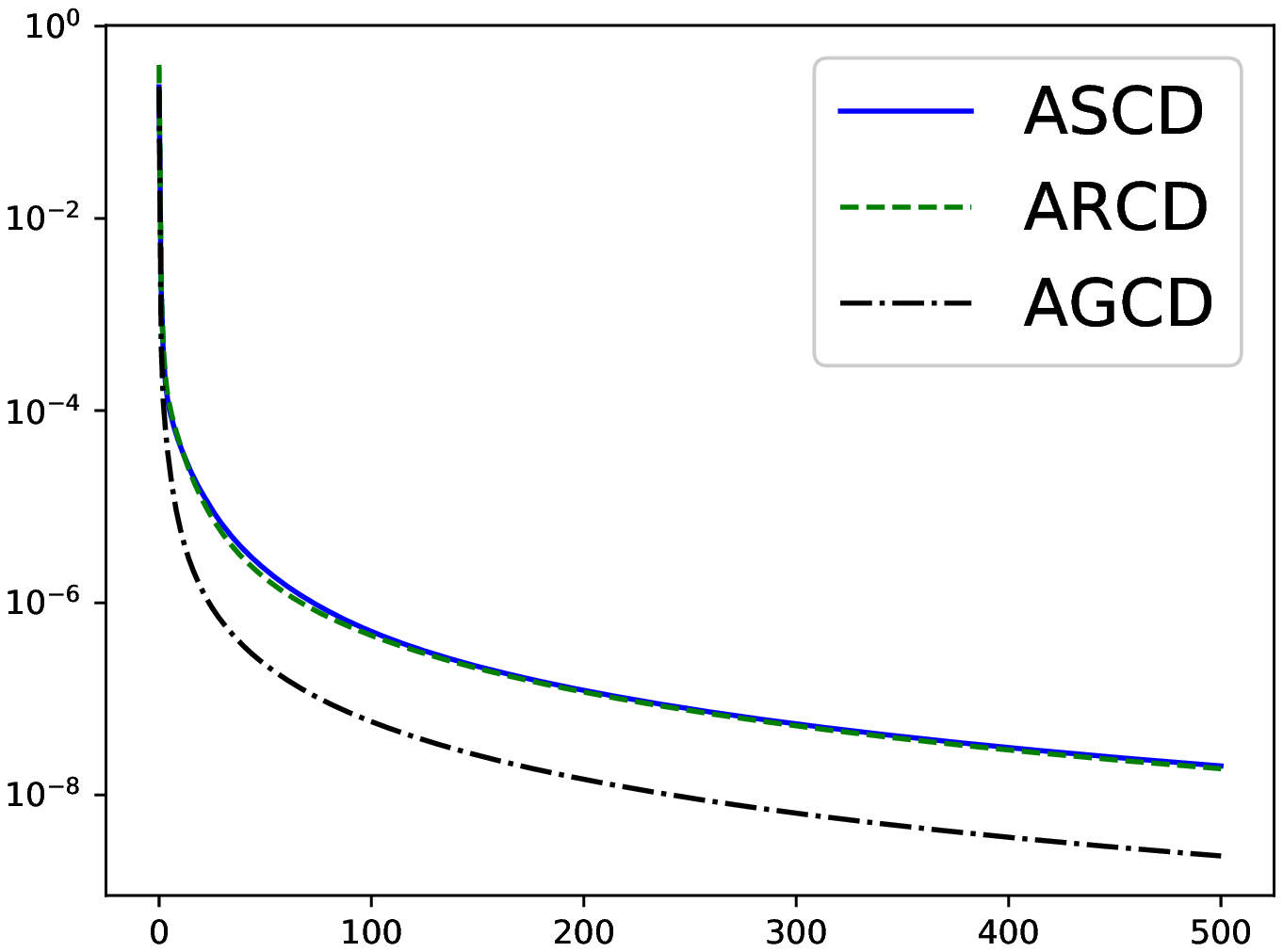} \\

\bottomrule
\end{tabular}
\caption{Plots showing the optimality gap versus run-time (in seconds) for the logistic regression instances w1a and a1a in LIBSVM, solved by ASCD, ARCD and AGCD. }
\label{Fig:LogisticR}
\end{figure}

Here we consider solving instances of the logistic regression loss minimization problem:
\begin{equation}
f^* \ := \ \min_{\beta \in \mathbb{R}^p} f(\beta):=\frac{1}{n} \sum_{i=1}^n \log(1+\exp(-y_i\beta^T x_i)) \ ,
\end{equation}
using ASCD, ARCD and AGCD, where $\{x_i, y_i\}$ is the feature-response pair for the $i$-th data point and $y_i \in \{-1,1\}$.    Although the loss function $f(\beta)$ is not in general strongly convex, it is essentially locally strongly convex around the optimum but with unknown strong convexity parameter $\bar \mu$. And although we do not know the local strong convexity parameter $\bar \mu$, we can still run the strongly convex algorithm (Algorithm Framework \ref{al:strframe}) by assigning a value of $\bar \mu$ that is hopefully close to the actual value.  Using this strategy, we solved a large number of logistic regression instances from LIBSVM \cite{chang2011libsvm}.  Figure \ref{Fig:LogisticR} shows the optimality gap versus time (in seconds) for solving two of these instances, namely w1a and a1a, which were chosen here because the performance of the algorithms on these two instances is representative of others in LIBSVM.  In each plot, the vertical axis is the objective value optimality gap $f(\beta^k)-f^*$ in log scale, and the horizontal axis is the running time in seconds. Each column corresponds to a different assigned value of the local strong convexity parameter $\bar \mu$.  The right-most column in the figure uses the assignment $\bar\mu=0$, in which case the algorithms are implemented as in the non-strongly convex case (Algorithm Framework \ref{al:nonframe}).

Here we see in Figure \ref{Fig:LogisticR} that AGCD always has superior performance as compared to either ASCD and ARCD.  In the relevant range of optimality gaps ($\le 10^{-9}$), ASCD typically outperforms ARCD for smaller values of the assigned strong convexity parameter $\bar \mu$.  However, the performance of ASCD and ARCD are essentially the same when no strong convexity is presumed.

%The following are some observations:
%\begin{itemize}
%\item AGCD is the best method out of the three methods. For strongly convex framework, ASCD wins ARCD for in most of the cases, especially at the beginning to run the algorithm. For non-strongly convex framework, ASCD is tied with ARCD.
%\item With a proper $\mu$ for strongly convex framework, the convergence will be faster.
%\end{itemize}

Last of all, we attempt to estimate the parameter $\gamma$ that arises in Technical Condition \ref{cond:nonstr} for AGCD in several of the datasets in SVMLIB.   Although for small $k$, the ratio between $\sum_{i=0}^k\frac{1}{\theta_i}\langle \nabla f(y^i), z^i-x^*\rangle$ and $ \sum_{i=0}^k\frac{n}{\theta_i} \nabla_{j_i} f(y^i)(z^i_{j_i}-x^*_{j_i})$ can fluctuate widely, this ratio stabilizes after a number of iterations in all of our numerical tests.  From Technical Condition \ref{cond:nonstr}, we know that $\gamma$ is the upper bound of such ratio for all $k\ge K$ for some large enough value of $K$. Table \ref{tab:gamma} presents the observed values of $\gamma$ for all $K\ge \bar{K}:=5,000$.  Recalling from Theorem \ref{thm:nonstrongtech} that the $\gamma$ value represents how much better AGCD can perform compared with ARCD in terms of computational guarantees, we see from Table \ref{tab:gamma} that AGCD should outperform ARCD for these representative instances -- and indeed this is what we observe in practice in our tests.

\begin{table}
\centering
\caption{Largest observed values of $\gamma$ for five different datasets in LIBSVM for $k \ge \bar{K}:=5000$.}
\begin{tabular}{|r|c|c|c|c|c|}
\hline
Dataset: & w1a & a1a & heart & madelon & rcv1 \\
\hline
$\gamma$: & 0.25 & 0.17 & 0.413 & 0.24 & 0.016 \\
\hline
\end{tabular}
\label{tab:gamma}
\end{table}

\appendix
\section{Appendix}

%\subsection{An Intuitive Explanation on Why Direct Extension of AGCD May Fail although ARCD Works \color{red}{(Put this part in appendix?)}}
%Here we only consider the non-strongly convex case, while this explanation can be extended to strongly convex case as well.
%
%At first, let us look at Accelerated Gradient Descent (AGD)
%
%\begin{algorithm}
%\caption{Accelerated Gradient Descent (non-strongly convex case)}\label{al:AGDnon}
%\begin{algorithmic}
%\STATE {\bf Initialize.}  Initialize with $x^0$, set $z^0 \leftarrow x^0$, and define the sequence $\{\theta_k\}$ as follows: $\theta_{0}=1$, and $\tfrac{1-\theta_{k}}{\theta_{k}^{2}}=\tfrac{1}{\theta_{k-1}^{2}}$ recursively for $k=1, 2, \ldots$. Assume $f(\cdot)$ is $L$-smooth for known and given $L$. \\
%
%$ \ $
%
%At iteration $k$ :\\
%\STATE  {\bf Perform Updates.}  
%
%Define $y^{k}:=(1-\theta_{k})x^{k}+\theta_{k}z^{k}$, \\\medskip
%
%Choose $j_{k}^{1} \in \arg\max_{i}\tfrac{1}{\sqrt{L_i}}|\nabla_{i}f(y^{k})|$,\\\medskip
%
%Compute $x^{k+1}:=y^{k}-\tfrac{1}{L_{j_k^1}}\nabla_{j_{k}^{1}}f(y^{k}) e_{j_k^1}$, \\\medskip
%
%Sample $j_{k}^{2} \sim {\mathcal U}[1,\cdots,n]$, \\\medskip
%
%Compute $z^{k+1}:=z^{k}-\tfrac{1}{nL_{j_{k}^2}\theta_{k}}\nabla_{j_{k}^{2}}f(y^{k}) e_{j_k^2}$ \ . \\\medskip
%
%\end{algorithmic}
%\end{algorithm}\medskip
%
%
%
%$ \ $
%For accelerated method, 

\subsection{Proof of Theorem \ref{thm:strong}}

In order to prove Theorem 3.1, we first prove the following three lemmas:\medskip

\begin{lem}\label{lem:geometry}
\[
a^{2}\|x^{*}-z^{k}\|_{L}^{2}+b\|x^{*}-y^{k}\|_{L}^{2}=(a^{2}+b)\|x^{*}-u^{k}\|_{L}^{2}+\tfrac{a^{2}b}{a^{2}+b}\|y^{k}-z^{k}\|_{L}^{2} \ .
\]
\end{lem}
{\bf Proof:}
\begin{equation}
\begin{array}{cl}
 & (a^{2}+b)\|x^{*}-u^{k}\|_{L}^{2}+\frac{a^{2}b}{a^{2}+b}\|y^{k}-z^{k}\|_{L}^{2}\\ \\ 
= & (a^{2}+b)\left(\|x^{*}\|_{L}^{2}-2\langle x^{*},u^{k}\rangle_{L}+\|u^{k}\|_{L}^{2}\right)+\frac{a^{2}b}{a^{2}+b}\|y^{k}-z^{k}\|_{L}^{2}\\ \\
= & (a^{2}+b)\|x^{*}\|_{L}^{2}-2\langle x^{*},a^{2}z^{k}+by^{k}\rangle_{L}+\frac{1}{a^{2}+b}\|a^{2}z^{k}+by^{k}\|_{L}^{2}+\frac{a^{2}b}{a^{2}+b}\|y^{k}-z^{k}\|_{L}^{2}\\ \\
= & (a^{2}+b)\|x^{*}\|_{L}^{2}-2\langle x^{*},a^{2}z^{k}+by^{k}\rangle_{L}+a^{2}\|z^{k}\|_{L}^{2}+b\|y^{k}\|_{L}^{2}\\ \\
= & a^{2}\|x^{*}-z^{k}\|_{L}^{2}+b\|x^{*}-y^{k}\|_{L}^{2} \ ,
\end{array}
\end{equation}
where the second equality utilizes $u^{k}=\frac{a^{2}}{a^{2}+b}z^{k}+\frac{b}{a^{2}+b}y^{k}$
and the other equalities are just mathematical manipulations. \qed\medskip

\begin{lem}\label{lem:expectation} Define $t^{k+1}:=u^{k}-\frac{a}{a^{2}+b}\frac{1}{n}\bl^{-1} \nabla f(y^{k})$, then
\[
\|x^*-t^{k+1}\|_{L}^{2}-\|x^*-u^{k}\|_{L}^{2}=n E_{j_k^2}\left[\|x^*-z^{k+1}\|_{L}^{2}-\|x^*-u^{k}\|_{L}^{2}\right]
\]
\end{lem}
{\bf Proof:}
\begin{equation}
\begin{array}{lcl}
\|x^*-t^{k+1}\|_{L}^{2}-\|x^*-u^{k}\|_{L}^{2} & = &2\langle x^*-u^{k},u^{k}-t^{k+1}\rangle_{L}+\|u^{k}-t^{k+1}\|_{L}^{2}\\ \\
 & = &2n E_{j_k^2}\left[\langle x^*-u^{k},u^{k}-z^{k+1}\rangle_{L}+\|u^{k}-z^{k+1}\|_{L}^{2}\right]\\ \\
 & = &n E_{j_k^2}\left[\|x^*-z^{k+1}\|_{L}^{2}-\|x^*-u^{k}\|_{L}^{2}\right],
\end{array}
\end{equation}
where the second equality is from the relationship of $t^{k+1}=u^k-\frac{a}{a^{2}+b}\frac{1}{n}\bl^{-1}\nabla f(y^{k})$
and $z^{k+1}=u^k-\frac{a}{a^{2}+b}\frac{1}{n L_{j_k^2}}\nabla_{j_k^2} f(y^{k})e_{j_k^2}$, and the first and third equations are just rearrangement. \qed\medskip

\begin{lem}\label{lem:ab_relation}
\[
a^{2}\le(1-a)(a^{2}+b)\ .
\]
\end{lem}

{\bf Proof:} Remember that $b=\frac{\mu a}{n^{2}}$ and $a>0$, thus the above inequality
is equivalent to 
\[
a\le(1-a)(a+\tfrac{\mu}{n^{2}})\ .
\]
Substituting $a=\tfrac{\sqrt{\mu}}{n+\sqrt{\mu}}$, the above inequality
becomes
\[
\tfrac{\sqrt{\mu}}{n}\le\tfrac{\sqrt{\mu}}{n+\sqrt{\mu}}+\tfrac{\mu}{n^{2}} \ .
\]
We furnish the proof by noting $\tfrac{\sqrt{\mu}}{n}-\tfrac{\sqrt{\mu}}{n+\sqrt{\mu}}=\tfrac{\mu}{n(n+\sqrt{\mu})} \le \tfrac{\mu}{n^2}$. \qed

{\bf Proof of Theorem \ref{thm:strong}:} Recall that $t^{k+1}=u^{k}-\frac{a}{a^{2}+b}\frac{1}{n}\bl^{-1}\nabla f(y^{k})$,
then it is easy to check that
\[
t^{k+1}=\arg\min_{z} \ a\langle\nabla f(y^{k}),z-z^{k}\rangle+\frac{n}{2}a^{2}\|z-z^{k}\|_{L}^{2}+\frac{n}{2}b\|z-y^{k}\|_{L}^{2}
\]
by writing the optimality conditions of the right-hand side.

We have
\begin{equation}\label{eq:group1}
\begin{array}{cl}
& f(x^{k+1})-f(y^{k}) \\ \\
 \le &\langle\nabla f(y^{k}),x^{k+1}-y^{k}\rangle+\frac{1}{2}\|x^{k+1}-y^{k}\|_{L}^{2}\\ \\
  = &-\frac{1}{2L_{j_{k}^{1}}}\left(\nabla_{j_{k}^{1}}f(y^{k})\right)^{2}\\ \\
  \le &-\frac{1}{2n}\|\nabla f(y^{k})\|_{L^{-1}}^{2}\\ \\
  \le & a\langle\nabla f(y^{k}),t^{k+1}-z^{k}\rangle+\frac{n}{2}a^{2}\|t^{k+1}-z^{k}\|_{L}^{2}\\ \\
  = &a\langle\nabla f(y^{k}),x^*-z^{k}\rangle+\frac{n}{2}a^{2}\|x^*-z^{k}\|_{L}^{2}-\frac{n}{2}a^{2}\|x^*-t^{k+1}\|_{L}^{2}+\frac{n}{2}b\|x^*-y^{k}\|_{L}^{2} \\ \\
  &-\frac{n}{2}b\|y^{k}-t^{k+1}\|_{L}^{2}-\frac{n}{2}b\|t^{k+1}-x^*\|_{L}^{2}\\ \\
  \le & a\langle\nabla f(y^{k}),x^*-z^{k}\rangle+\frac{n}{2}a^{2}\|x^*-z^{k}\|_{L}^{2}-\frac{n}{2}a^{2}\|x^*-t^{k+1}\|_{L}^{2}+\frac{n}{2}b\|x^*-y^{k}\|_{L}^{2}-\frac{n}{2}b\|t^{k+1}-x^*\|_{L}^{2}\\ \\
  = &a\langle\nabla f(y^{k}),x^*-z^{k}\rangle+\frac{n}{2}(a^{2}+b)\left(\|x^*-u^{k}\|_{L}^{2}-\|x^*-t^{k+1}\|_{L}^{2}\right)+\frac{n}{2}\frac{a^{2}b}{a^{2}+b}\|y^{k}-z^{k}\|_{L}^{2}\\ \\
  = &a\langle\nabla f(y^{k}),x^*-z^{k}\rangle+\frac{n^{2}}{2}(a^{2}+b) E_{j_{k}^{2}}\left[\|x^*-u^{k}\|_{L}^{2}-\|x^*-z^{k+1}\|_{L}^{2}\right]+\frac{n}{2}\frac{a^{2}b}{a^{2}+b}\|y^{k}-z^{k}\|_{L}^{2}\\ \\
  \le &a\langle\nabla f(y^{k}),x^*-z^{k}\rangle+\frac{n^{2}}{2}(a^{2}+b) E_{j_{k}^{2}}\left[\|x^*-u^{k}\|_{L}^{2}-\|x^*-z^{k+1}\|_{L}^{2}\right]+\frac{n^{2}}{2}\frac{a^{2}b}{a^{2}+b}\|y^{k}-z^{k}\|_{L}^{2}\\ \\
  = &a\langle\nabla f(y^{k}),x^*-z^{k}\rangle+\frac{n^{2}}{2}\left((a^{2}+b)\|x^*-u^{k}\|_{L}^{2}+\frac{a^{2}b}{a^{2}+b}\|y^{k}-z^{k}\|_{L}^{2}\right)-\frac{n^{2}}{2}(a^{2}+b)\EE_{j_{k}^{2}}\left[\|x^*-z^{k+1}\|_{L}^{2}\right]\\ \\
  = &a\langle\nabla f(y^{k}),x^*-z^{k}\rangle+\frac{n^{2}}{2}\left(a^{2}\|x^*-z^{k}\|_{L}^{2}+b\|x^*-y^{k}\|_{L}^{2}\right)-\frac{n^{2}}{2}(a^{2}+b) E_{j_{k}^{2}}\left[\|x^*-z^{k+1}\|_{L}^{2}\right],
\end{array}
\end{equation}
where the first inequality is due to coordinate-wise smoothness, the
first equality utilizes $x^{k+1}=y^{k}-\frac{1}{L_{j_{k}^{1}}}\nabla_{j_{k}^{1}}f(y^{k})e_{j_k^1}$,
the second inequality follows from the fact that $j_{k}^{1}$ is the
greedy coordinate of $\nabla f(y^{k})$ in the $\| \cdot\|_{L^{-1}}$ norm, the third
inequality follows from the basic inequality $\|v\|_L^2 + \|w\|_{L^{-1}}^2 \ge 2\langle v,w \rangle$ for all $v,w$, the second
equality is from Three Point Property by noticing 
\[
t^{k+1}=\arg\min_{z} \ a\langle\nabla f(y^{k}),z-z^{k}\rangle+\frac{n}{2}a^{2}\|z-z^{k}\|_{L}^{2}+\frac{n}{2}b\|z-y^{k}\|_{L}^{2} \ ,
\]
the third equality follows from Lemma 3.1, and the fourth and sixth equalities
each utilize Lemma 3.2.

On the other hand, by strong convexity we have 
\begin{equation}\label{eq:group2}
\begin{array}{lcl}
f(y^{k})-f(x^*) & \le &\langle\nabla f(y^{k}),y^{k}-x^*\rangle-\frac{1}{2}\mu\|y^{k}-x^*\|_{L}^{2}\\ \\
 & = &\langle\nabla f(y^{k}),y^{k}-z^{k}\rangle+\langle\nabla f(y^{k}),z^{k}-x^*\rangle-\frac{1}{2}\mu\|y^{k}-x^*\|_{L}^{2}\\ \\
 & = &\frac{1-a}{a}\langle\nabla f(y^{k}),x^{k}-y^{k}\rangle+\langle\nabla f(y^{k}),z^{k}-x^*\rangle-\frac{1}{2}\mu\|y^{k}-x^*\|_{L}^{2}\\ \\
 & \le &\frac{1-a}{a}\left(f(x^{k})-f(y^{k})\right)+\langle\nabla f(y^{k}),z^{k}-x^*\rangle-\frac{1}{2}\mu\|y^{k}-x^*\|_{L}^{2},
\end{array}
\end{equation}
where the second equality uses the fact that $y^{k}=(1-a)x^{k}+az^{k}$
and the last inequality is from the gradient inequality.

By rearranging \eqref{eq:group2}, we obtain
\begin{equation}\label{eq:group2rearange}
f(y^{k})-f(x^*)\le(1-a)\left(f(x^{k})-f(x^*)\right)+a\langle\nabla f(y^{k}),z^{k}-x^*\rangle-\frac{1}{2}\mu a\|y^{k}-x^*\|_{L}^{2} \ .
\end{equation}

Notice that $b=\frac{\mu a}{n^{2}}$ and $a^{2}\le(1-a)(a^{2}+b)$
following from Lemma 3.3. Thus summing up \eqref{eq:group1} and \eqref{eq:group2rearange} leads to
\begin{equation}
\begin{array}{cl}
& E_{j_{k}^{2}}\left[f(x^{k+1})-f(x^*)+\frac{n^{2}}{2}(a^{2}+b)\|z^{k+1}-x^*\|_{L}^{2}\right] \\ \\
 \le &(1-a)\left(f(x^{k})-f(x^*)\right)+\frac{n^{2}}{2}a^{2}\|z^{k}-x^*\|_{L}^{2}\\ \\
 \le &(1-a)\left(f(x^{k})-f(x^*)+\frac{n^{2}}{2}(a^{2}+b)\|z^{k}-x^*\|_{L}^{2}\right) \ ,
\end{array}
\end{equation}
which furnishes the proof using a telescoping series. \qed

\section{More Material on the Numerical Experiments}\label{sec:more}
\subsection{Implementation Detail}\label{sec:Imple}
To be consistent with the notation in statistics and machine learning we use $p$ to denote the dimension of the variables in the optimization problems describing linear and logistic regression.  Then the per-iteration computation cost of AGCD and ASCD is dominated by three computations: (i) $p$-dimensional vector operations (such as in computing $y^k$ using $x^k$ and $z^k$), (ii) computation of the gradient $\nabla f(\cdot)$, and (iii) computation of the maximum (weighted) magnitude coordinate of the gradient $\nabla f(\cdot)$.  \cite{lee2013efficient} proposed an efficient way to avoid (i) by changing variables.  Distinct from the dual approaches discussed in \cite{lee2013efficient}, \cite{lin2015accelerated}, and \cite{allen2016even}, here we only consider the primal problem in the regime $n>p$, and therefore the cost of (ii) dominates the cost of (i) in these cases. For this reason in our numerical experiments we use the simple implementation of ARCD proposed by Nesterov \cite{nesterov2012efficiency} and which we adopt for AGCD and ASCD as well. We note that both the randomized methods and the greedy methods can take advantage of the efficient calculations proposed in \cite{lee2013efficient} as well.

For the linear regression experiments we focused on synthetic problem instances with different condition numbers $\kappa$ of the matrix $X^T X$ and where $X$ is dense. In this case the cost for computation (i) is $O(p)$.  And by taking advantage of the coordinate update structure, we can implement (ii) in $O(p)$ operations by pre-computing and storing $X^T X$ in memory, see \cite{nesterov2012efficiency} and \cite{lin2015accelerated} for further details. The cost of (iii) is simply $O(p )$.

The data $(X,y)$ for the linear regression problems is generated as follows.  For a given number of samples $n$ and problem dimension $p$ (in the experiments we used $n=200$ and $p=100$), we generate a standard Gaussian random matrix $\bar X\in \mathbb{R}^{p\times n}$ with each entry drawn $\sim N(0,1)$. In order to generate the design matrix $X$ with fixed condition number $\kappa$, we first decompose of $\bar X$ as $\bar X=U^T \bar D V$. Then we rescale the diagonal matrix $\bar D$ of singular values linearly to $D$ such that the smallest singular value of $D$ is $\frac{1}{\sqrt{\kappa}}$ and the largest singular value of $D$ is $1$. We then compute the final design  matrix ${X}=U^T{D} V$ and therefore the condition number of  ${X}^T {X}$ becomes $\kappa$.  We generate the response vector $y$ using the linear model $y \sim N(X\beta^*, \sigma^2)$, with true model $\beta^*$ chosen randomly by a Gaussian distribution as well.  For the cases with finite $\kappa$, we are able to compute the strong convexity parameter $\mu$ exactly because the objective function is quadratic, and we use that $\mu$ to implement our Algorithm Framework for strongly convex problems (Algorithm Framework 2). When $\kappa = \infty$, we instead use the smallest positive eigenvalue of $X^TX$ to compute $\mu$.

For the logistic regression experiments, the cost of (ii) at each iteration of AGCD and ASCD can be much larger than $O(p)$ because there is no easy way to update the full gradient $\nabla f(\cdot)$. For these problems we have 
\begin{equation}
\nabla f(\beta) = -\frac{1}{n} X^T w(\beta) \end{equation}
where $X$ is the sample matrix with $x_i$ composing the $i$-th row, and $w(\beta)_i:=\frac{1}{1+\exp(y_i\beta^T x_i)}$. Notice that calculating $w(\beta)$ can be done using a rank $1$-update with cost $O(n)$. But calculating the matrix-vector product $X^T w(\beta)$ will cost $O(np)$, which dominates the cost of (i) and/or (iii). However, in the case when $X$ is a sparse matrix with density $\rho$, the cost can be decreased to $O(\rho np)$. 

\subsection{Comparing the Algorithms using Running Time and the Number of Iterations}

\begin{figure}
  \begin{minipage}[b]{0.4\textwidth}
   \begin{flushright}
    \includegraphics[width=50mm,height=50mm]{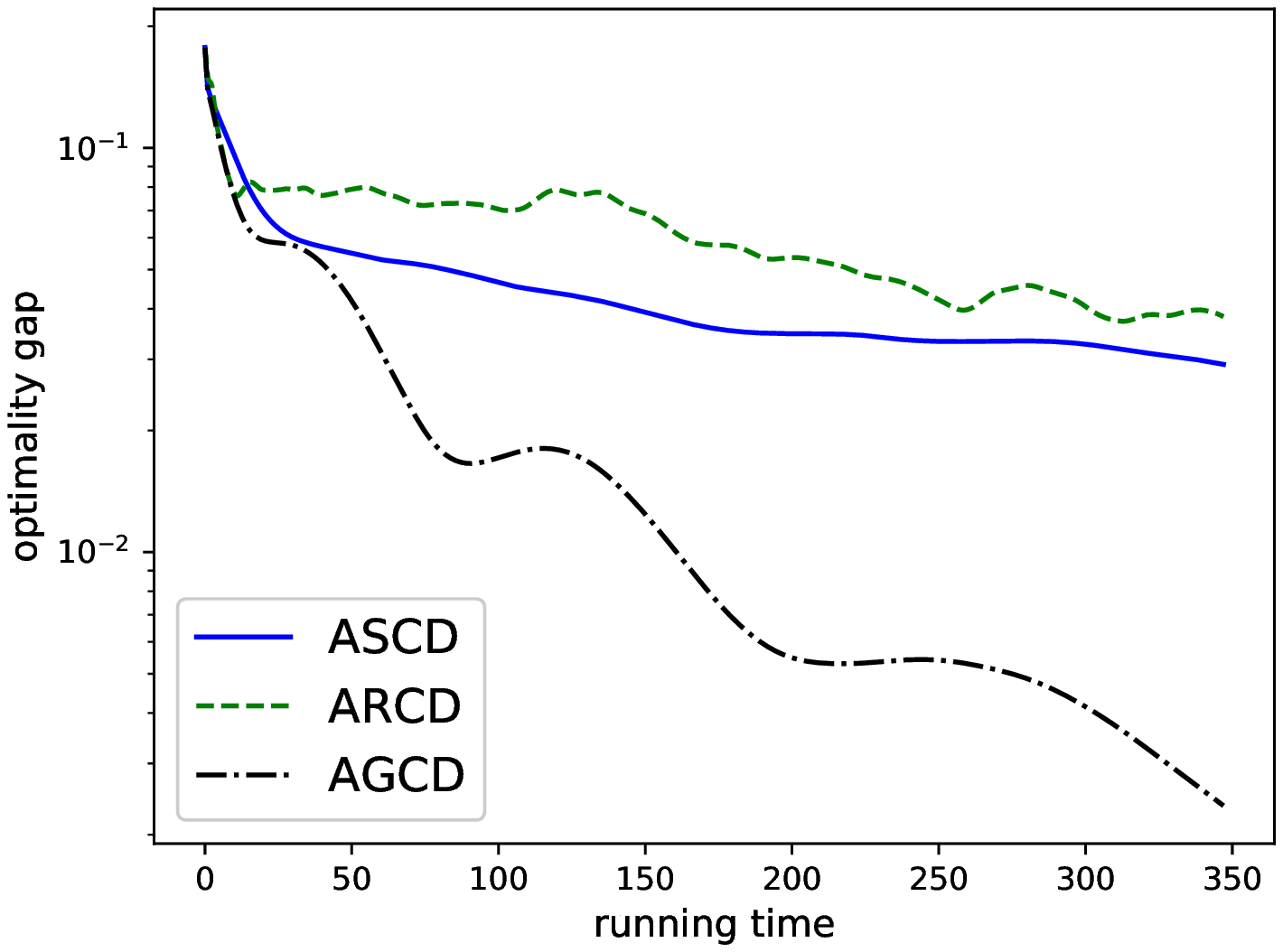}
   \end{flushright}
  \end{minipage}
  \hfill
  \begin{minipage}[b]{0.4\textwidth}
  \begin{flushleft}
    \includegraphics[width=50mm,height=50mm]{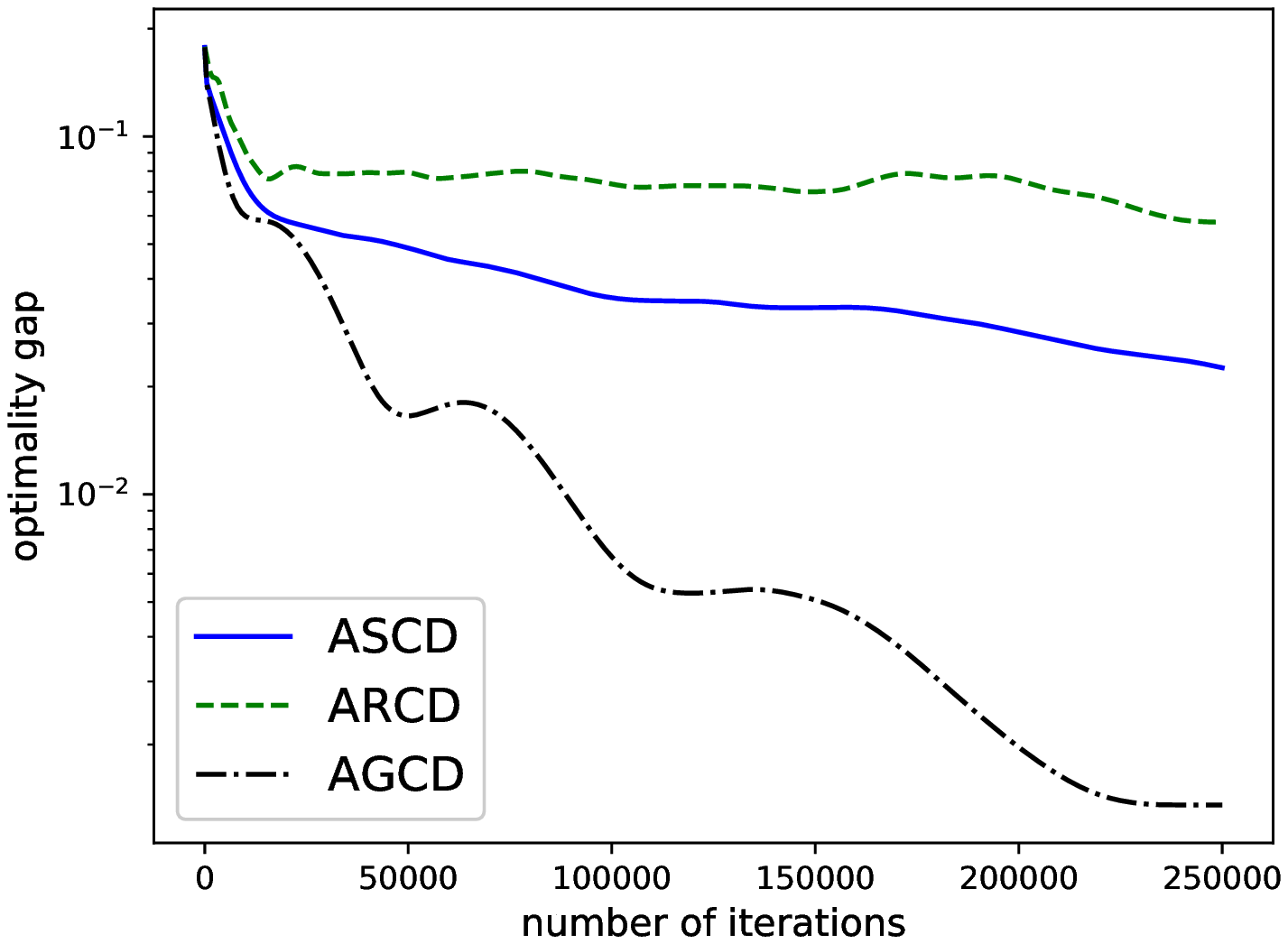}
  \end{flushleft}
  \end{minipage}
      \caption{Plots showing the optimality gap versus run-time (in seconds) on the left and versus the number of iterations on the right, for the logistic regression instance madelon with $\bar \mu=10^{-7}$, solved by ASCD, ARCD and AGCD.}
      \label{iterations}
\end{figure}

Figure \ref{iterations} shows the optimality gap versus running time (seconds) in the left plot and and versus the number of iterations  in the right plot, logistic regression problem using the dataset madelon in LIBSVM \cite{chang2011libsvm}, with $\bar\mu=10^{-7}$. Here we see that AGCD and ASCD are vastly superior to ARCD in term of the number of iterations, but not nearly as much in terms of running time, because one iteration of AGCD or ASCD can be more expensive than an iteration of ARCD.

\subsection{Comparing Accelerated Method with Non-Accelerated Method}
\begin{figure}
\center
  
    \includegraphics[width=50mm,height=50mm]{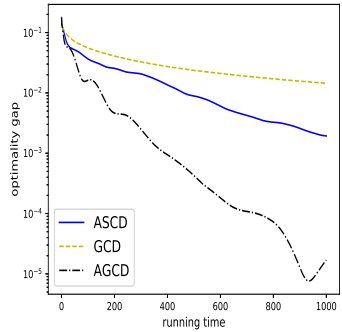}
      \caption{Plots showing the optimality gap versus run-time (in seconds) for the logistic regression instance madelon with $\bar \mu=10^{-6}$, solved by ASCD, GCD and AGCD.}
      \label{GCD}
\end{figure}

Figure \ref{GCD} shows the optimality gap versus running time (seconds) with GCD, ASCD and AGCD for logistic regression problem using the dataset madelon in LIBSVM \cite{chang2011libsvm}, with $\bar\mu=10^{-6}$. Here we see that ASCD and AGCD are superior to non-accelerated GCD.

\subsection{Numerical Results for Logistic Regression with Other Datasets}\label{sec:more_dataset}

\begin{figure}
\centering
\begin{tabular}{c|M{30mm}M{30mm}M{30mm}M{30mm}}
\toprule
Dataset & $\bar\mu=10^{-3}$ & $\bar\mu=10^{-5}$ & $\bar\mu=10^{-7}$ & $\bar\mu=0$ \\
\midrule
heart & \includegraphics[width=35mm,height=35mm]{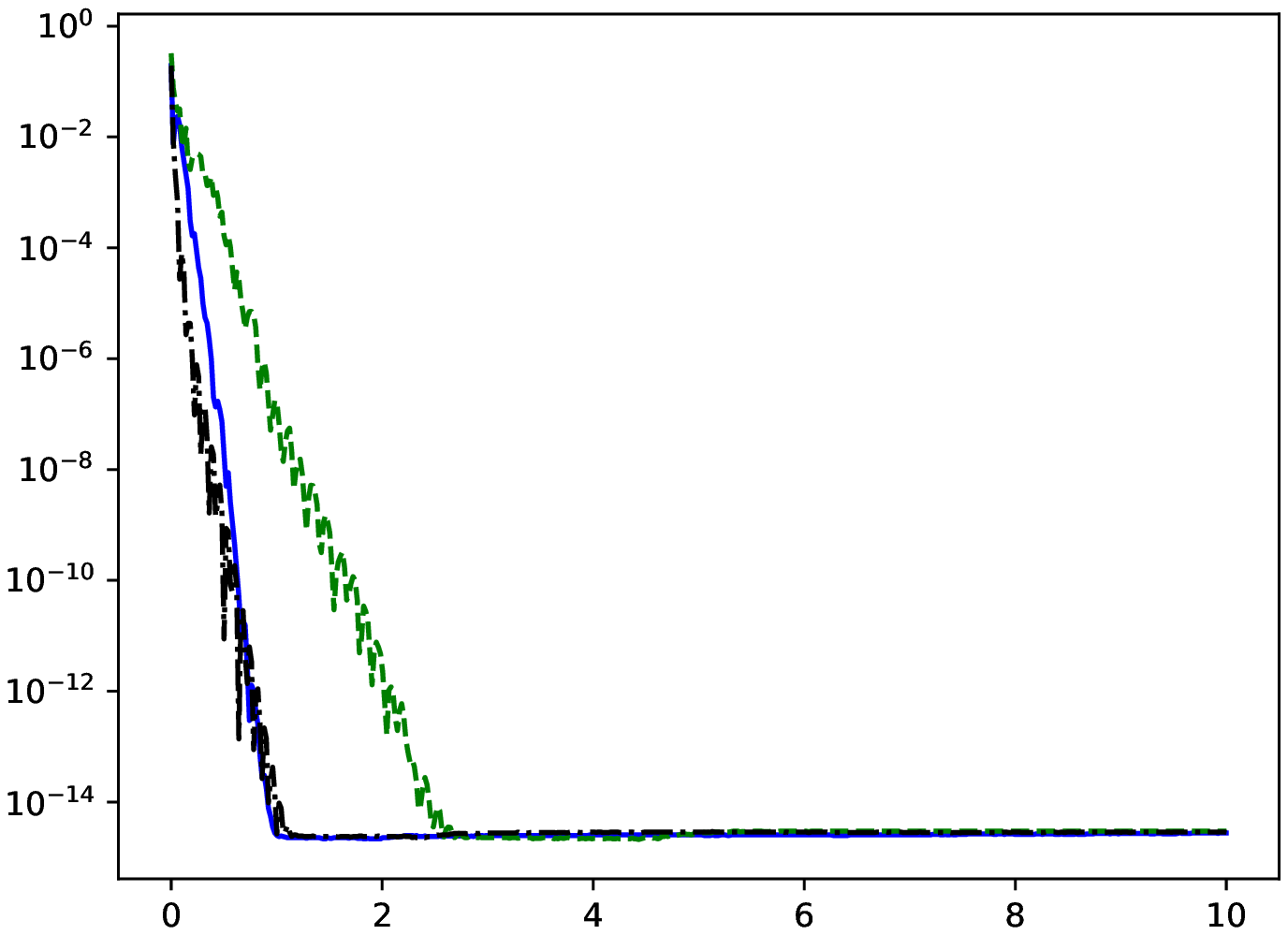} & \includegraphics[width=35mm,height=35mm]{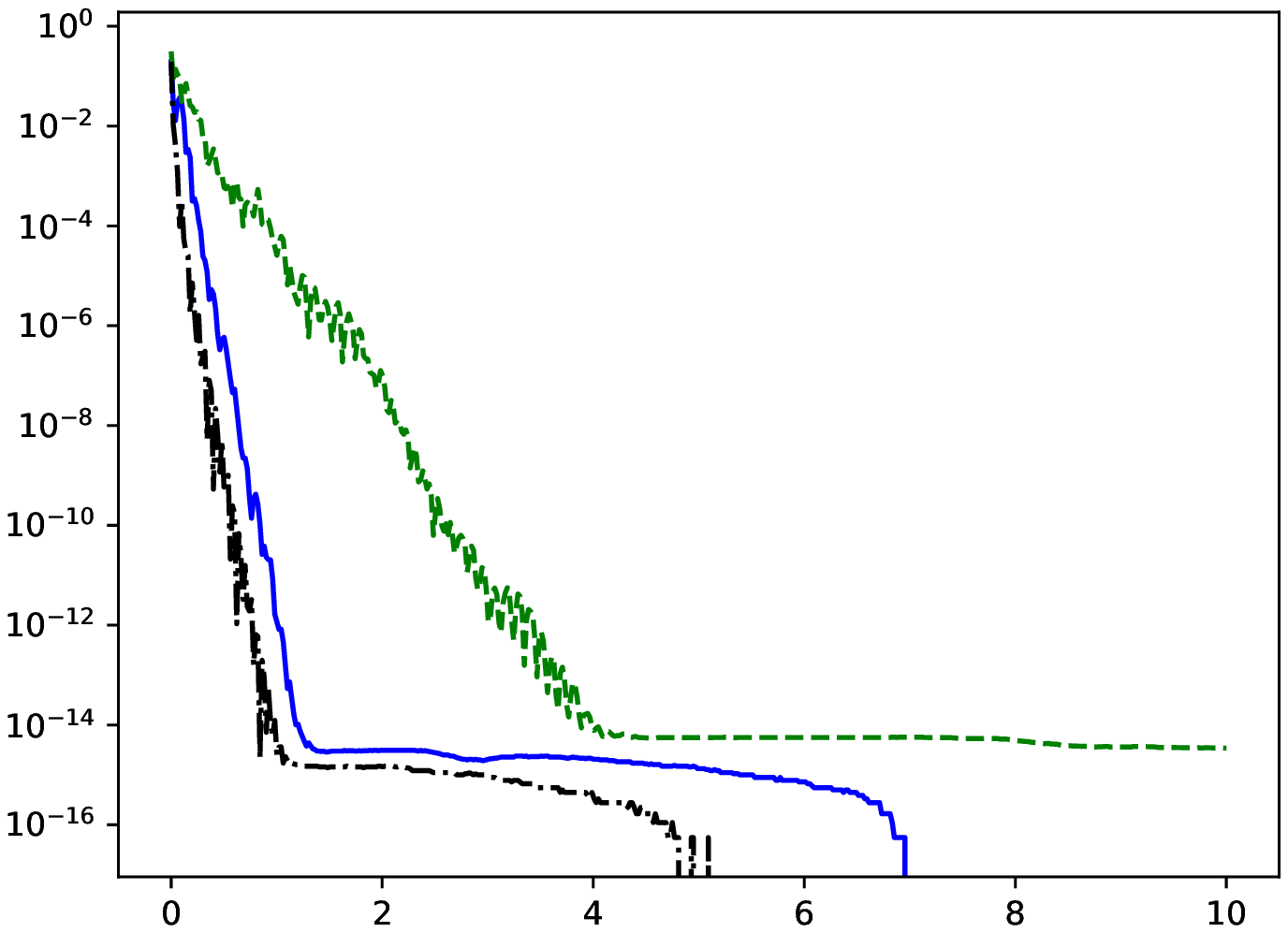} & \includegraphics[width=35mm,height=35mm]{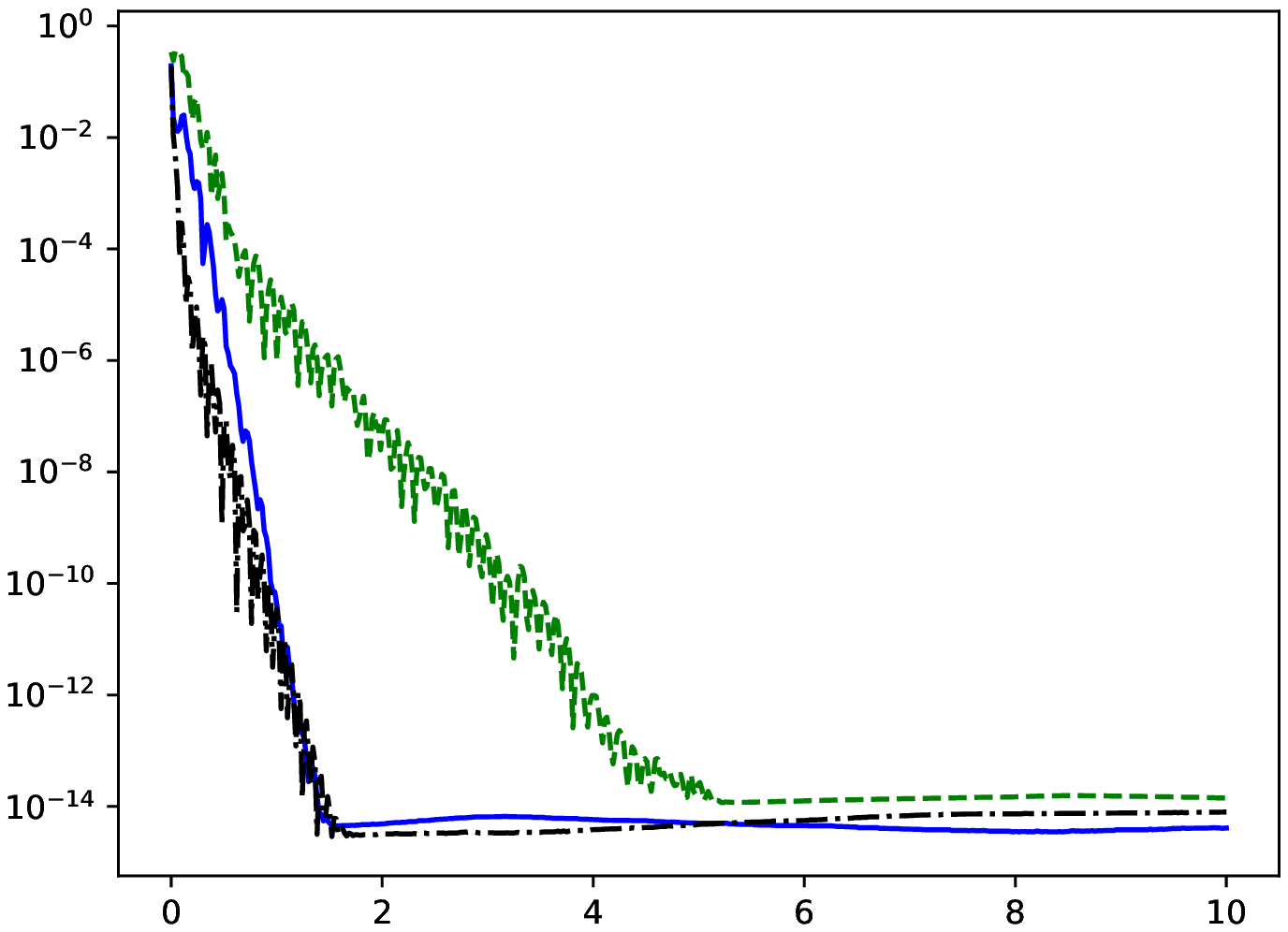} & \includegraphics[width=35mm,height=35mm]{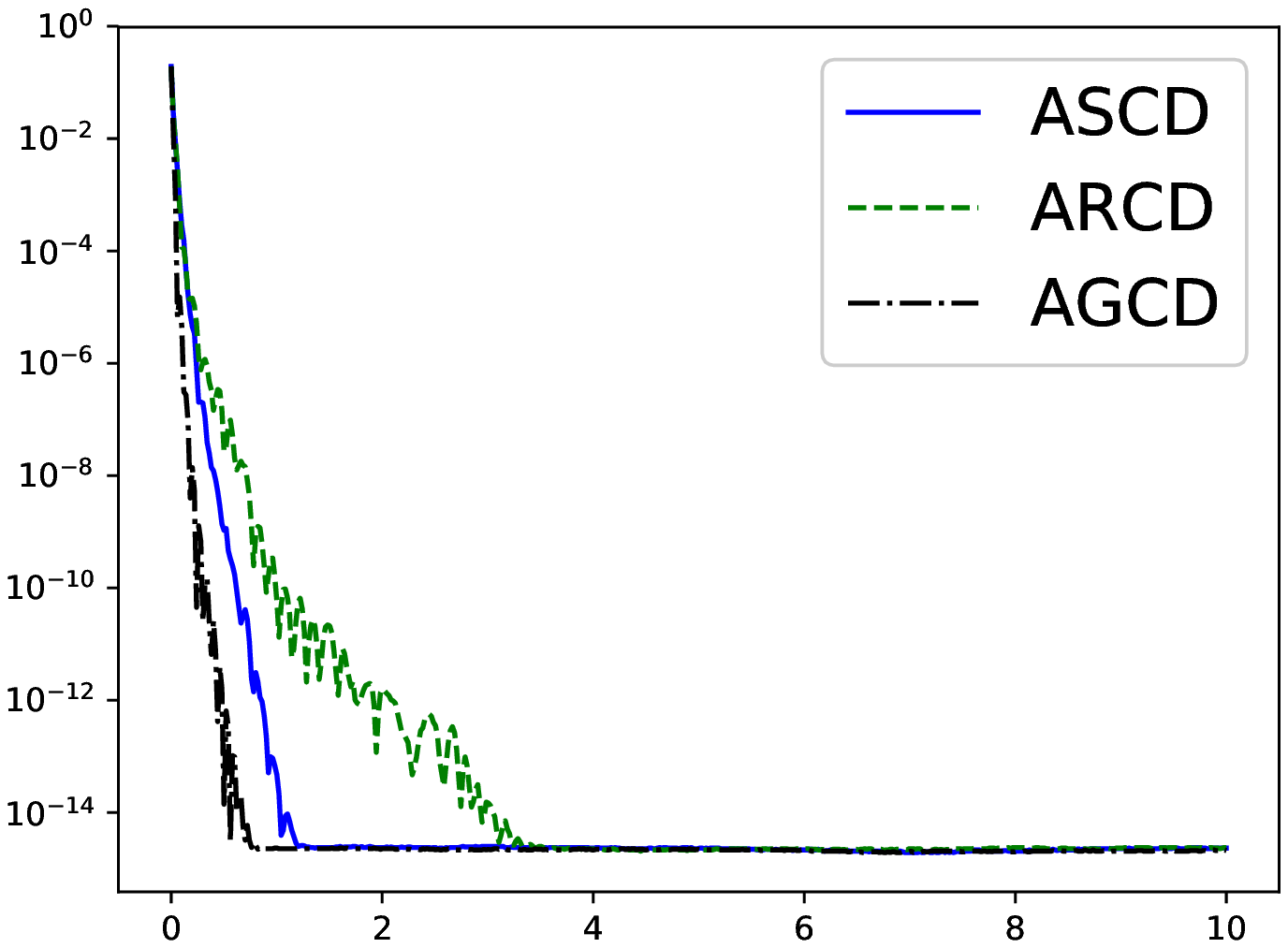} \\
madelon & \includegraphics[width=35mm,height=35mm]{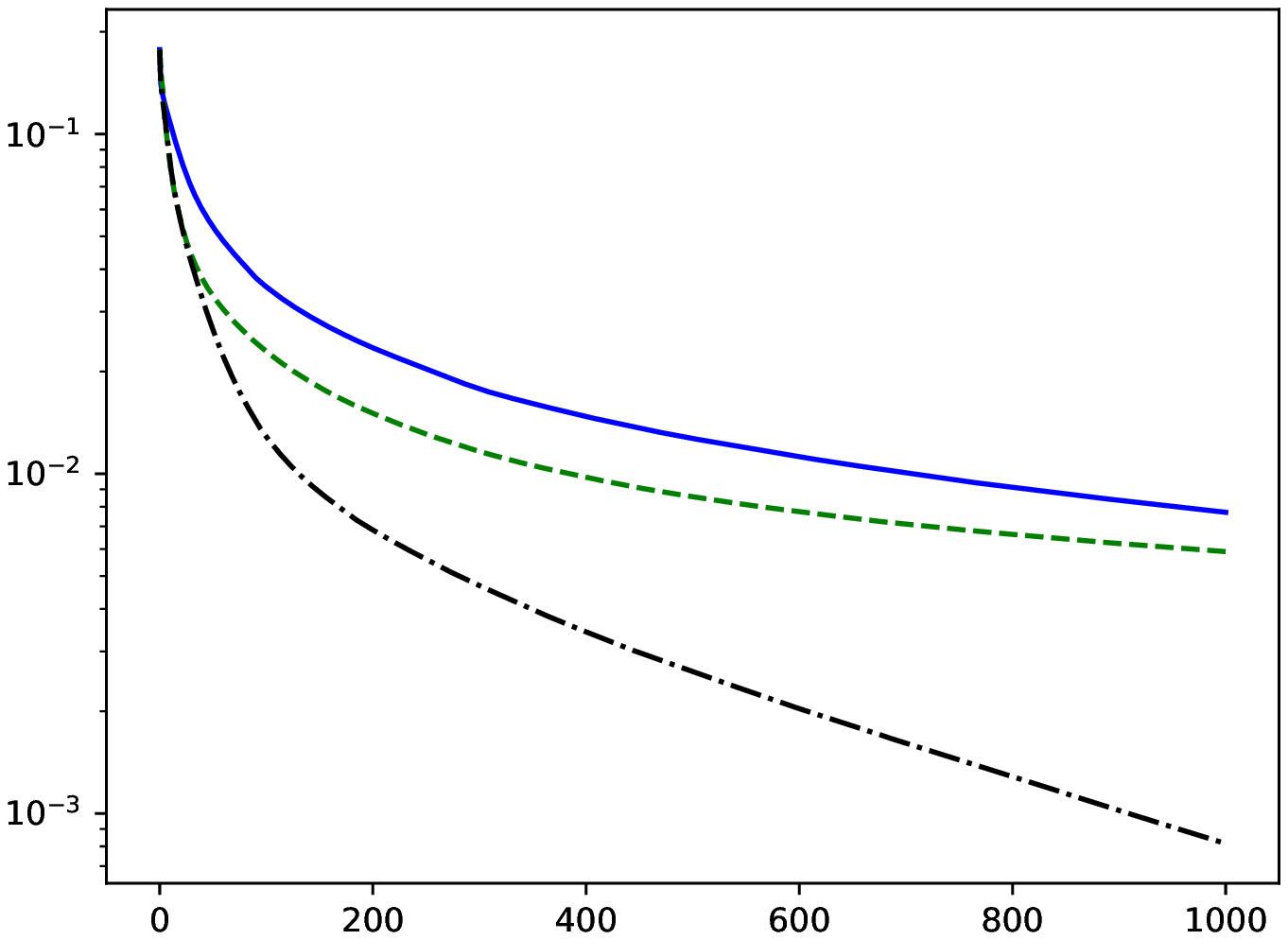} & \includegraphics[width=35mm,height=35mm]{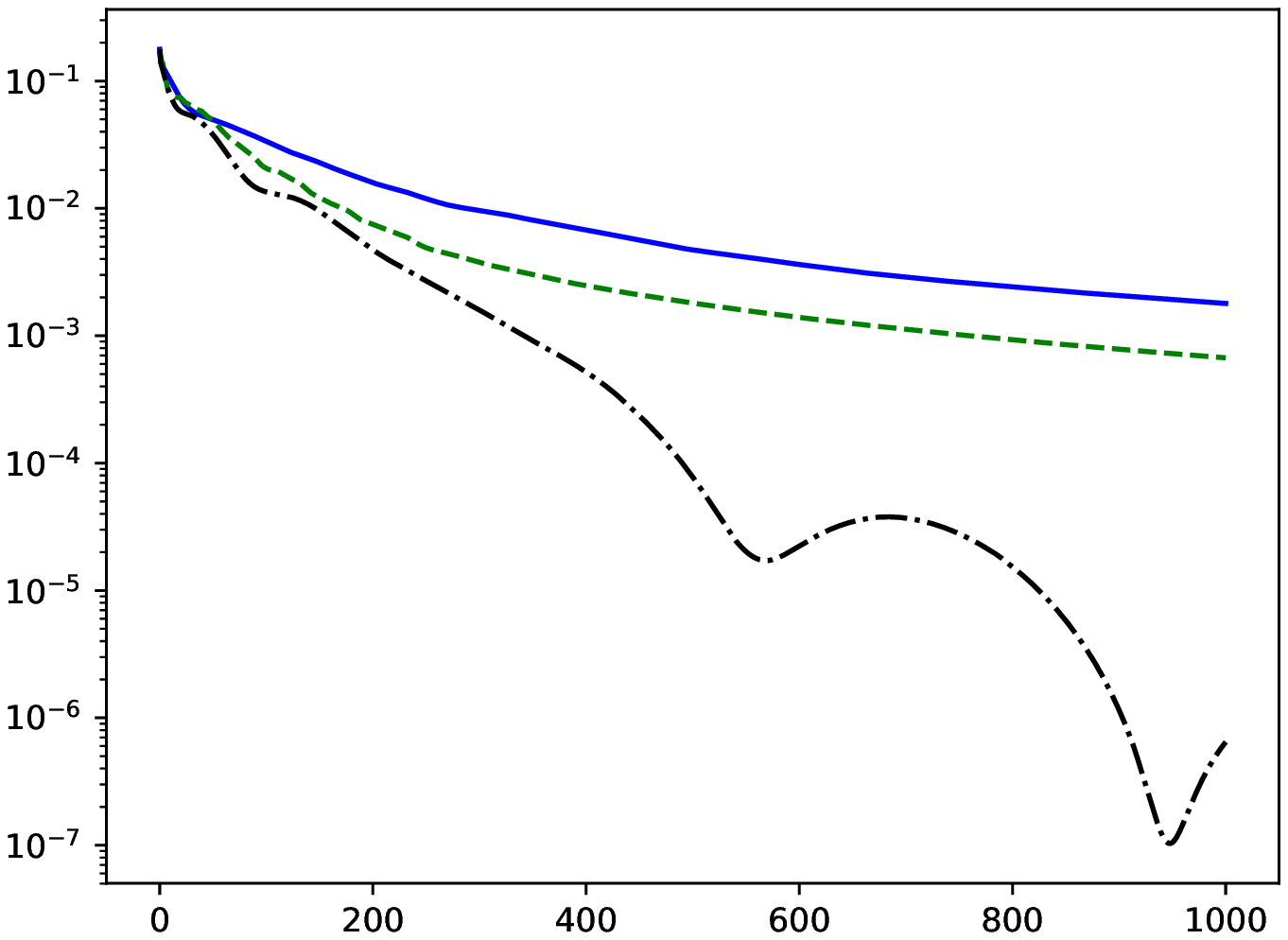} & \includegraphics[width=35mm,height=35mm]{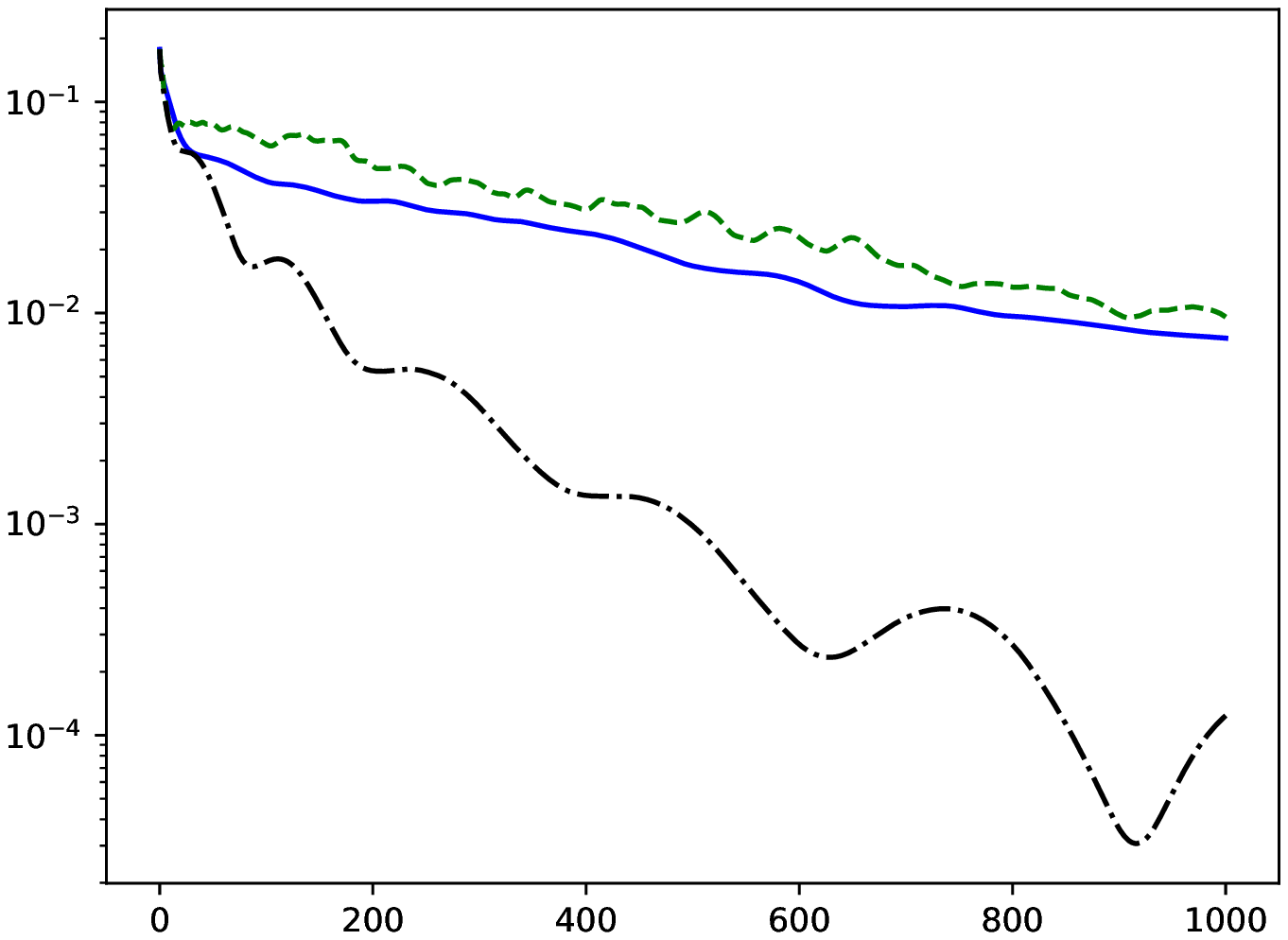} & \includegraphics[width=35mm,height=35mm]{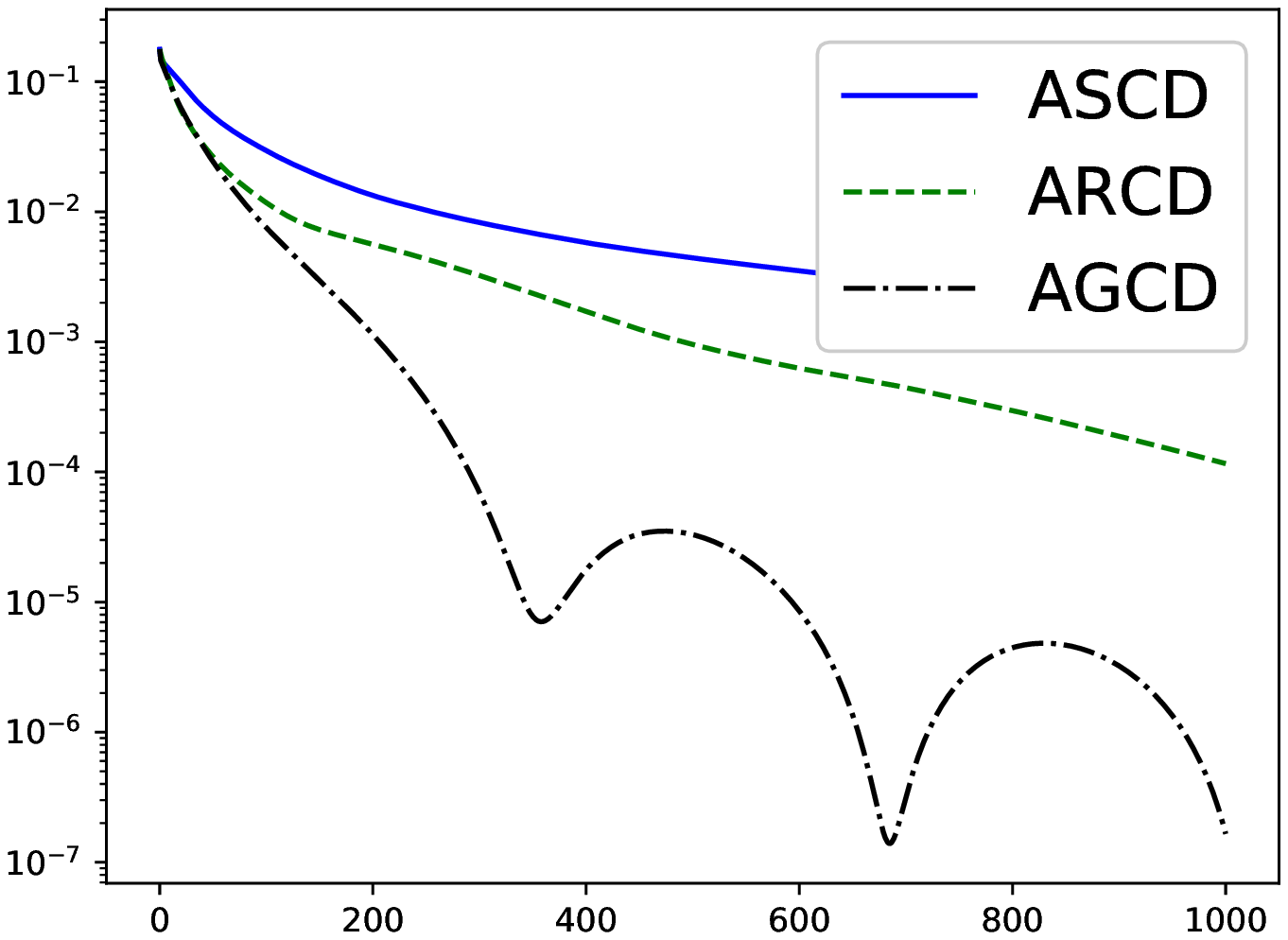} \\
rcv1 & \includegraphics[width=35mm,height=35mm]{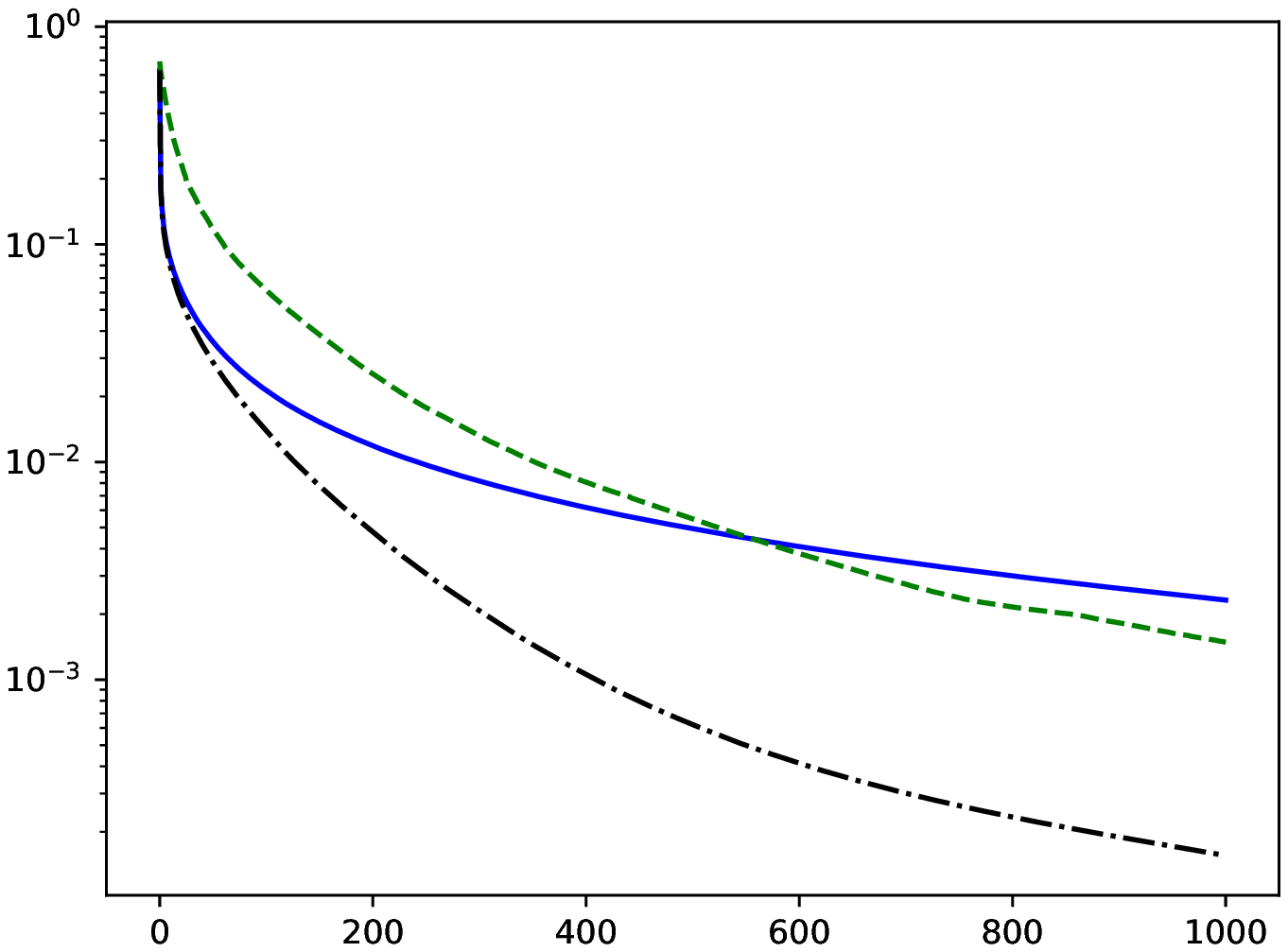} & \includegraphics[width=35mm,height=35mm]{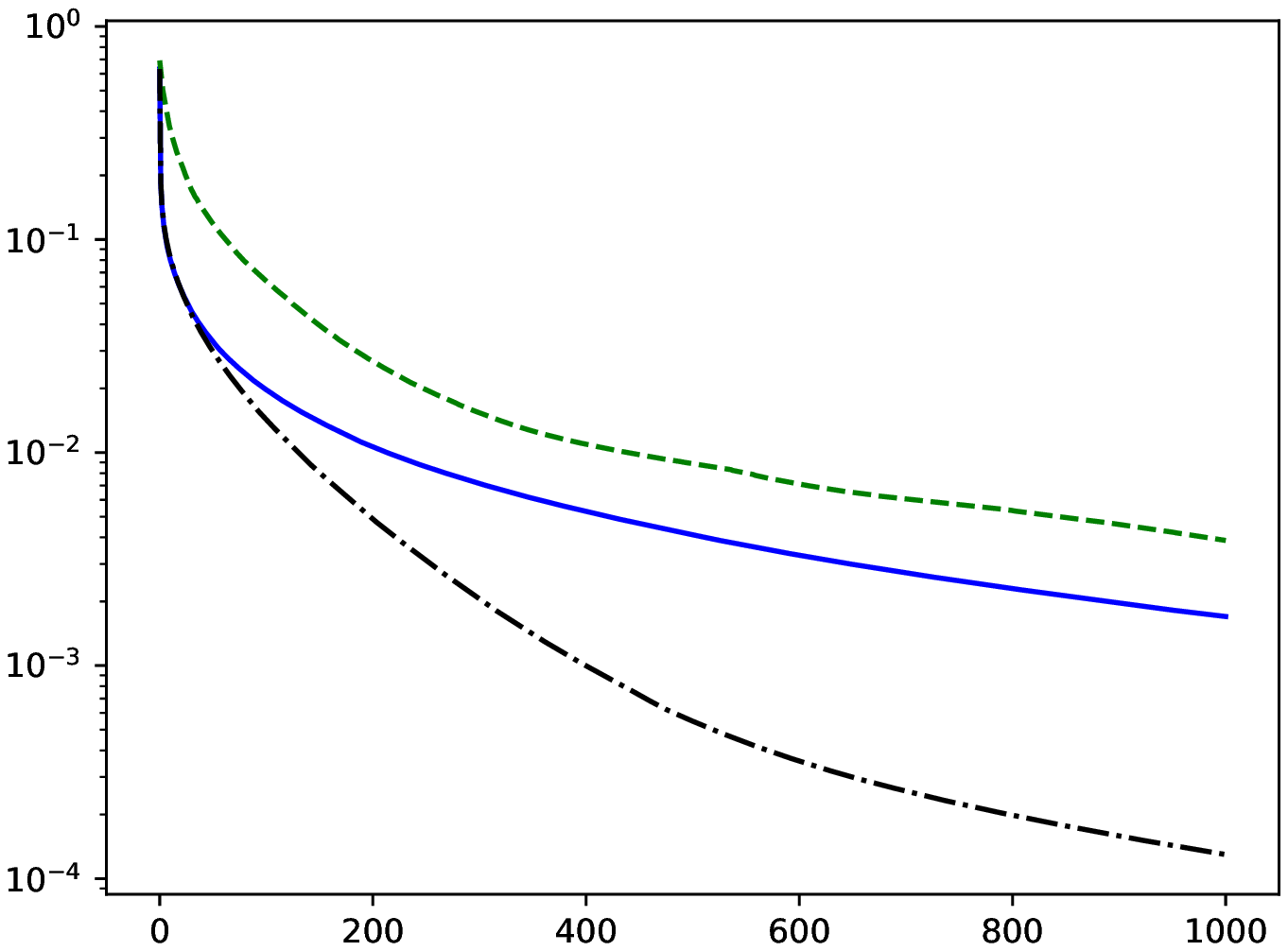} & \includegraphics[width=35mm,height=35mm]{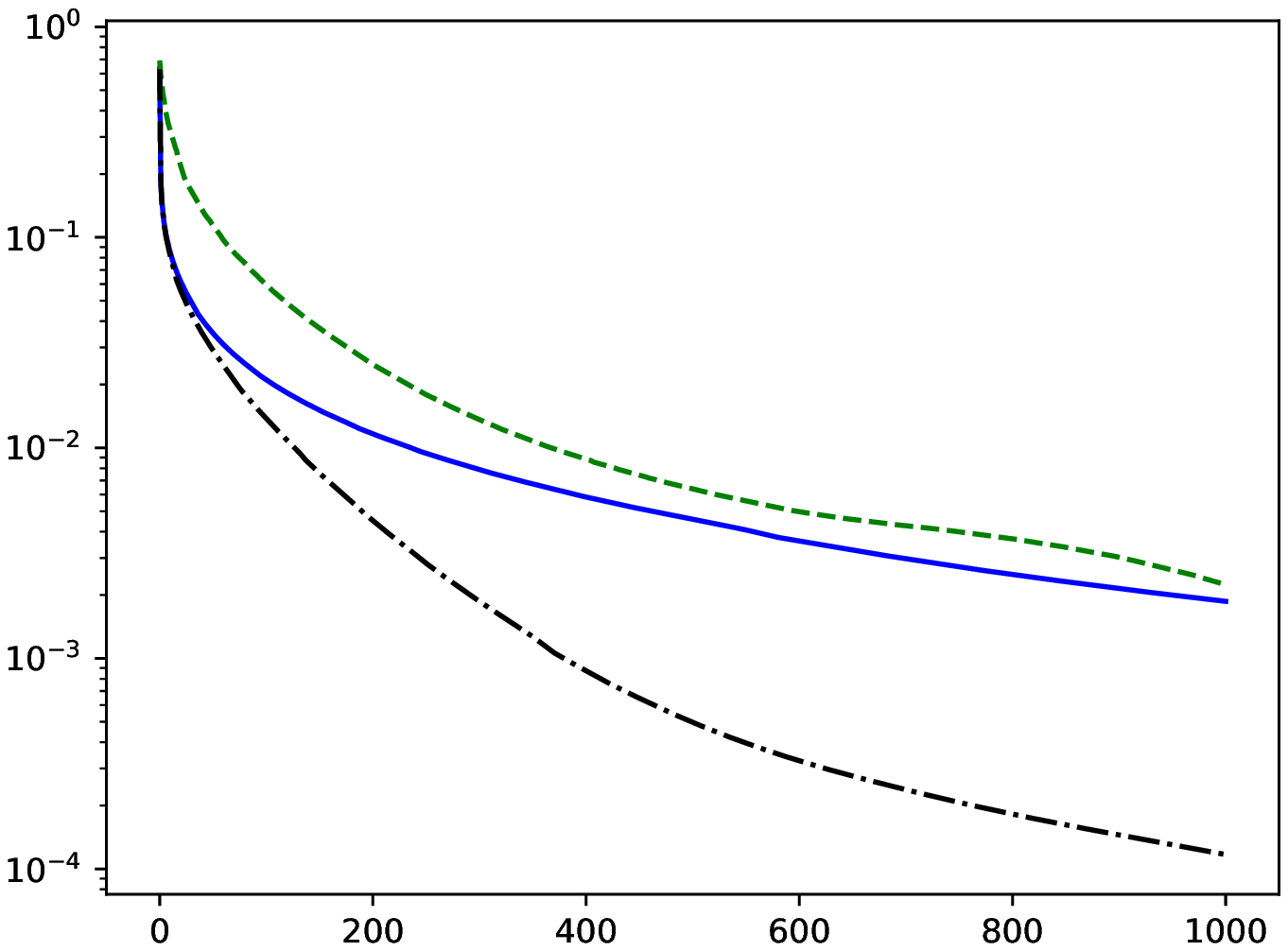} & \includegraphics[width=35mm,height=35mm]{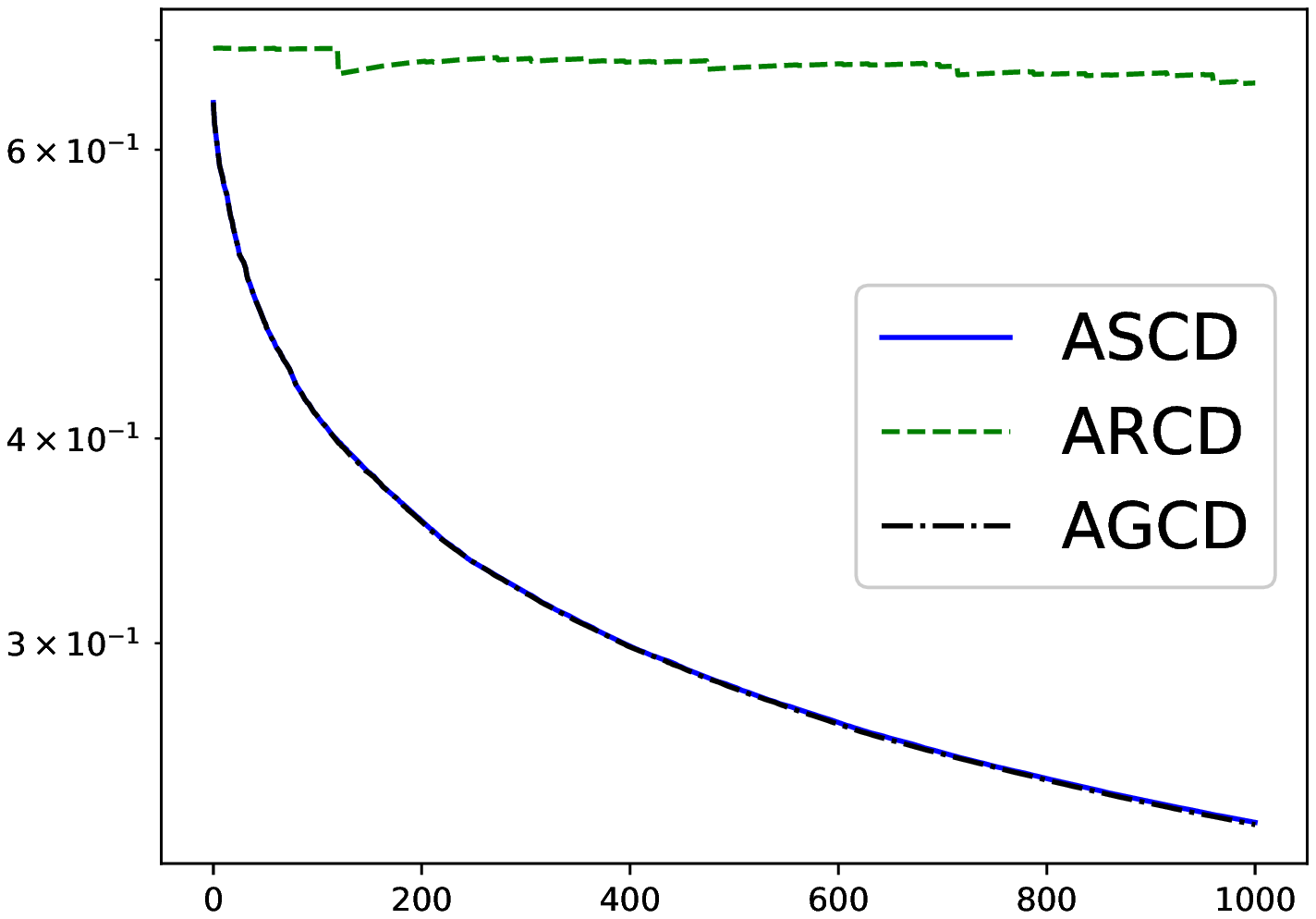}  \\
url & \includegraphics[width=35mm,height=35mm]{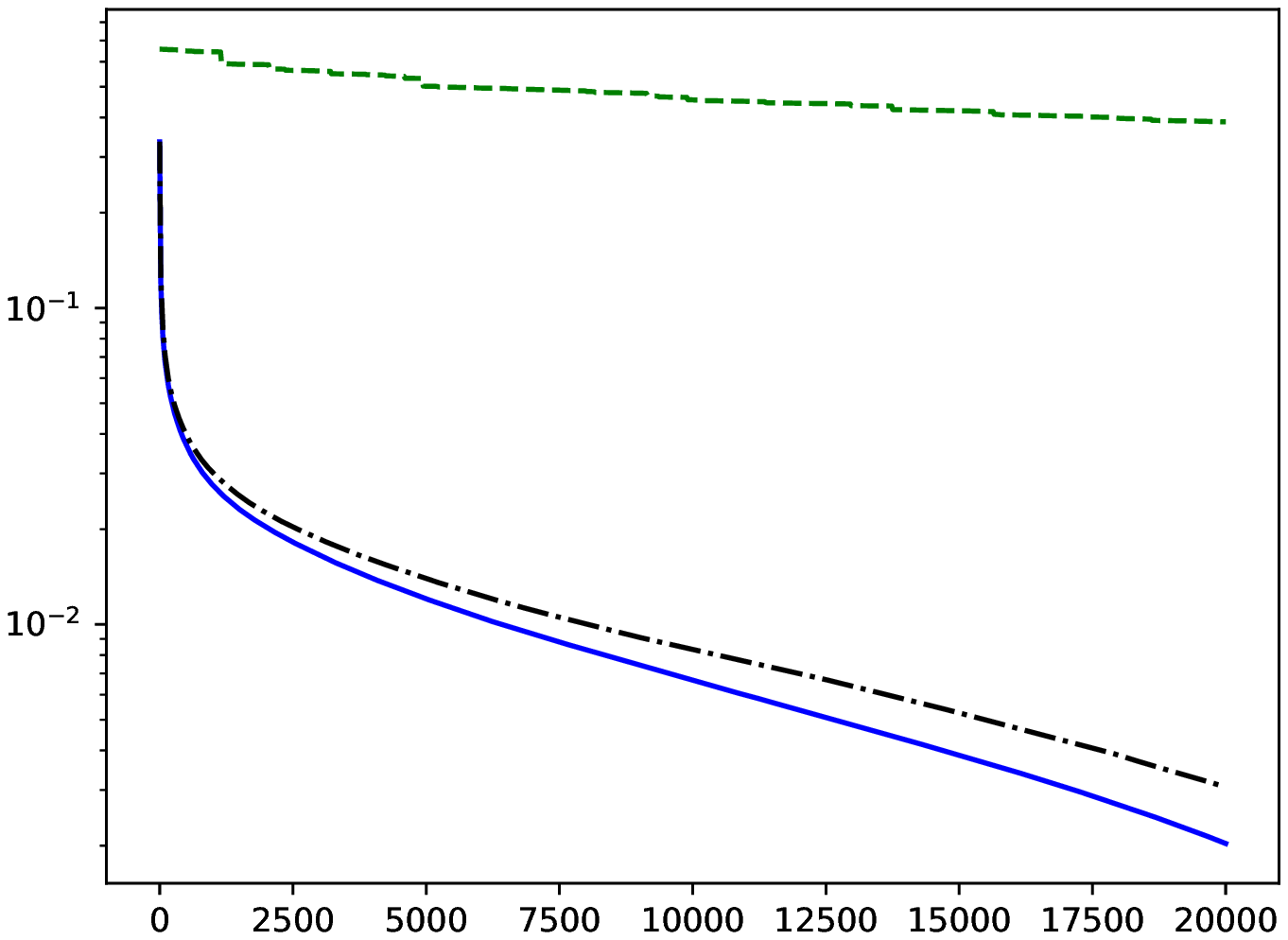} & \includegraphics[width=35mm,height=35mm]{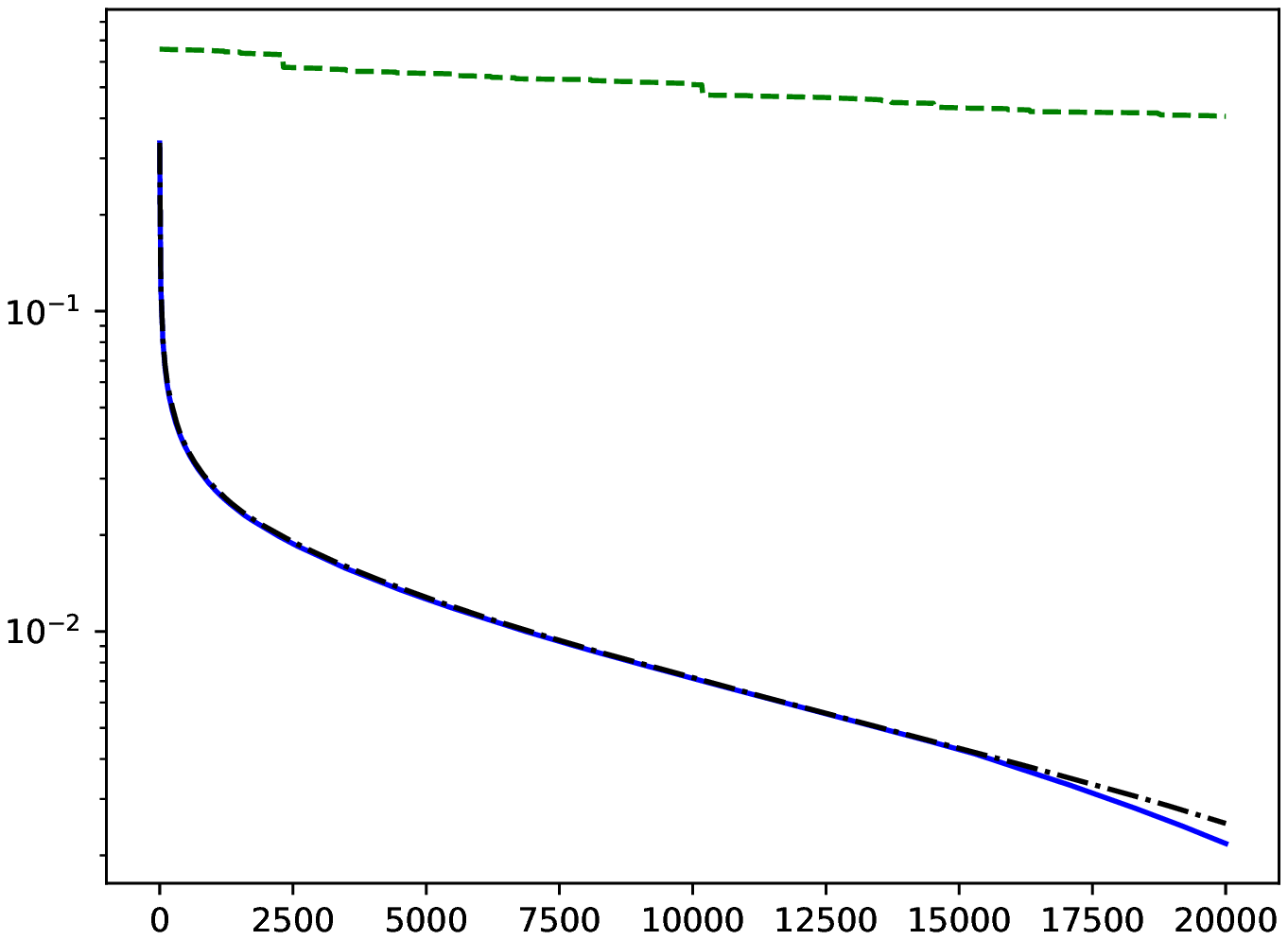} & \includegraphics[width=35mm,height=35mm]{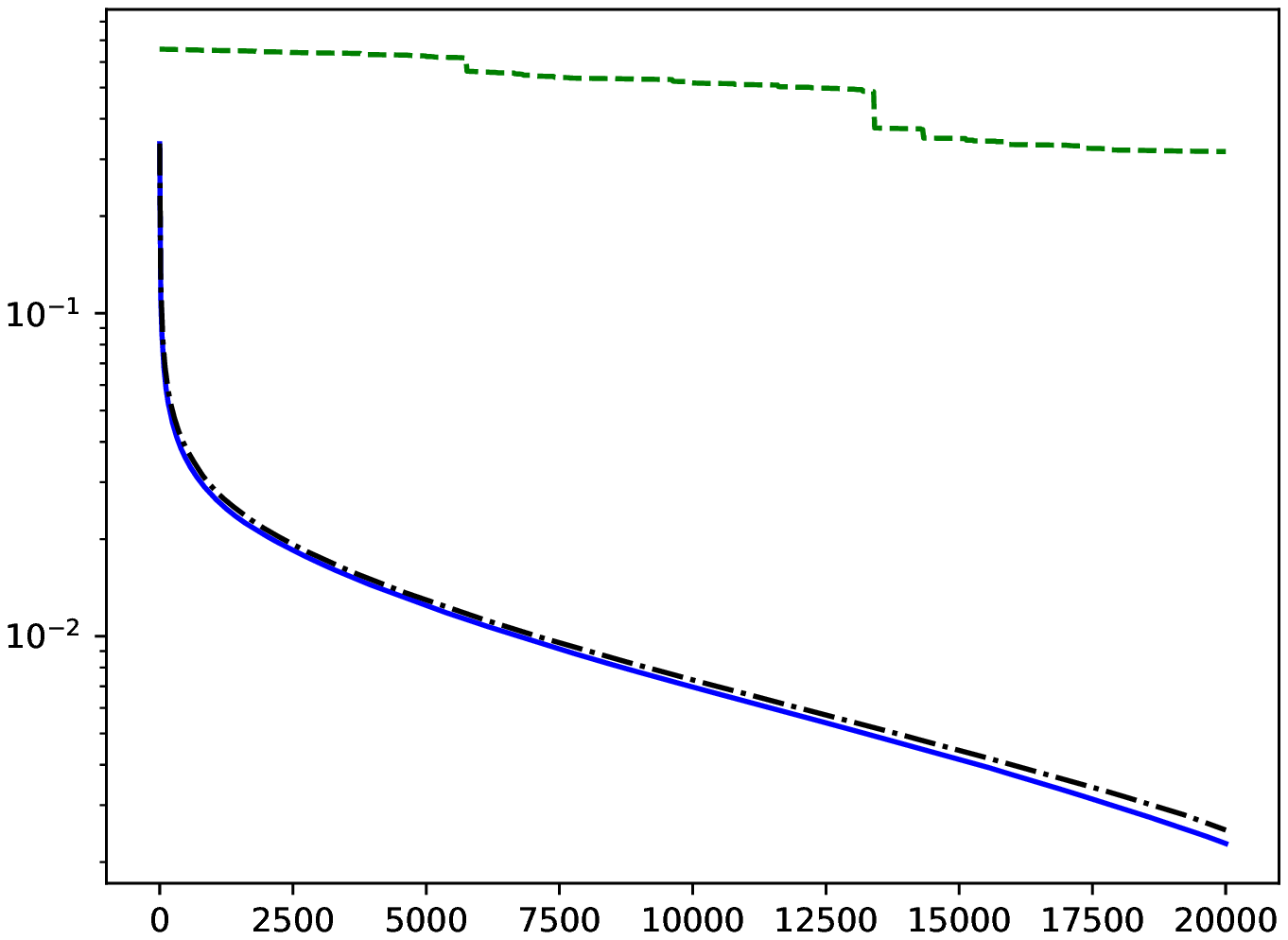} & \includegraphics[width=35mm,height=35mm]{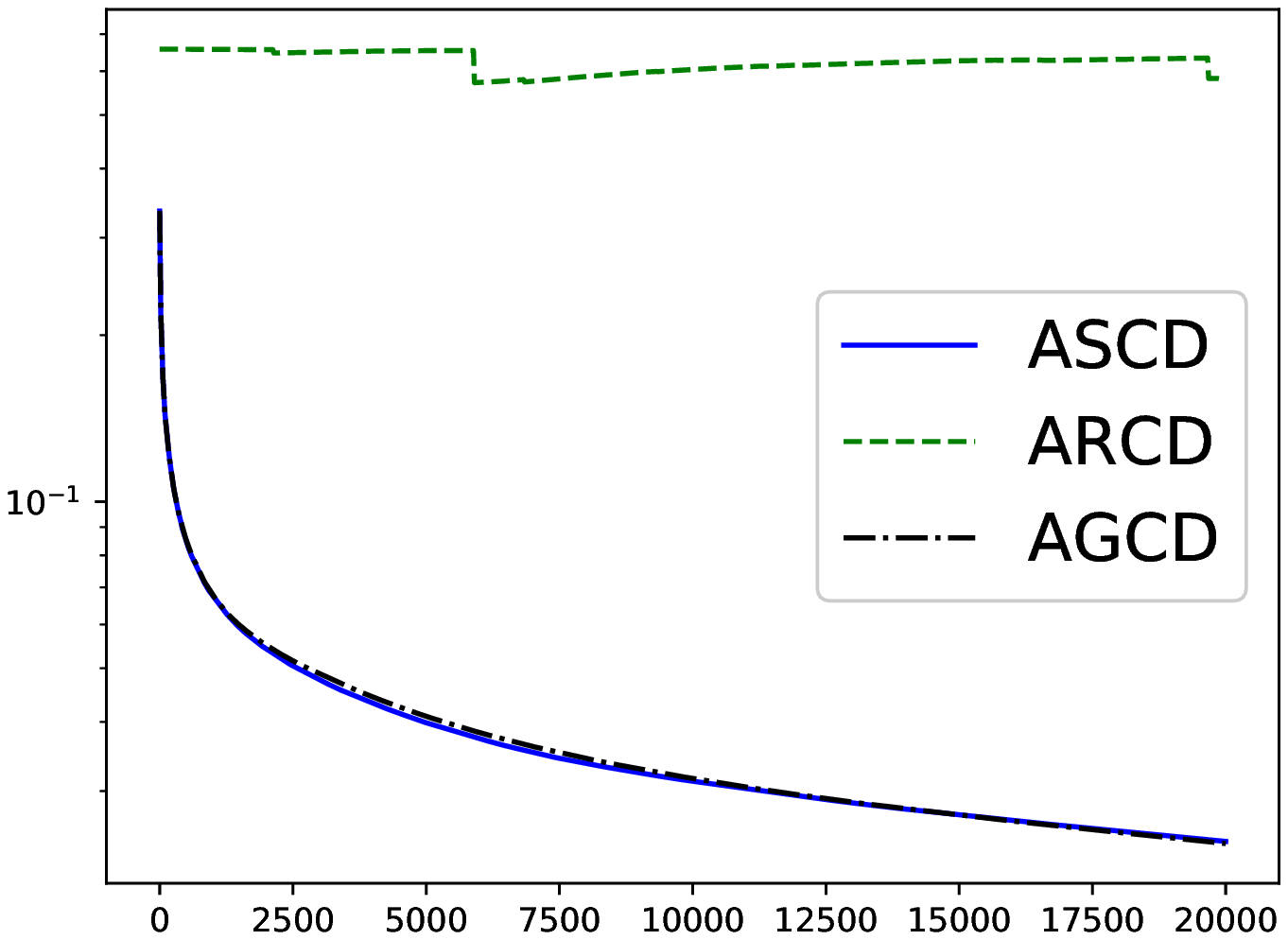} \\
\bottomrule
\end{tabular}
\caption{Plots showing the optimality gap versus run-time (in seconds) for some other logistic regression instances in LIBSVM, solved by ASCD, ARCD and AGCD.}\label{Fig:LogisticRmore}
\end{figure}

We present numerical results for logistic regression problems for several other datasets in LIBSVM solved by ASCD, ARCD and AGCD in Figure \ref{Fig:LogisticRmore}.  Here we see that AGCD always has superior performance as compared to ASCD and ARCD, and ASCD outperforms ARCD in most of the cases.

\bibliographystyle{amsplain}
\bibliography{LF-papers}

\providecommand{\bysame}{\leavevmode\hbox to3em{\hrulefill}\thinspace}
\providecommand{\MR}{\relax\ifhmode\unskip\space\fi MR }
\providecommand{\MRhref}[2]{%
  \href{http://www.ams.org/mathscinet-getitem?mr=#1}{#2}
}
\providecommand{\href}[2]{#2}
\begin{thebibliography}{10}

\bibitem{allen2014linear}
Zeyuan Allen-Zhu and Lorenzo Orecchia, \emph{Linear coupling: An ultimate
  unification of gradient and mirror descent}, arXiv preprint arXiv:1407.1537
  (2014).

\bibitem{allen2016even}
Zeyuan Allen-Zhu, Zheng Qu, Peter Richtarik, and Yang Yuan, \emph{Even faster
  accelerated coordinate descent using non-uniform sampling}, International
  Conference on Machine Learning, 2016.

\bibitem{beck2013convergence}
Amir Beck and Luba Tetruashvili, \emph{On the convergence of block coordinate
  descent type methods}, SIAM journal on Optimization \textbf{23} (2013),
  no.~4, 2037--2060.

\bibitem{bertsekas1989parallel}
Dimitri Bertsekas and John Tsitsiklis, \emph{Parallel and distributed
  computation: numerical methods}, vol.~23, Prentice hall Englewood Cliffs, NJ,
  1989.

\bibitem{bubeck2015geometric}
S{\'e}bastien Bubeck, Yin~Tat Lee, and Mohit Singh, \emph{A geometric
  alternative to nesterov's accelerated gradient descent}, arXiv preprint
  arXiv:1506.08187 (2015).

\bibitem{chang2011libsvm}
Chih-Chung Chang and Chih-Jen Lin, \emph{Libsvm: a library for support vector
  machines}, ACM transactions on intelligent systems and technology (TIST)
  \textbf{2} (2011), no.~3, 27.

\bibitem{fercoq2015accelerated}
Olivier Fercoq and Peter Richtarik, \emph{Accelerated, parallel, and proximal
  coordinate descent}, SIAM Journal on Optimization \textbf{25} (2015), no.~4,
  1997--2023.

\bibitem{frostig2015regularizing}
Roy Frostig, Rong Ge, Sham Kakade, and Aaron Sidford, \emph{Un-regularizing:
  approximate proximal point and faster stochastic algorithms for empirical
  risk minimization}, International Conference on Machine Learning, 2015.

\bibitem{gurbuzbalaban2017cyclic}
Mert Gurbuzbalaban, Asuman Ozdaglar, Pablo~A Parrilo, and Nuri Vanli,
  \emph{When cyclic coordinate descent outperforms randomized coordinate
  descent}, Advances in Neural Information Processing Systems, 2017,
  pp.~7002--7010.

\bibitem{hu2017dissipativity}
Bin Hu and Laurent Lessard, \emph{Dissipativity theory for {N}esterov's
  accelerated method}, arXiv preprint arXiv:1706.04381 (2017).

\bibitem{joachims1998making}
Thorsten Joachims, \emph{Advances in kernel methods}, MIT Press, Cambridge, MA,
  USA, 1999, pp.~169--184.

\bibitem{lee2013efficient}
Yin~Tat Lee and Aaron Sidford, \emph{Efficient accelerated coordinate descent
  methods and faster algorithms for solving linear systems}, Proceedings of the
  2013 IEEE 54th Annual Symposium on Foundations of Computer Science
  (Washington, DC, USA), FOCS '13, IEEE Computer Society, 2013, pp.~147--156.

\bibitem{lin2015universal}
Hongzhou Lin, Julien Mairal, and Zaid Harchaoui, \emph{A universal catalyst for
  first-order optimization}, Advances in Neural Information Processing Systems,
  2015.

\bibitem{lin2015accelerated}
Qihang Lin, Zhaosong Lu, and Lin Xiao, \emph{An accelerated randomized proximal
  coordinate gradient method and its application to regularized empirical risk
  minimization}, SIAM Journal on Optimization \textbf{25} (2015), no.~4,
  2244--2273.

\bibitem{locatello2018revisiting}
Francesco Locatello, Anant Raj, Sai~Praneeth Reddy, Gunnar R{\"a}tsch, Bernhard
  Sch{\"o}lkopf, Sebastian~U Stich, and Martin Jaggi, \emph{On matching pursuit
  and coordinate descent}, ICML 2018 - Proceedings of the 35th International
  Conference on Machine Learning (2018).

\bibitem{lu2015complexity}
Zhaosong Lu and Lin Xiao, \emph{On the complexity analysis of randomized
  block-coordinate descent methods}, Mathematical Programming \textbf{152}
  (2015), no.~1-2, 615--642.

\bibitem{luo1992convergence}
Zhi-Quan Luo and Paul Tseng, \emph{On the convergence of the coordinate descent
  method for convex differentiable minimization}, Journal of Optimization
  Theory and Applications \textbf{72} (1992), no.~1, 7--35.

\bibitem{luo1993error}
\bysame, \emph{Error bounds and convergence analysis of feasible descent
  methods: a general approach}, Annals of Operations Research \textbf{46}
  (1993), no.~1, 157--178.

\bibitem{mazumder2012graphical}
Rahul Mazumder and Trevor Hastie, \emph{The graphical lasso: new insights and
  alternatives}, Electronic Journal of Statistics \textbf{6} (2012), 2125.

\bibitem{NemirovskyYudin83}
A.~S. Nemirovsky and D.~B. Yudin, \emph{Problem complexity and method
  efficiency in optimization}, Wiley, New York, 1983.

\bibitem{nesterov2012efficiency}
Yu~Nesterov, \emph{Efficiency of coordinate descent methods on huge-scale
  optimization problems}, SIAM Journal on Optimization \textbf{22} (2012),
  no.~2, 341--362.

\bibitem{nesterov1983method}
Yurii Nesterov, \emph{A method of solving a convex programming problem with
  convergence rate ${O}(1/k^2)$}, Soviet Mathematics Doklady, vol.~27, 1983,
  pp.~372--376.

\bibitem{nutini2015coordinate}
Julie Nutini, Mark Schmidt, Issam Laradji, Michael Friedlander, and Hoyt
  Koepke, \emph{Coordinate descent converges faster with the
  {G}auss-{S}outhwell rule than random selection}, International Conference on
  Machine Learning, 2015, pp.~1632--1641.

\bibitem{zeng2008fast}
John~C. Platt, \emph{Advances in kernel methods}, MIT Press, Cambridge, MA,
  USA, 1999, pp.~185--208.

\bibitem{richtarik2014iteration}
Peter Richtarik and Martin Takac, \emph{Iteration complexity of randomized
  block-coordinate descent methods for minimizing a composite function},
  Mathematical Programming \textbf{144} (2014), no.~1-2, 1--38.

\bibitem{song2017accelerated}
Chaobing Song, Shaobo Cui, Yong Jiang, and Shu-Tao Xia, \emph{Accelerated
  stochastic greedy coordinate descent by soft thresholding projection onto
  simplex}, Advances in Neural Information Processing Systems, 2017,
  pp.~4841--4850.

\bibitem{su2016differential}
Weijie Su, Stephen Boyd, and Emmanuel~J Candes, \emph{A differential equation
  for modeling {N}esterov’s accelerated gradient method: theory and
  insights}, Journal of Machine Learning Research \textbf{17} (2016), no.~153,
  1--43.

\bibitem{sun2016worst}
Ruoyu Sun and Yinyu Ye, \emph{Worst-case complexity of cyclic coordinate
  descent: $ {O} (n^{2}) $ gap with randomized version}, arXiv preprint
  arXiv:1604.07130 (2016).

\bibitem{tsengaccelerated}
P.~Tseng, \emph{On accelerated proximal gradient methods for convex-concave
  optimization}, Tech. report, May 21, 2008.

\bibitem{wilson2016lyapunov}
Ashia~C Wilson, Benjamin Recht, and Michael~I Jordan, \emph{A {L}yapunov
  analysis of momentum methods in optimization}, arXiv preprint
  arXiv:1611.02635 (2016).

\bibitem{you2016asynchronous}
Yang You, Xiangru Lian, Ji~Liu, Hsiang-Fu Yu, Inderjit~S Dhillon, James Demmel,
  and Cho-Jui Hsieh, \emph{Asynchronous parallel greedy coordinate descent},
  Advances in Neural Information Processing Systems, 2016, pp.~4682--4690.

\end{thebibliography}

\end{document}